\documentclass[preprint,12pt,sort&compress]{elsarticle} 
\pdfoutput=1
\usepackage{epsfig} 
\usepackage{graphicx}
\usepackage[ansinew]{inputenc}
\usepackage{epstopdf}
\usepackage{amsmath}
\usepackage{psfrag}
\usepackage{amsthm}
\usepackage{amssymb}
\usepackage{float}
\usepackage{multirow}
\usepackage{lineno}
\usepackage{setspace}
\usepackage{flushend}
\usepackage{array}
\usepackage{styles/tabu}
\usepackage{bm}
\usepackage{color}
\newcolumntype{L}[1]{>{\raggedright\let\newline\\\arraybackslash\hspace{0pt}}m{#1}}
\newcolumntype{C}[1]{>{\centering\let\newline\\\arraybackslash\hspace{0pt}}m{#1}}
\newcolumntype{R}[1]{>{\raggedleft\let\newline\\\arraybackslash\hspace{0pt}}m{#1}}
\usepackage[figuresright]{rotating}
\usepackage{subfigure}
\usepackage[title,titletoc,toc]{appendix}
\usepackage{mathrsfs}
\usepackage{amsfonts}
\usepackage{svgcolor}
\usepackage{color}
\usepackage{styles/ieeetrantools}
\usepackage{times}
\usepackage{mathptmx}
\usepackage{booktabs}
\usepackage{amsthm}
\usepackage{styles/numcompress}
\usepackage{url}

\theoremstyle{lemma}

\theoremstyle{definition}
\newtheorem*{remark}{Remark}

\theoremstyle{statement}
\newtheorem*{statement}{Statement}

\definecolor{azul1}{RGB}{0,153,153}
\definecolor{azul2}{RGB}{102,178,255}
\definecolor{azul3}{RGB}{76,0,153}

\definecolor{rosa1}{RGB}{204,0,102}
\definecolor{rosa2}{RGB}{255,178,102}
\definecolor{rosa3}{RGB}{255,102,102}

\journal{ISA Transactions}

\begin{document}

\begin{frontmatter}

\title{Optimal control analysis and Practical NMPC applied to refrigeration systems \footnote{© 2020. This manuscript version is made available under the CC-BY-NC-ND 4.0 license \url{https://creativecommons.org/licenses/by-nc-nd/4.0/}. The link to the formal publication is \url{https://doi.org/10.1016/j.isatra.2020.07.041}}}

\author[ula]{G. Bejarano\corref{mycorrespondingauthor}}
\cortext[mycorrespondingauthor]{Corresponding author}
\ead{gbejarano@uloyola.es}

\author[us]{M. G. Ortega}
\author[ufsc]{J. E. Normey-Rico}
\author[us]{F. R. Rubio}    

\address[ula]{Departamento de Ingeniería, Escuela Técnica Superior de Ingeniería, Universidad Loyola Andalucía. Avda. de las Universidades s/n, 41704 Dos Hermanas, Sevilla (España)}
\address[us]{Departamento de Ingeniería de Sistemas y Automática, Escuela Técnica Superior de Ingeniería, Universidad de Sevilla. Camino de los Descubrimientos s/n, 41092 Sevilla (España)}
\address[ufsc]{Departamento de Automação e Sistemas, Universidade Federal de Santa Catarina, 88040-900, Florianópolis, SC, Brasil}

\begin{abstract}

{
	This work is focused on optimal control of mechanical compression refrigeration systems. A reduced-order state-space model based on the \emph{moving boundary} approach is proposed for the canonical cycle, which eases the controller design. The optimal cycle (that satisfying the cooling demand while maximizing efficiency) is defined by three variables, but only two inputs are available, therefore the controllability of the proposed model is studied. It is shown through optimization simulations how optimal cycles for a range of the cooling demand turn out not to be achieved by keeping the degree of superheating to a minimum. The Practical NMPC and a well-known feedback-plus-feedforward strategy from the literature are compared in simulation, both showing trouble in reaching the optimal cycle, which agrees with the controllability study.
}
\end{abstract}

\begin{keyword}

 Refrigeration system \sep 
 Process control \sep 
 Global optimization \sep 
 Controllability study \sep
 Model predictive control

\end{keyword}

\end{frontmatter}

\begin{table*}
	\scalebox{0.7}[0.7]{
	\tabulinesep=0.5mm
	\begin{tabu} {L{1.8cm} L{8cm} L{1.8cm} L{8cm}}
		\multicolumn{4}{l}{{\Large \it \bf Terminology}}\\
		\multicolumn{2}{l}{\it \bf Latin symbols} & \multicolumn{2}{l}{\it \bf Greek symbols} \\
		$\bm A$ & Linear dynamic matrix & $\alpha$ & Heat transfer coefficient [W m\textsuperscript{-2} K\textsuperscript{-1}] \\
		$A$ & Heat transfer perimeter [m] & $\bar{\gamma}$ & Mean void fraction \\
		$A_v$ & Expansion valve opening [\%] & $\varepsilon$ & Heat transfer effectiveness \\
		$a$ & Compressor constant [W] & $\zeta$ & Normalised heat exchanger length \\
		$\bm B$ & Linear input matrix & $\bm \Lambda$ & Projection matrix \\
		$b$ & Compressor constant & $\rho$ & Density [kg $ \text{m}^{\text{-3}}$] \\
		$\bm{\mathbb{C}}$ & Linear controllability matrix & $\bm \phi$ & Measurable variable set \\
		$C$ & Thermal capacity ratio & $\bm \chi$ & Optimization variable set \\
		$COP$ & Coefficient of Performance & $\bm \psi$ & Error state vector \\
		$CR$ & Compression ratio & & \\
		$c$ & Compressor constant [m\textsuperscript{3}] & & \\
		$c_{eev}$ & Expansion valve coefficient [m\textsuperscript{2}] & \multicolumn{2}{l}{\it \bf Subscripts and superscripts} \\ 
		$c_p$ & Specific heat at constant pressure [J kg\textsuperscript{-1} K\textsuperscript{-1}] & $c$ & condenser \\ 
		$c_v$ & Specific heat at constant volume [J kg\textsuperscript{-1} K\textsuperscript{-1}] & $cntrl$ & control \\
		$d$ & Element defining the reduced error model & $comp$ & compressor \\
		$\bm e$ & Steady-state model iterative error & $cycle$ & mechanical compression cycle \\
		$\bm f$ & Force function & $e$ & evaporator \\ 
		$\bm G$ & Jacobian matrix & $FB$ & feedback \\
		$g$ & Heat exchanger function & $FB+FF$ & feedback-plus-feedforward \\
		$h$ & Specific enthalpy [J kg\textsuperscript{-1}] & $FF$ & feedforward \\
		$J$ & Objective function & $f$ & saturated liquid \\
		$\bm K$ & Controller gain matrix & $forc$ & forced response \\
		$L$ & Heat exchanger length [m] & $fr$ & free response \\ 
		$\dot{m}$ & Mass flow [kg s\textsuperscript{-1}] & $g$ & saturated vapour \\
		$\bm{NLF}$ & Non-linear function & $he$ & heat exchanger \\
		$NTU$ & Number of Transfer Units & $in$ & inlet \\
		$N$ & Compressor speed [Hz] & $inv$ & inverse \\ 
		$N_c$ & Control horizon & $is$ & isentropic \\
		$N_p$ & Prediction horizon & $max$ & maximum \\ 
		$n$ & Element defining the error model & $mid$ & mid-point \\
		$P$ & Pressure [bar] & $min$ & minimum \\
		$\bm Q$ & Error weighting parameter & $opt$ & optimal \\
		$\dot{Q}$ & Cooling power [W] & $out$ & outlet \\ 
		$q$ & Vapour quality & $PNMPC$ & Practical Non-linear MPC \\
		$\bm R$ & Control weighting parameter & $pst$ & past and current values \\
		$S_t$ & Compressor parameter [m\textsuperscript{3}] & $predict$ & predicted values \\
		$s$ & Specific entropy [J kg\textsuperscript{-1} K\textsuperscript{-1}] & $ref$ & reference \\
		$T$ & Temperature [K] & $sec$ & secondary fluid \\
		$T_{SC}$ & Degree of subcooling [K] & $sh$ & superheated vapour section \\ 
		$T_{SH}$ & Degree of superheating [K] & $surr$ & compressor surroundings \\
		$t$ & Time & $sc$ & subcooled liquid section \\
		$UA$ & Overall heat transfer coefficient [W K\textsuperscript{-1}] & $tot$ & total phase change \\
		$\bm u$ & Input vector & $tp$ & two-phase section \\
		$u$ & Distance to the controllability subspace & $trnsf$ & heat transfer \\
		$V_R$ & Heat exchanger volume [$ \text{m}^{\text{3}}$] & $valve$ & expansion valve \\
		$v$ & Specific volume [m\textsuperscript{3} kg\textsuperscript{-1}] & & \\
		$\dot{W}$ & Compressor power [W] & & \\	
		$\bm w$ & Input and disturbance vector & & \\ 
		$\bm x$ & State vector & & \\
		$\bm y$ & Output vector & & \\
		$\bm Z$ & Coefficient matrix & & \\
	\end{tabu}}
\end{table*}


\section{Introduction} \label{Introduccion}

Mechanical compression refrigeration systems represent the leading technology in the market of industrial and residential cooling applications \cite{kalliski2019support}. Energy consumption of Heating, Ventilating, and Air Conditioning (HVAC systems) is reported to represent approximately 30\% of total energy consumed around the world \cite{jahangeer2011numerical}. In commercial and domestic applications, HVAC are claimed to be responsible for more than 40\% of overall electricity consumption \cite{kalkan2012solar}, whereas about 20\% of total energy consumption in buildings in Europe is devoted to mechanical cooling \cite{InternationalEnergyAgency2018}. Among commercial buildings, the energy intensity of a medium-sized supermarket is estimated to be over 500 kWh/m\textsuperscript{2} a year in the USA \cite{US_EPA2}, and around 60\% is related to refrigeration systems \cite{suzuki2011analysis}. These facts, as well as the pressing shortage of non-renewable sources of energy and the difficult progress of sustainable technologies \cite{tomasrodriguez2019modelado}, involve that the development of novel control techniques for refrigeration systems achieving higher energy efficiency becomes a priority. 

\begin{figure}[h]
	\centering
	\centerline{\includegraphics[width=8cm,angle=0]{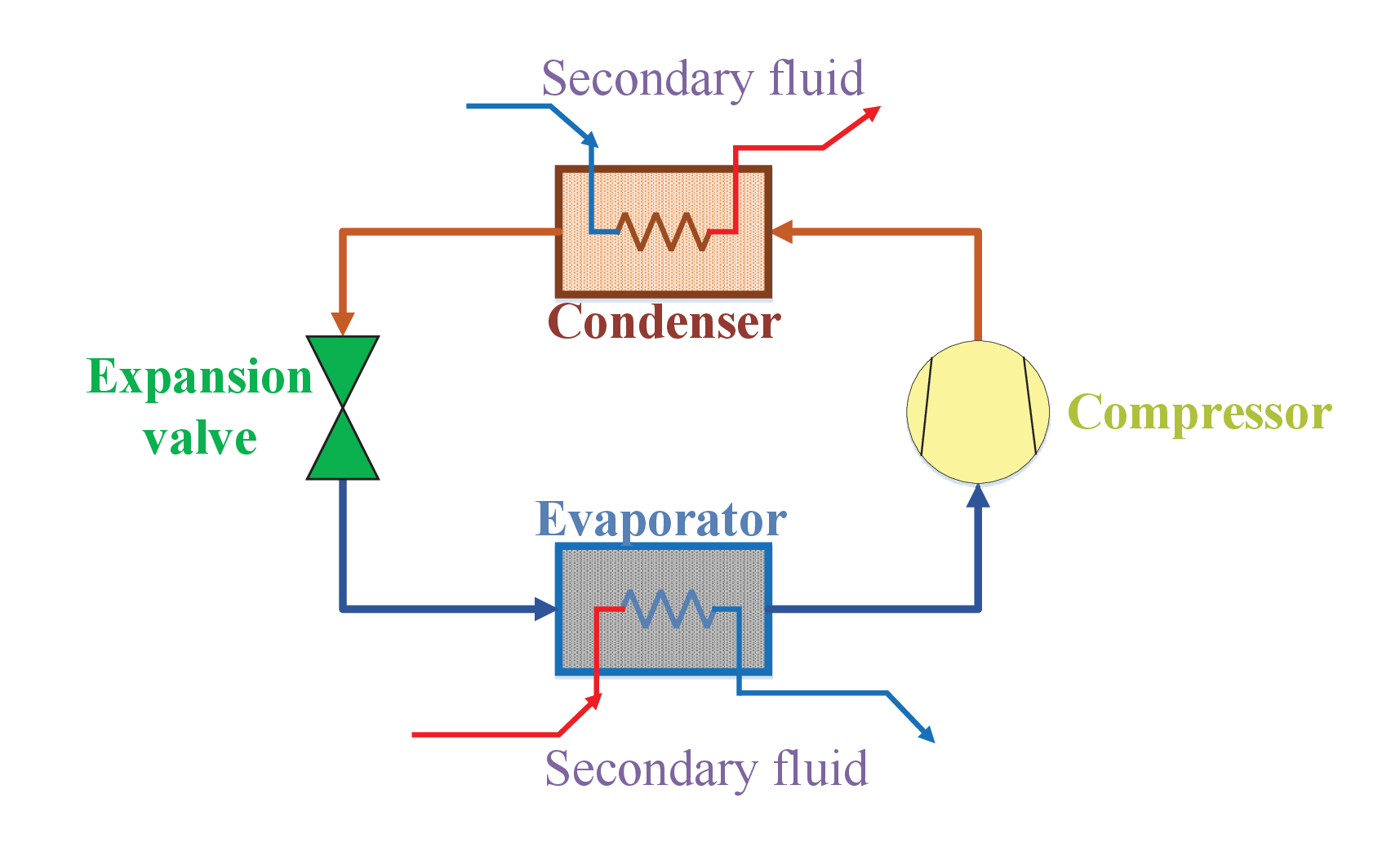}}
	\caption{Mechanical compression system} 
	\label{fig_sist_refrigeracion}
\end{figure}

A canonical mechanical compression refrigeration cycle is illustrated in Fig. \ref{fig_sist_refrigeracion}, including its main components: the evaporator, the condenser, the electronic expansion valve, and the variable-speed compressor. After compression, the refrigerant in vapour state at high temperature and pressure enters the condenser and transfers heat to the secondary fluid, which is usually air. Then, the refrigerant, still at high pressure but in liquid or two-phase state, expands at the expansion valve and then it enters the evaporator, where it removes heat from the secondary fluid due to its low pressure and temperature. As a result, the refrigerant evaporates and becomes superheated vapour. Finally, it enters the compressor, where its pressure and temperature are once again increased.

Let consider that a continuous flow enters the evaporator as the secondary fluid, and that a certain cooling power is intended to be generated at the evaporator to reduce the temperature of this fluid, whose mass flow is assumed to be defined externally and to act as a measurable disturbance for the refrigeration cycle. Likewise, the inlet temperature of the secondary fluid to be cooled is also assumed as a measurable disturbance. Therefore, the refrigeration cycle is intended to provide a certain temperature gradient to the secondary fluid, since the inlet temperature and the mass flow are externally determined. Notice that this is not the conventional case studied in the refrigeration literature, where the cooling demand is intended to be satisfied merely as a thermal power without imposing any constraint on the evaporator secondary mass flow, which turns out to be an additional manipulated input. The same assumption is considered for the condenser, and consequently only the valve opening $A_v$ and the compressor speed $N$ are assumed to be the manipulated variables. 

While the main objective of the refrigeration cycle is to meet the user-defined cooling demand $\dot{Q}_e$, a secondary goal would be to do it with utmost efficiency. The latter is usually rated in refrigeration by the Coefficient of Performance ($COP$), defined for the canonical cycle as shown in Eq. \eqref{eq_COP}, where $\dot{W}_{comp}$ is the compression power supplied at the compressor.

\begin{equation}
	\begin{aligned}
		COP = \dfrac{\dot{Q}_e}{\dot{W}_{comp}} \\
	\end{aligned}
	\label{eq_COP}
\end{equation}

\begin{figure*}[h]
	\centering
	\centerline{\includegraphics[width=13cm,angle=0]{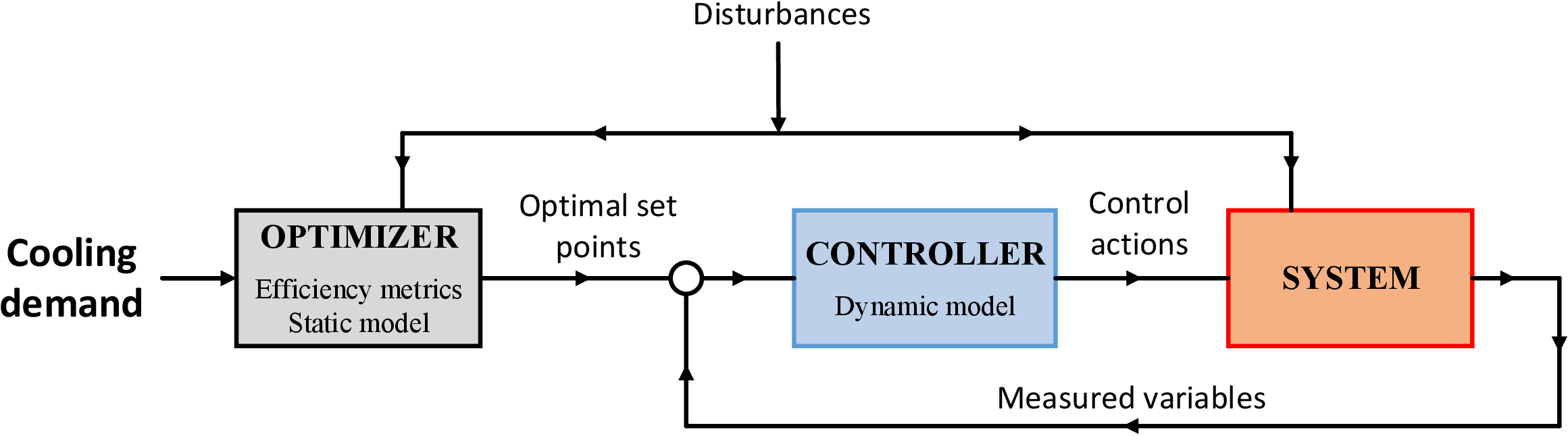}}
	\caption{Classical optimization and control structure} 
	\label{fig_sistema_control}
\end{figure*}

Let a control strategy such as that illustrated in Fig. \ref{fig_sistema_control}, where three key elements must be studied: modelling, optimization, and control design. Concerning modelling, the steady-state features of the refrigeration cycle must be modelled to address the optimization stage, but the dynamic behaviour must be also consider when designing the controller. Regarding optimization, the optimal cycle can be only achieved if suitable references are provided to the controller; these set points are intended to be computed by a model-based optimization stage, where the steady-state model is an essential component to compute accurate optimal operating points. Finally, the controller itself must be designed in such a way that the intended optimal set points are realised. In this work, these three elements (modelling, optimization, and control design) are analysed and discussed. 
 
Regarding dynamic modelling, special attention has been paid to the heat exchangers, since their dynamics are usually dominant. One of the most successful methods is the \emph{moving boundary} approach (MB) \cite{Liang2010}, which divides the heat exchanger into several sections according to the refrigerant states, whose lengths are state variables. Simulation issues arising when any section disappears or reappears have been addressed in the \emph{switched moving boundary} model (SMB) \cite{BINLI}. It has been shown that, despite greatly reducing the computational load, the \emph{moving boundary} methodology achieves similar accuracy to finite-volume formulation \cite{pangborn2015comparison}. Even though, the model order and numerical complexity turn out to be high enough to prevent this technique from being used within model-based control techniques. In this work, the intrinsic dynamics of all SMB model state variables and typical volumes of heat exchangers in the market are all taken into account to propose a simplified third-order state-space model of the refrigeration cycle, whose computational efficiency greatly outperforms the original SMB model. Moreover, while in \cite{bejarano2017difficulty,bejarano2017suboptimal} only a few model representations have been dealt with, in this work a more general model is provided, since all evaporator configurations are tackled. Thus, a more robust dynamic model is provided, able to correctly describe the system dynamics for a wide variety of cooling demand. This reduced-order model is one of the main contributions of this work, since its simplicity and robustness are key factors to ease its use within advanced control strategies.

Concerning optimization, the evaporator efficiency mainly depends on heat transfer. It is well-known that it is much greater when phase change happens, thus if the refrigerant was two-phase along the whole evaporator, energy efficiency would be maximal. However, this is not achievable in practice, since a minimum degree of superheating ($T_{SH}$) must be held at the compressor intake to ensure that it operates safely. This temperature difference is usually preferred to be small, since the higher $T_{SH}$, the lower $COP$. Nevertheless, potential energy savings using optimal references for the cycle pressures have been shown not to be negligible \cite{jain2015exergy}. While the conventional set point selection merely opts for a constant low reference for $T_{SH}$, some works have addressed optimization procedures to compute the optimal operating point of the refrigeration while meeting a certain cooling demand. As an example, an exergy-based global optimization is implemented in \cite{jain2012framework}, which computes not only the key variables of the optimal cycle and the references to be imposed to the controller, but also the steady-state values of the control inputs corresponding to that optimal cycle, suitable to be used in a feedforward control strategy. In this work, a similar procedure is applied, but in this case the $COP$ is considered as energy efficiency metrics. However, unlike other works, only two manipulated variables are considered to be available ($N$ and $A_v$), while the secondary flow variables are assumed as disturbances. Through the optimization procedure, the minimum number of variables that unequivocally define the cycle is shown to be three, which agrees with other related works \cite{jensen2007optimal,jain2011thermodynamics}. The non-linear steady-state model stemming from the reduced-order model previously developed is included within the optimizer. The optimization procedure also allows to impose other operating constraints, such as operating pressures and the minimum desired $T_{SH}$. A previous work by the authors has posed the optimization problem with similar objectives and steady-state models \cite{bejarano2017difficulty}, but imposing such technological and operating constraints makes the obtained results closer to realistic points and introduces some relevant differences. Simulation results featuring the obtained optimal points for a wide range of cooling demand draw attention to potential efficiency increase which cannot be achieved by always holding low $T_{SH}$ for all cooling powers, while the influence of the technological and operating constraints is also highlighted.

Regarding control, the main problems to be faced when regulating refrigeration systems are high coupling and strong non-linearity. The conventional approach is based on imposing two set points: the cooling demand, usually expressed as the outlet temperature of the evaporator secondary fluid $T_{e,sec,out}$, and a low reference for $T_{SH}$ to ensure high energy efficiency and safe compressor operation. However, optimal operation is not achieved nor even intended. Nevertheless, there are other control strategies where optimal cycles are at least aimed. On the one hand, one example of model-based strategy is studied by Jain \cite{jain2013thermodynamics}, where a multivariable feedback-plus-feedforward controller is proposed to achieve the optimal operating point previously computed by the offline global optimization described in the previous paragraph \cite{jain2012framework}. On the other hand, extremum-seeking control (ESC) strategies implement a model-free online optimization \cite{burns2015realtime}. However, none of the previous control strategies actually achieves globally optimal cycles, but they turn out to attain suboptimal solutions. 

In this work, using the previously derived model and the results of the global optimization, the control problem is in depth analysed. It turns out to be an underactuated system, since only two manipulated variables are available, while up to three variables are necessary to completely define the optimal cycle. Furthermore, other works have already shown that, even adding more controllable inputs to the problem, the control problem continues to be underactuated, unless the topology is modified \cite{jensen2007optimal}. Then, the control problem is simplified considering only $N$ and $A_v$ as manipulated variables. The controllability of the reduced-order model previously proposed is studied, using both the classical approach and a numerical graphic study. In order to highlight the difficulties arising when trying to globally optimally control these systems, a novel technique called Practical Non-linear Model Predictive Control (PNMPC) \cite{plucenio2007practical} is applied. This approach, whose recursive feasibility, stability, and robustness have partially been studied in some related works \cite{allan2017inherent,gesser2018robust}, allows to simplify the adaptation of well-known linear Model Predictive Control (MPC) strategies to highly non-linear systems. The steady-state performance and dynamic behaviour achieved by applying the proposed controller are analysed and compared to those provided by the multivariable feedback-plus-feedforward controller previously mentioned \cite{jain2013thermodynamics}, providing some justifications of the controllability issues pointed out by the previous studies.

The remainder of the work is organised as follows: the original SMB model and the state reduction which leads to the proposed reduced-order model are described in Section \ref{Modelling}. The global model-based optimization is described in detail in Section \ref{Optimization}, where some simulation results showing potential efficiency increase for a certain range of cooling demand are also presented. Section \ref{Controllability} studies the controllability of the system using the reduced-order model described in Section \ref{Modelling}. In Section \ref{Control} some control results are provided and discussed, where the controllability problems arise and define the achievable energy efficiency. Finally, Section \ref{Conclusions} draws the main conclusions and expresses some future work.


\section{System modelling} \label{Modelling}

The SMB approach is taken as a starting point when modelling the dynamic behaviour of the heat exchangers \cite{BINLI}. This approach consists in dividing the heat exchanger length into some sections, according to the refrigerant state: subcooled liquid, two-phase fluid, and superheated vapour. The lengths of such sections are time-varying, since they depend on the control inputs and the remaining system states themselves. The original formulation by Li and Alleyne \cite{BINLI} models the heat transfer considering also the wall temperatures of all sections as state variables, but if the \emph{effectiveness-NTU} method \cite{bergman2011fundamentals} is applied, a smaller set of 10 state variables is enough to describe the dynamics of the evaporator and the condenser, assuming that the dynamics of the expansion valve and the compressor are fast enough to be considered negligible in comparison with those of the heat exchangers. These state variables are indicated in Eq. \eqref{eq_vector_estados}, where $\bm x_e$ is the evaporator state vector and $\bm x_c$ refers to the condenser state vector. According to \cite{BINLI}, the condenser state variables are the average enthalpy of the refrigerant at the superheated vapour section $h_{c,sh}$, the condensation pressure $P_c$, the average enthalpy of the refrigerant at the subcooled liquid section $h_{c,sc}$, the normalised lengths of the vapour $\zeta_{c,sh}$ and the two-phase $\zeta_{c,tp}$ sections, and finally the refrigerant void fraction $\bar{\gamma}_c$. In the case of the evaporator, the state vector $\bm{x}_e$ is made up of the normalised length of the two-phase section $\zeta_{e,tp}$, the evaporation pressure $P_e$, the average enthalpy of the refrigerant at the superheated vapour section $h_{e,sh}$, and the refrigerant void fraction $\bar{\gamma}_e$.

\begin{equation} 
	\begin{aligned}
		\bm x_c &= [h_{c,sh} \;\; P_c \;\; h_{c,sc} \;\; \zeta_{c,sh} \;\; \zeta_{c,tp} \;\; \bar{\gamma}_c]^T \in  \mathbb{R}^6 \\
		\bm x_e &= [\zeta_{e,tp} \;\; P_e \;\; h_{e,sh} \;\; \bar{\gamma}_e]^T \in \mathbb{R}^4 \\
		\bm x_{cycle} &= [\bm x_c^T \;\; \bm x_e^T]^T \in \mathbb{R}^{10} \\
	\end{aligned}
	\label{eq_vector_estados}
\end{equation}

The original SMB model proposes a uniform state vector for every heat exchanger, regardless of the model representation or \emph{mode}, defined by the presence or absence of the different sections corresponding to subcooled liquid, two-phase fluid, and superheated vapour. The highly non-linear dynamics of the heat exchangers are expressed in Eq. \eqref{eq_estructura_ecuaciones}, where each \emph{mode} has its own coefficient matrix $\bm Z_i(\bm x_i,\bm w_i) \;\; \forall i=e,c$ and so-called forcing function $\bm f_i (\bm x_i,\bm w_i) \;\; \forall i=e,c$. The expressions of $\bm Z_c$, $\bm f_c$, $\bm Z_e$, and $\bm f_e$ for every condenser/evaporator \emph{mode} may be computed from the original differential equations based on mass and energy balances by removing the wall temperature states \cite{BINLI}. Both heat exchanger models are interconnected through the input/disturbance vectors $\bm w_c$ and $\bm w_e$ shown in Eq. \eqref{eq_input_intercambiadores}, that include the refrigerant mass flows and  enthalpy at every heat exchanger inlet, and the inlet conditions of the corresponding secondary fluid. Please notice that the refrigerant mass flows and enthalpies are computed using the steady-state models of the compressor and the expansion valve, as well as the state vectors $\bm x_c$ and $\bm x_e$. The keen reader is referred to the original work \cite{BINLI} for further details.
 
\begin{equation} 
	\begin{aligned}
		\bm Z_c(\bm x_c,\bm w_c) \; \bm {\dot{x}}_c = \bm f_c (\bm x_c,\bm w_c)\\
		\bm Z_e(\bm x_e,\bm w_e) \; \bm {\dot{x}}_e = \bm f_e (\bm x_e,\bm w_e)\\
	\end{aligned}
	\label{eq_estructura_ecuaciones}
\end{equation}

\begin{equation}
	\begin{aligned}
		\bm{w}_c &= [\dot{m}_{c,sec} \ \  T_{c,sec,in} \ \  \dot{m}_{c,in} \ \  \dot{m}_{c,out} \ \  h_{c,in}]^T\\
		\bm{w}_e &= [\dot{m}_{e,sec} \ \  T_{e,sec,in} \ \  \dot{m}_{e,in} \ \  \dot{m}_{e,out} \ \  h_{e,in}]^T\\
	\end{aligned}
	\label{eq_input_intercambiadores}
\end{equation}

As mentioned in Section \ref{Introduccion}, although the order of the state-space model proposed by Li and Alleyne \cite{BINLI} is much reduced with respect to finite-volume formulation \cite{pangborn2015comparison}, the complexity of a tenth-order model is still enough, preventing this model from being used in advanced control strategies. Therefore, some justified assumptions are made in this work to further reduce the order of the SMB original model, thus simplifying its direct application to control.

First, internal volume of evaporators is getting more and more reduced in the market when compared to condensers, due to limited space, increasing costs of raw material, and progressive replacement of environment-unfriendly refrigerants \cite{Rasmussen2011}. For instance, plate and micro channel heat exchangers are currently dominating the market. Since the lower the volume, the faster the dynamics of the heat exchanger, the dynamics of low-volume evaporators are much faster than the ones of typically-sized condensers. It involves that the former could be neglected and the evaporator could be considered static, becoming the intrinsic dynamics of the condenser the dominant ones. Although this assumption is not general, since it depends on the specific volume ratio between the condenser and evaporator, new systems already incorporate low-volume but high-efficiency evaporators, thus it is expected that such an assumption increases its validity over time. Moreover, since European directives state that chlorofluorocarbons (CFC) refrigerants must be replaced by hydrofluorocarbons (HFC) in new refrigeration facilities \cite{mota2015analysis}, whose energy efficiency is much higher, the evaporator volume has no longer to be high enough to achieve typical cooling demands.

Moreover, apart from the specific volume ratio, the refrigerant density is also a key factor to the intrinsic dynamics of each heat exchanger. In the original SMB model, these are defined mainly by the heat exchanger volume and the average refrigerant density at every section, in addition to its the normalised length. In the case of superheated vapour, density is much smaller than in the case of subcooled liquid, which implies that the intrinsic dynamics of the former sections are expected to be much faster than those of the latter sections. Concerning two-phase fluid, density is determined by vapour quality. As the condenser usually presents a subcooled liquid section, whereas the evaporator only includes two-phase (with not insignificant average vapour quality) and superheated vapour sections, it may involve that the dynamics of the evaporator are faster than the ones of the condenser. Furthermore, refrigerant pressure has also an important influence on density, in such a way that the lower the pressure, the lower the density. It can be checked that, for instance for a typical refrigerant such as R404a (an HFC widely used in refrigeration applications), characteristic densities at typical condenser pressures are about one order of magnitude greater than characteristic densities at typical evaporator pressures. All these factors support the assumption that the dynamics of the evaporator are much faster than the ones of condenser and therefore could be assumed as negligible.

Secondly, concerning the refrigerant inlet and outlet mass flows for a given heat exchanger (for instance the condenser), it is well known that mass equilibrium dynamics are much faster than heat transfer ones. As a result, $\dot{m}_{c,in}$ and $\dot{m}_{c,out}$ become equal at a time scale much faster than that related to heat transfer, and thus a single refrigerant mass flow can be considered.

Given this, the variables exhibiting the dominant dynamics are those shown in Eq. \eqref{eq_estado_reducido_modo1}, all of them related to the condenser.

\begin{equation}
	\begin{aligned}
	\bm x_{c} &= [P_c \;\; \zeta_{c,sc} \;\; h_{c,sc} ]^T\\
	\end{aligned}
	\label{eq_estado_reducido_modo1}
\end{equation}    

Among the condenser state variables, the intrinsic dynamics of $h_{c,sh}$ and $\bar{\gamma}_c$ are the fastest, since the density at the superheated vapour and the two-phase sections is much greater than that at the subcooled liquid section. Therefore, the intrinsic dynamics of the subcooled liquid section turn out to be dominant and that is why the state variables $\zeta_{c,sh}$ and $\zeta_{c,tp}$ have been replaced by the new state variable $\zeta_{c,sc}$ that corresponds to the length of the subcooled liquid section. Notice that this state vector does not depend on the evaporator \emph{mode}, since the latter is statically modelled and its particular configuration with or without $T_{SH}$ has no effect on the dominating system dynamics, beyond the steady-state specific features of the whole cycle.

Therefore, considering the assumptions detailed above, the reduced state vector of the refrigeration cycle is the one shown in Eq. \eqref{eq_estado_reducido_modo1}, while the differential equations are displayed in full in Appendix A. The assumptions described above allow to simplify the original SMB model while describing the dominant dynamics with affordable accuracy and much lower computational load, as will be shown in some simulations included below. 

\begin{figure*}[h] 
	\centering
	\subfigure[Manipulated inputs and refrigerant flow]{
		\includegraphics[width=7.0cm,angle=0] {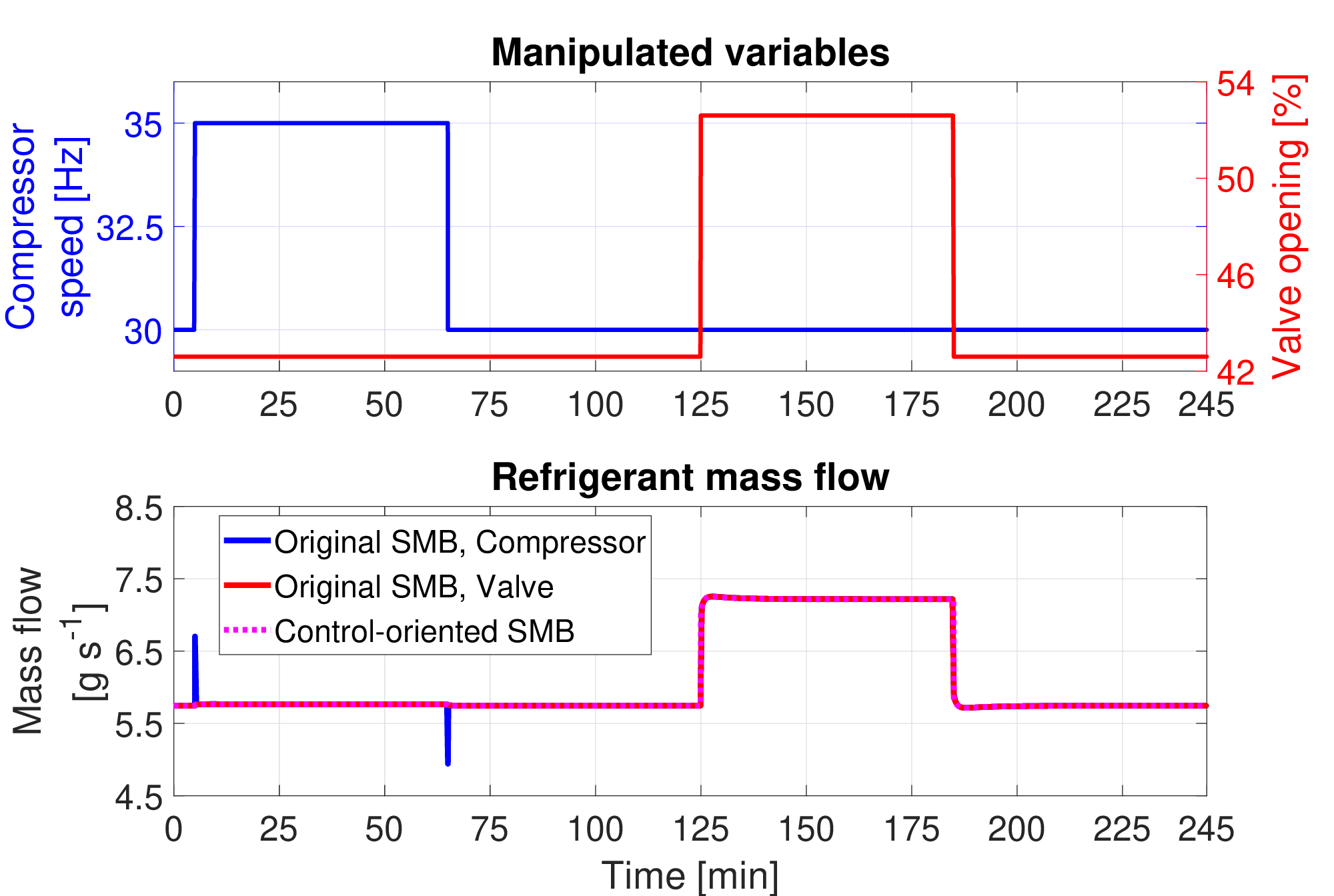}
	}\subfigure[Operating pressures]{
		\includegraphics[width=7.0cm,angle=0] {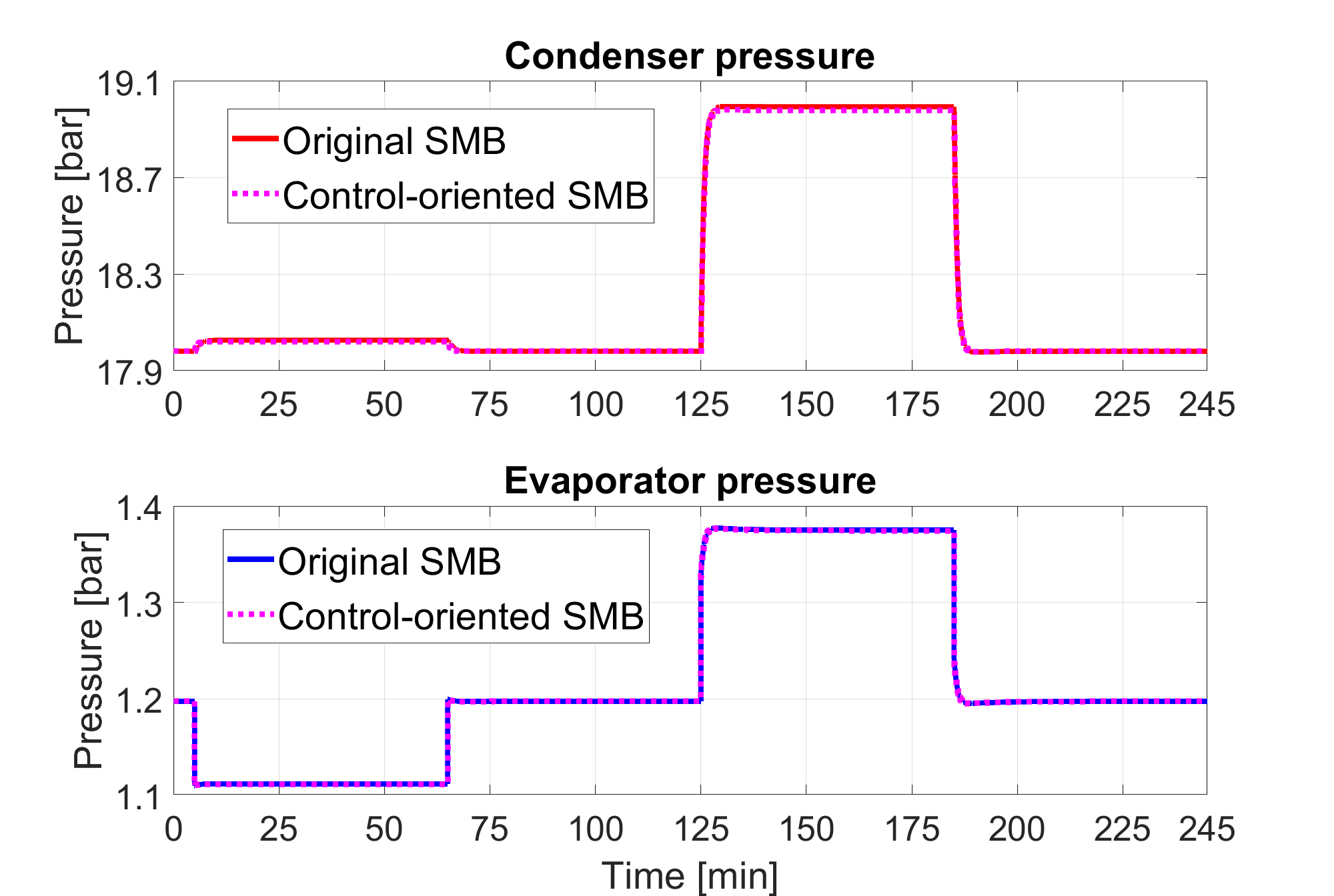}
	}
	\subfigure[Degrees of subcooling and superheating]{
		\includegraphics[width=7.0cm,angle=0] {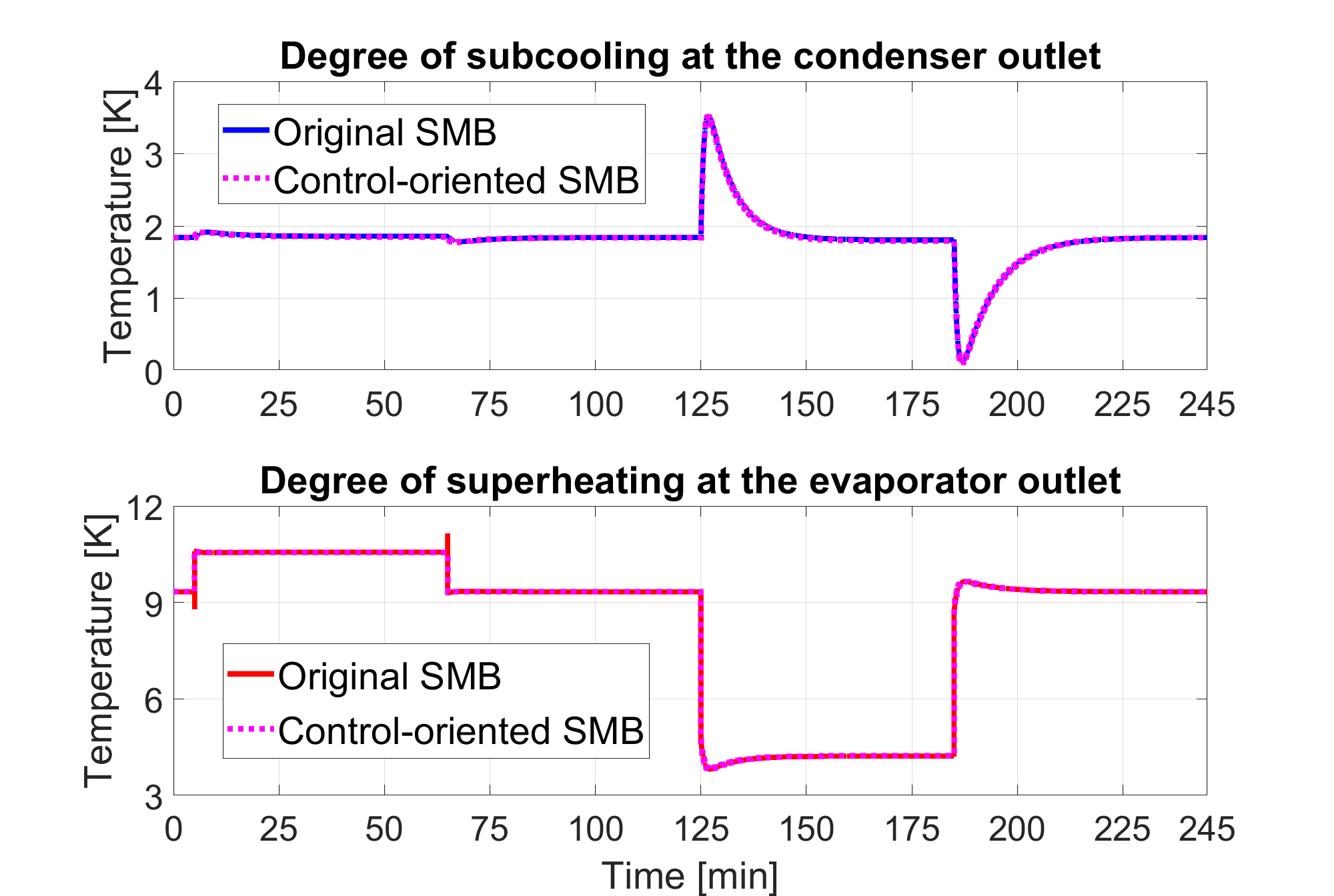}
	}\subfigure[Outlet temperature of the evaporator secondary fluid and Coefficient of Performance]{
		\includegraphics[width=7.0cm,angle=0] {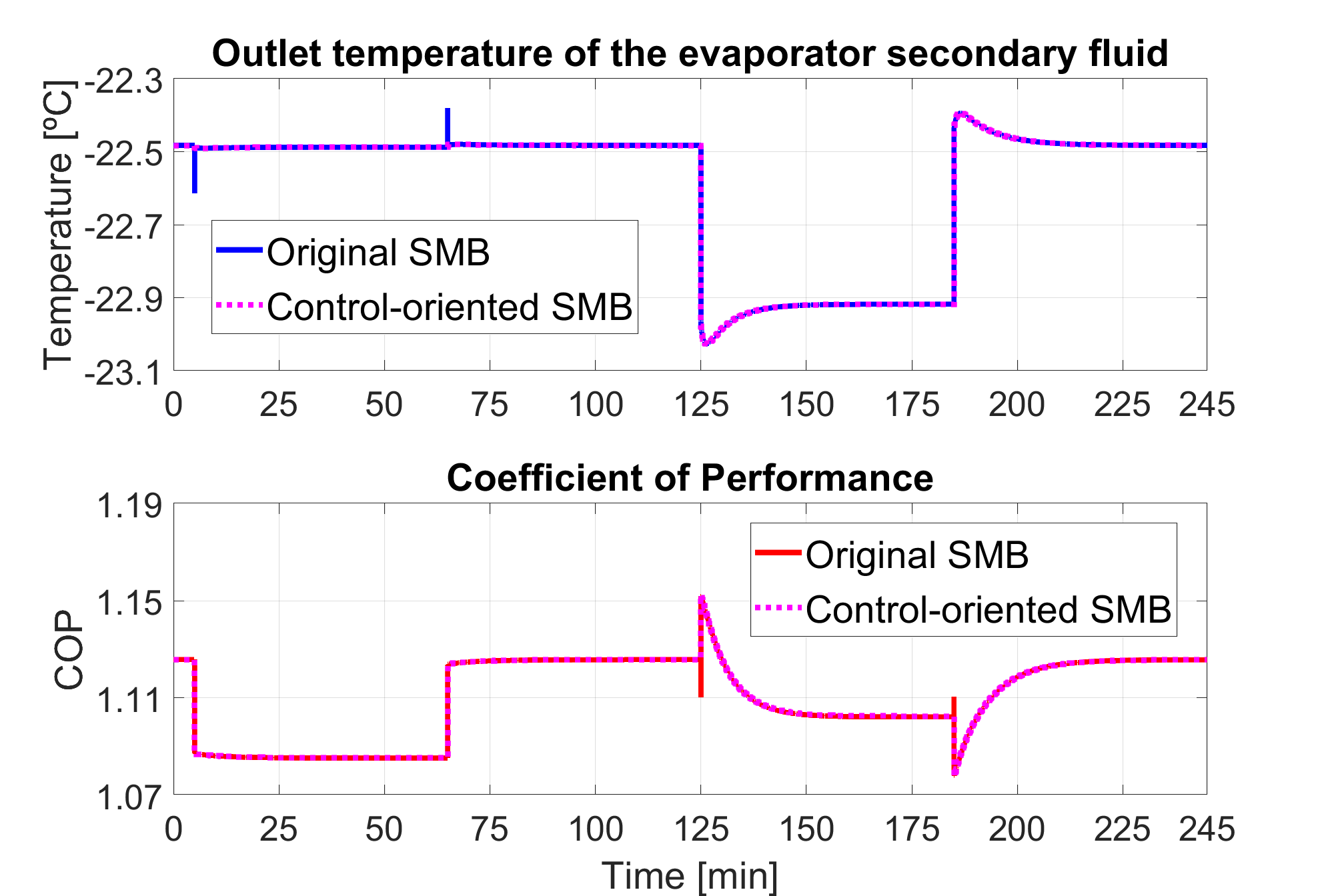}
	}
	\caption{Open-loop simulation results when simple input profiles are applied to both the original SMB model and the reduced-order one} 
	\label{fig_comp_modo1}
\end{figure*}

Fig. \ref{fig_comp_modo1} compares the features of the original SMB model and the reduced order one for simple input profiles. Step changes are applied on $N$ at $t$ = 5 min and $t$ = 65 min, while the expansion valve opening is held constant, whereas step changes are applied on $A_v$ at $t$ = 125 min and $t$ = 185 min, for constant $N$.

As appreciated in Fig. \ref{fig_comp_modo1}, the response of the original SMB model for all cycle variables when a step change is applied on the manipulated inputs is made up of a very fast transient and some slower dynamics before steady state. The fastest transient is mainly caused by the dynamics of the evaporator, while the slower transient turns out to be due to the dominant dynamics of the condenser. Since the proposed reduced-order model does not consider the dynamics of the evaporator, it is very likely that modelling errors are only noticeable just immediately after the step changes. However, the proposed model is still able to describe with reasonable accuracy the dominant dynamics due to the condenser slow states described in Eq. \eqref{eq_estado_reducido_modo1}, as well as the steady-state performance of the complete cycle. Concerning computational load, the computation time of the reduced-order model turns out to be around 2\% of the computation time of the original SMB model for the simulation shown in Fig. \ref{fig_comp_modo1}, while the achieved accuracy is high enough to ease the consideration of this reduced-order dynamic model when addressing model-based control.


\section{Global optimization} \label{Optimization}

\subsection{Problem formulation} \label{Optimization_ProblemStatement}

\begin{figure}[h] 
	\centering
	\centerline{\includegraphics[width=8cm]{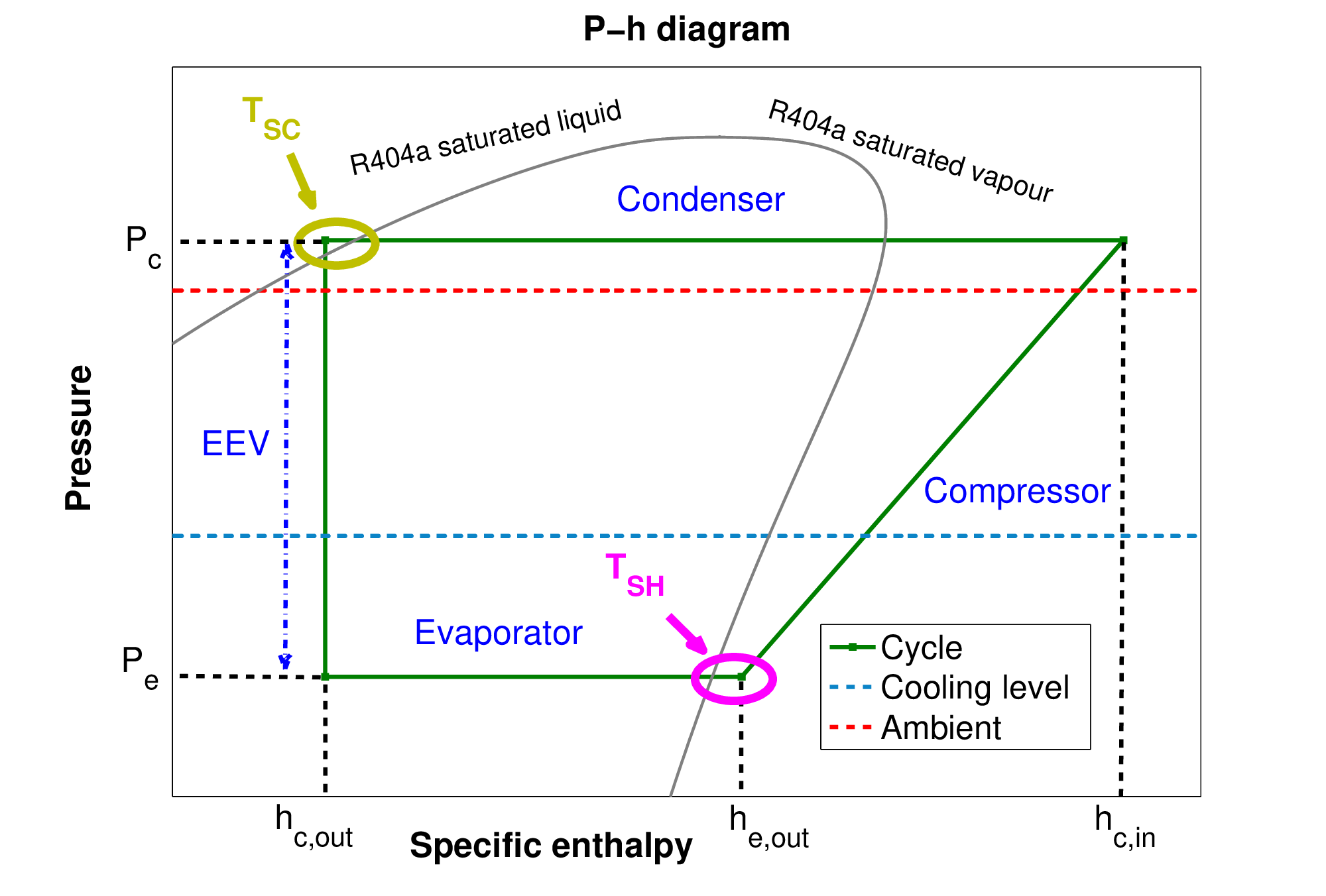}}
	\caption{Pressure-specific enthalpy chart of a canonical mechanical compression cycle} 
	\label{fig_ph_diagram}
\end{figure}

It was stated in Section \ref{Introduccion} that there are two main approaches when it comes to ensuring high energy efficiency performance of mechanical compression refrigeration systems. They can be described through the pressure-specific enthalpy chart (P-h diagram) of a canonical cycle as that shown in Fig. \ref{fig_ph_diagram}. The conventional set point selection merely opts for a constant low reference for $T_{SH}$, in addition to the set point devoted to cooling demand (for instance $T_{e,sec,out}$). This approach is safe for the compressor and it has been shown to lead to high efficiency performance for a certain range of cooling demand. However, if the objective is to ensure maximum energy efficiency for the whole range of cooling demand, a global optimization has been proposed to compute the optimal cycle which meets the cooling demand while maximizing a given energy efficiency metrics, for a known facility and some input conditions of the secondary fluids. An optimization framework also allows to impose some technological constraints or to merely limit the operating range according to manufacturers' recommendations, as well as imposing a given desired value for the minimum admissible $T_{SH}$. 

Anyway, before addressing the optimization, it is necessary to find out the decision variables. If only $N$ and $A_v$ are assumed to be manipulated inputs, thus considering the variables related to the secondary fluids as disturbances, it has been shown that the minimum set which defines a given refrigeration cycle is made up of three variables \cite{jensen2007optimal}. Therefore, the decision variable set could be defined as indicated in Eq. \eqref{eq_1C1R_conjunto_minimo}.

\begin{equation}
	\begin{aligned}
	\bm \chi_{cycle} = [h_{e,out} \;\; h_{c,out} \;\; \dot{m} ] ^T \\
	\end{aligned}
	\label{eq_1C1R_conjunto_minimo}
\end{equation}

This particular choice comes handy when computing the cooling power provided by the cycle according to Eq. \eqref{eq_potencia_frigorifica}, since it turns out to depend only on the decision variables, provided that $h_{e,in}$ is assumed to match $h_{c,out}$, as shown in Fig. \ref{fig_ph_diagram}.  

\begin{equation}
	\begin{aligned}
	\dot{Q_e} = \dot{m}(h_{e,out}-h_{e,in}) = \dot{m}(h_{e,out}-h_{c,out})\\
	\end{aligned}
	\label{eq_potencia_frigorifica}
\end{equation}

As mentioned above, the first and main aim of the refrigeration cycle is to meet the cooling demand. Since this is an external requirement, a constraint between the three decision variables is imposed by Eq. \eqref{eq_potencia_frigorifica}, which leads to a two degree-of-freedom (dof) optimization. The remaining constraints to be imposed may include the following ones:

\begin{itemize}
	\item Maximum and minimum working pressures and pressure ratio.
	\item Minimum temperature gradient between the refrigerant and the secondary fluid at the heat exchangers, according to manufacturers' recommendations.
	\item Minimum admissible $T_{SH}$, to shield the compressor from refrigerant liquid droplets.
	\item Maximum and minimum values of $N$ and $A_v$.
\end{itemize}

\subsection{Non-linear static model}

A non-linear static model of the complete cycle must be included within the optimizer to define if a given decision variable set involves a feasible cycle, considering the constraints described in Section \ref{Optimization_ProblemStatement}. The steady-state model is made up of single static models of all components of the cycle, detailed by Rodríguez \emph{et al.} \cite{rodriguez2017parameter} in a work devoted to identification. The input/output variables, equations to be applied to compute the outputs from the inputs, and the order of implementation when considering the decision variable set shown in Eq. \eqref{eq_1C1R_conjunto_minimo} are described in Eq. \eqref{eq_modelo_valv_implicit}--\eqref{eqs_modelo_cond_implicit}. In particular, as far as the heat exchangers are concerned, the \emph{effectiveness-NTU} method is applied \cite{bergman2011fundamentals}. Moreover, the thermal capacity ratio $C$ in case of phase change is assumed to be zero, since the thermal capacity $\dot{m}c_p$ of a phase-changing fluid is greater than that of any secondary fluid and indeed assumed to be infinite \cite{bergman2011fundamentals}. The experimental plant located at the Department of Systems Engineering and Automatic Control of the University of Seville, Spain \cite{bejarano2015multivariable} is considered as the test-bed, where one brazed-plate evaporator and the air-cooled condenser are considered. All thermodynamic properties of the refrigerant and secondary fluids are computed using the \emph{CoolProp} tool \cite{CoolProp}. Indeed, thermodynamic properties are usually computed from other two independent thermodynamic variables. For instance, $\rho_{c,out} = \rho(P_c,h_{c,out})$ means that the density at the condenser outlet $\rho_{c,out}$ is computed through a certain fluid-specific thermodynamic function $\rho = \rho(P,\ h)$, implemented within the \emph{CoolProp} tool, from the two following known thermodynamic variables: the condenser pressure $P_c$ and the specific enthalpy at the condenser outlet $h_{c,out}$. Similar expressions are applied for other thermodynamic functions used in Eq. \eqref{eq_modelo_valv_implicit}--\eqref{eqs_modelo_cond_implicit}, such as $v = v(P,\ h)$, $c_{v} = c_v(P,\ q)$, $c_{p} = c_p(P,\ q)$, $s = s(P,\ h)$, $h = h(P,\ s)$, $T = T(P,\ q)$, $T = T(P,\ h)$, and $h = h(P,\ q)$.

\subsubsection*{Expansion valve}
\noindent Inputs: $P_e$, $P_c$, $\dot{m}$, $h_{c,out}$ \\
Output: $\bm y_{valve} = A_v$ \\
Parameters: $c_{eev}$ \\
Equations to be applied and order of implementation: \\
\begin{equation} 
	\begin{aligned}
		\rho_{c,out} &= \rho(P_c,h_{c,out}), \quad \quad A_v &= \frac{\dot{m}}{c_{eev} \ \sqrt{2 \rho_{c,out}\ (P_c - P_e)}}. \\
	\end{aligned} 
	\label{eq_modelo_valv_implicit}
\end{equation} 

\subsubsection*{Compressor}
\noindent Inputs: $P_e$, $P_c$, $\dot{m}$, $h_{e,out}$, $T_{surr}$ \\
Outputs: $\bm y_{comp} = [h_{c,in}, \; N, \; \dot{W}_{comp}]$ \\
Parameters: $a$, $b$, $c$, $S_t$, $UA_{comp}$ \\
Equations to be applied and order of implementation: \\
\begin{equation} 
	\begin{aligned}
		&v_{e,out} = v(P_e,\ h_{e,out}), \qquad c_{v,e,g} = c_v(P_e,\ q=1), \qquad c_{p,e,g} = c_p(P_e,\ q=1), \\
		&s_{e,out} = s(P_e,\ h_{e,out}), \qquad h_{c,in,is} = h(P_c,\ s_{e,out}), \qquad T_{c} = T(P_c,\ q=1), \\
		&h_{c,g} = h(P_c,q=1), \qquad \;\; c_{p,c,g} = c_p(P_c,q=1), \\
		&N = \dot{m} \, \dfrac{v_{e,out}}{S_t  - c \   \left(\left(\dfrac{P_c}{P_e}\right)^{\frac{c_{v,e,g}}{c_{p,e,g}}} - 1\right)}, \qquad \dot{W}_{comp} = a + b \  \dot{m} \  (h_{c,in,is} - h_{e,out}), \\
		&T_{c,in,is} = T_{c} + \dfrac{h_{c,in,is} - h_{c,g}}{c_{p,c,g}}, \qquad h_{c,in} = h_{e,out} + \dfrac{\dot{W}_{comp} - U A_{comp} (T_{c,in,is} - T_{surr})}{\dot{m}}. \\
	\end{aligned}
	\label{eqs_modelo_compressor_implicit}
\end{equation}

\subsubsection*{Evaporator}
\noindent Inputs: $P_e$, $\dot{m}$, $h_{e,out}$, $\dot{m}_{e,sec}$, $T_{e,sec,in}$ \\
Outputs: $\bm y_{e} = [h_{e,in}, \; T_{SH}, \; T_{e,sec,out}, \;	\zeta_{e,tp}]$ \\
Parameters: $A_{e,trnsf}$, $L_e$, $\alpha_{e,sh}$, $\alpha_{e,tp}$ \\
Equations to be applied and order of implementation: \\
\begin{equation} 
	\begin{aligned}
		&h_{e,g} = h(P_e,\ q=1), \qquad \qquad T_{e} = T(P_e,\ q=1), \\
		&\dot{Q}_{e,sh} =\dot{m} \, (h_{e,out}-h_{e,g}), \qquad \; C_{e,sh} = \dfrac{(\dot{m} c_p)_{min}}{(\dot{m} c_p)_{max}}, \\
		&\left.
			\begin{matrix}
				\begin{aligned}
					UA_{e,sh} &= \alpha_{e,sh} \ (1 - \zeta_{e,tp}) \ A_{e,trnsf} \ L_e \\
					NTU_{e,sh} &= \dfrac{UA_{e,sh}}{(\dot{m} c_p)_{min}} \\
					\varepsilon_{e,sh} &= g_{he}(NTU_{e,sh}, C_{e,sh}) \\
					\dot{Q}_{e,sh}^{'} &= \varepsilon_{e,sh}(\dot{m} c_p)_{min}(T_{e,sec,in} - T_e) \\
				\end{aligned}	
			\end{matrix}
		\right\} \begin{tabular}{l}
					\text{Iterative resolution on $\zeta_{e,tp}$:} \\
		 			\text{a feasible initial value of } \\
		 			\text{$\zeta_{e,tp}$ is assumed, and then it is } \\
		 			\text{iteratively computed such as } \\
		 			\text{$\dot{Q}_{e,sh}^{'}$ matches $\dot{Q}_{e,sh}$} \\
		 		\end{tabular} \\
		&T_{e,sec,mid} = T_{e,sec,in} - \dfrac{\dot{Q}_{e,sh}}{\dot{m}_{e,sec} \ c_{p,e,sec}}, \\
		&UA_{e,tp} = \alpha_{e,tp} \ \zeta_{e,tp} \ A_{e,trnsf} \ L_e, \qquad \qquad NTU_{e,tp} = \dfrac{UA_{e,tp}}{(\dot{m} c_p)_{min}}, \\
		&C_{e,tp} = \dfrac{(\dot{m} c_p)_{min}}{(\dot{m} c_p)_{max}} = 0 \quad \text{by definition where phase change happens}, \\
		&\varepsilon_{e,tp} = g_{he}(NTU_{e,tp}, C_{e,tp}), \\
		&\dot{Q}_{e,tp} = \varepsilon_{e,tp}(\dot{m} c_p)_{min} (T_{e,sec,mid} - T_e), \qquad \qquad h_{e,in} = h_{e,g} - \dfrac{\dot{Q}_{e,tp}}{\dot{m}}, \\
		&T_{e,sec,out} = T_{e,sec,mid} - \dfrac{\dot{Q}_{e,tp}}{\dot{m}_{e,sec} \ c_{p,e,sec}}, \\
		&T_{e,out} = T(P_e, h_{e,out}), \qquad \qquad T_{SH} = T_{e,out} - T_e. \\
	\end{aligned}
	\label{eqs_modelo_evap_implicit}
\end{equation} 

\subsubsection*{Condenser}
\noindent Inputs: $P_c$, $\dot{m}$, $h_{c,out}$, $\dot{m}_{c,sec}$, $T_{c,sec,in}$ \\
Outputs: $\bm y_{c} = [h_{c,in}^{'}, \; T_{SC}, \; T_{c,sec,out}, \;	\zeta_{c,sh}, \; \zeta_{c,tp}]$ \\
Parameters: $A_{c,trnsf}$, $L_c$, $\alpha_{c,sh}$, $\alpha_{c,tp}$, $\alpha_{c,sc}$ \\
Equations to be applied and order of implementation: \\
\begin{equation}
	\begin{aligned}
		&h_{c,f} = h(P_c,\ q=0), \qquad \dot{Q}_{c,sc} =\dot{m}\,(h_{c,f}-h_{c,out}), \qquad C_{c,sc} = \dfrac{(\dot{m}c_p)_{min}}{(\dot{m}c_p)_{max}}, \\
		&\left.
			\begin{matrix}
				\begin{aligned}
					UA_{c,sc} &= \alpha_{c,sc}\ \zeta_{c,sc}\ A_{c,trnsf}\ L_c \\
					NTU_{c,sc} &= \dfrac{UA_{c,sc}}{(\dot{m}c_p)_{min}} \\
					\varepsilon_{c,sc} &= g_{he}(NTU_{c,sc},\ C_{c,sc}) \\
					\dot{Q}_{c,sc}^{'} &= \varepsilon_{c,sc}(\dot{m}c_p)_{min}(T_{c}-T_{c,sec,in}) \\
				\end{aligned}	
			\end{matrix}
		\right\} \begin{tabular}{l}
					\text{Iterative resolution on $\zeta_{c,sc}$:} \\
					\text{a feasible initial value of } \\
					\text{$\zeta_{c,sc}$ is assumed, and then it is } \\
					\text{iteratively computed such as } \\
					\text{$\dot{Q}_{c,sc}^{'}$ matches $\dot{Q}_{c,sc}$} \\
				\end{tabular} \\
		&T_{c,sec,sc,out}=T_{c,sec,in}+\dfrac{\dot{Q}_{c,sc}}{\dot{m}_{c,sec}\ \zeta_{c,sc} \ c_{p,c,sec}}, \\
		&T_{c,out}=T(P_c,\ h_{c,out}), \quad T_{SC}=T_{c} - T_{c,out}, \\
		&h_{c,g} = h(P_c,\ q=1) \qquad \dot{Q}_{c,tp} =\dot{m}\,(h_{c,g}-h_{c,f}), \\
		&C_{c,tp}=\dfrac{(\dot{m}c_p)_{min}}{(\dot{m}c_p)_{max}} = 0 \quad \text{by definition where phase change happens}, \\
		&\left.
			\begin{matrix}
				\begin{aligned}
					UA_{c,tp} &= \alpha_{c,tp}\ \zeta_{c,tp}\ A_{c,trnsf}\ L_c \\
					NTU_{c,tp} &= \dfrac{UA_{c,tp}}{(\dot{m}c_p)_{min}} \\
					\varepsilon_{c,tp} &= g_{he}(NTU_{c,tp},\ C_{c,tp}) \\
					\dot{Q}_{c,tp}^{'} &= \varepsilon_{c,tp}(\dot{m}c_p)_{min}(T_{c}-T_{c,sec,in}) \\
				\end{aligned}
			\end{matrix}
		\right\} \begin{tabular}{l}
					\text{Iterative resolution on $\zeta_{c,tp}$:} \\
					\text{a feasible initial value of } \\
					\text{$\zeta_{c,tp}$ is assumed, and then it is } \\
					\text{iteratively computed such as } \\
					\text{$\dot{Q}_{c,tp}^{'}$ matches $\dot{Q}_{c,tp}$} \\
				\end{tabular} \\
		&T_{c,sec,tp,out} = T_{c,sec,in}+\dfrac{\dot{Q}_{c,tp}}{\dot{m}_{c,sec}\ \zeta_{c,tp} \ c_{p,c,sec}}, \\
		&\zeta_{c,sh} = 1 - \zeta_{c,sc} - \zeta_{c,tp}, \qquad C_{c,sh}=\dfrac{(\dot{m}c_p)_{min}}{(\dot{m}c_p)_{max}}, \qquad UA_{c,sh}=\alpha_{c,sh}\ \zeta_{c,sh}\ A_{c,trnsf}\ L_c, \\
		&NTU_{c,sh}=\dfrac{UA_{c,sh}}{(\dot{m}c_p)_{min}}, \qquad \qquad \varepsilon_{c,sh}=g_{he}(NTU_{c,sh},\ C_{c,sh}), \\
		&\left.
			\begin{matrix}
				\begin{aligned}
					T_{c,in}^{'} &= T(P_c,h_{c,in}^{'}) \\
					\dot{Q}_{c,sh} &= \varepsilon_{c,sh}(\dot{m}c_p)_{min}(T_{c,in}^{'}-T_{c,sec,in}) \\
					\dot{Q}_{c,sh}^{'} &= \dot{m}(h_{c,in}^{'}-h_{c,g}) \\
				\end{aligned}
			\end{matrix}
		\right\} \begin{tabular}{l}
					\text{Iterative resolution on $h_{c,in}^{'}$:} \\
					\text{a feasible initial value of } \\
					\text{$h_{c,in}^{'}$ is assumed, and then it is } \\
					\text{iteratively computed such as } \\
					\text{$\dot{Q}_{c,sh}^{'}$ matches $\dot{Q}_{c,sh}$} \\
				\end{tabular} \\
		&T_{c,sec,sh,out}=T_{c,sec,in}+\dfrac{\dot{Q}_{c,sh}}{\dot{m}_{c,sec}\ \zeta_{c,sh} \ c_{p,c,sec}}, \\
		&T_{c,sec,out}=\zeta_{c,sh} \ T_{c,sec,sh,out}+\zeta_{c,tp} \ T_{c,sec,tp,out}+\zeta_{c,sc} \ T_{c,sec,sc,out}. \\
	\end{aligned}
	\label{eqs_modelo_cond_implicit}
\end{equation}

\subsection{Optimization procedure}

A non-linear programming solver is applied to solve the problem, specifically the MATLAB\textsuperscript{\textregistered} \emph{fmincon} tool \cite{MatlabOTB} and the Sequential Quadratic Programming (SQP) algorithm. The objective of the optimization is to maximize the $COP$, thus the cost function described in Eq. \eqref{eq_objective_optimization} is minimized. 

\begin{equation}
	J = - \ COP = \dfrac{h_{e,in} - h_{e,out}}{h_{c,in} - h_{e,out}}
	\label{eq_objective_optimization}
\end{equation}

\begin{figure*}[h] 
	\centering
	\centerline{\includegraphics[width=15cm,trim = 115 460 115 120,clip]{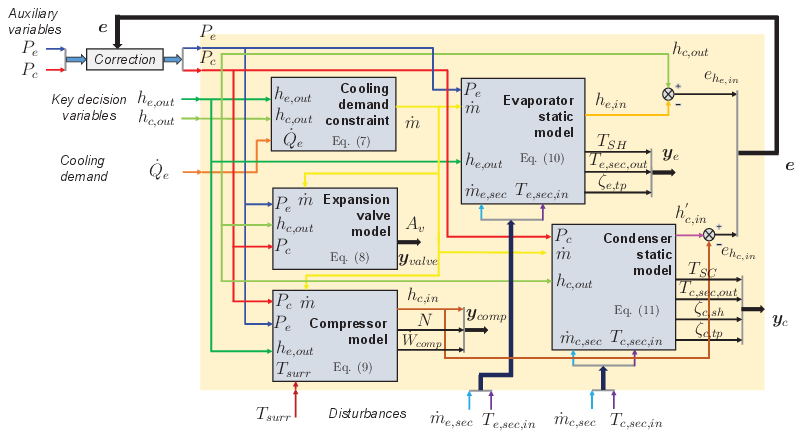}}
	\caption{Integration of the static model of all cycle components into the optimization procedure} 
	\label{fig_modelo_estatico_completo}
\end{figure*}

The procedure described in Fig. \ref{fig_modelo_estatico_completo} has been designed to integrate the resolution of the static model of the whole cycle into the optimizer. A certain variable set $\bm \chi_{cycle}$ is proposed as a feasible \emph{candidate} to the optimum, as long as the cooling demand is met. Given the system topology and the definition of the decision variable set, it is required to know the cycle pressures $P_e$ and $P_c$ to evaluate the models of the expansion valve, compressor, evaporator, and condenser. Therefore, although the conceptual decision variable vector is that defined in Eq. \eqref{eq_1C1R_conjunto_minimo}, the optimization procedure is actually carried out considering up to five decision variables: the three ones included in $\bm \chi_{cycle}$, and the cycle pressures $P_e$ and $P_c$. However, it must be imposed that $h_{e,in}$ matches $h_{c,out}$ to close the cycle, whereas the compressor and condenser models likely provide dissimilar values of $h_{c,in}$, that must also match, as shown in Eq. \eqref{eq_extra_equality_constraints}.

\begin{equation}
	\left.
		\begin{matrix}
			\begin{aligned}
				h_{e,in} &= h_{c,out} \\
				h_{c,in}^{'} &= h_{c,in} \\
			\end{aligned}
		\end{matrix}
	\right\}	
	\label{eq_extra_equality_constraints}
\end{equation}

Therefore, although there are nominally five degrees of freedom in the optimization, the equality constraints shown in Eq. \eqref{eq_extra_equality_constraints} reduce the number of actual decision variables to three, while Eq. \eqref{eq_potencia_frigorifica} still reduces the number of degrees of freedom to two. The optimization algorithm corrects the values of the cycle pressures according to the deviation of the equality constraints, called $e_{h_{e,in}}$ and $e_{h_{c,in}}$ in Fig. \ref{fig_modelo_estatico_completo}, in the search for a feasible solution that minimizes the objective function shown in Eq. \eqref{eq_objective_optimization}. 

Regarding the operating and technological constraints, they are detailed in Eq. \eqref{eq_constraints_optimization}, where $P_{e,max}$ corresponds to the maximum evaporation pressure ensuring a recommended minimum temperature gradient between heat transfer fluids at the evaporator. Similarly, $P_{c,min}$ refers to the minimum condensation pressure ensuring a recommended minimum temperature gradient between heat transfer fluids at the condenser. Eventually, $CR_{min}$ and $CR_{max}$ represent the minimum and maximum compression ratio imposed by the compressor manufacturer, respectively.

\begin{equation}
	\begin{array}{lllll}
	N_{min} & \leqslant & N & \leqslant & N_{max} \\
	A_{v,min} & \leqslant & A_v & \leqslant & A_{v,max} \\
	T_{SH,min} & \leqslant & T_{SH} & \leqslant & T_{SH,max} \\
	 & & P_e & \leqslant & P_{e,max} \\
	P_{c,min} & \leqslant & P_c & & \\ 
	CR_{min} & \leqslant & P_c/P_e & \leqslant & CR_{max} \\
	\end{array}
	\label{eq_constraints_optimization}
\end{equation}

If the optimization procedure fails, it implies that no feasible cycle is able to meet the cooling demand while satisfying all constraints. However, if it finds a solution, that means that, among all feasible cycles, this is the one achieving highest $COP$. Depending on the initial values of the decision variables, the convergence of the numerical algorithm cannot theoretically ensured. Nevertheless, the uniqueness of the cycle pressure solution for a given vector $\bm \chi_{cycle}$ is ensured, since a cycle is unequivocally represented by a three-variable set, as long as $N$ and $A_v$ are the manipulated variables. Therefore, if a given \emph{candidate} $\bm \chi_{cycle}$ is provided by the optimization procedure as the optimum, one can affirm that this cycle is undoubtedly feasible, being the pressures simultaneously obtained the unique pair featuring this particular cycle. 

\subsection{Optimization results}

\begin{figure*}[h]
	\centering
	\subfigure[Optimal $COP$ and $T_{SH}$]{
		\includegraphics[width=6.5cm] {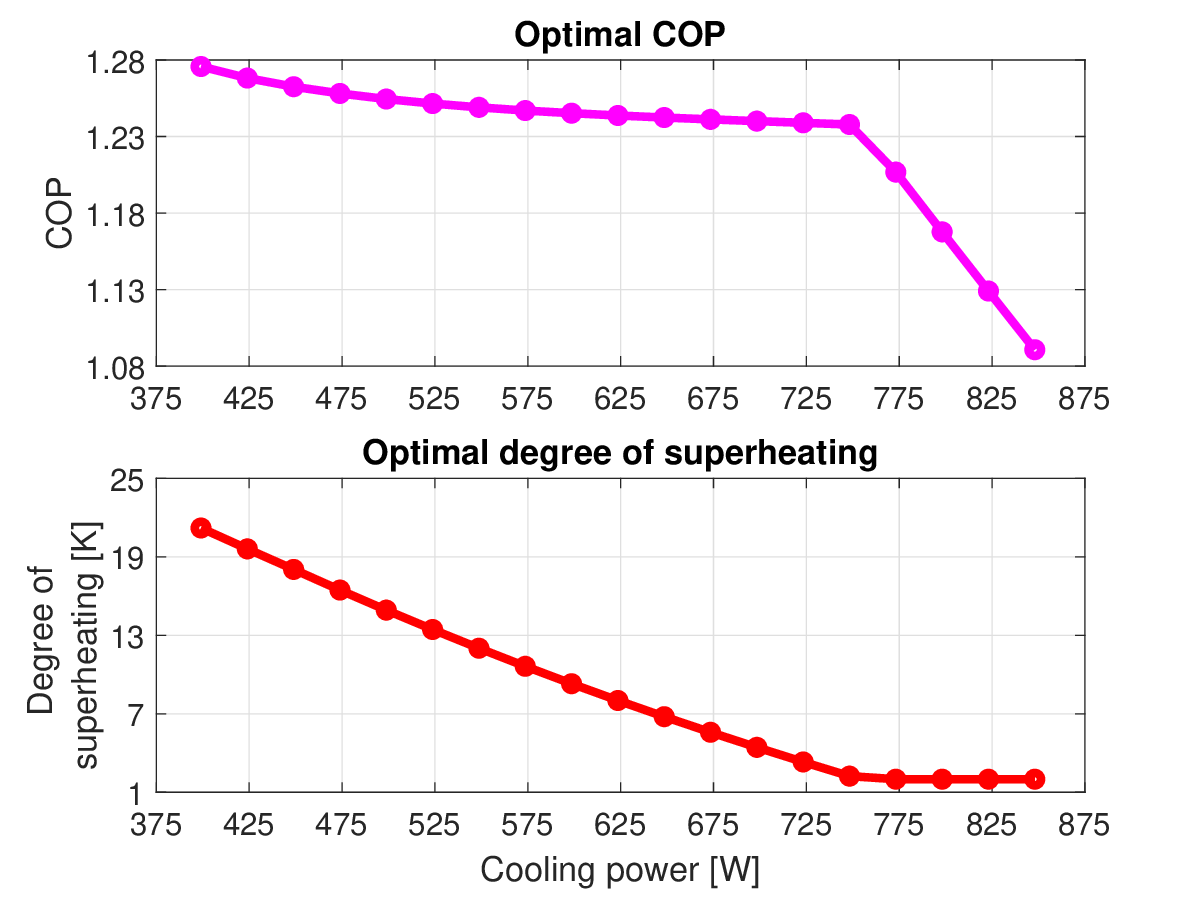}
		\label{fig_optimizacion_a}
	}\subfigure[Operating pressures]{
		\includegraphics[width=6.5cm] {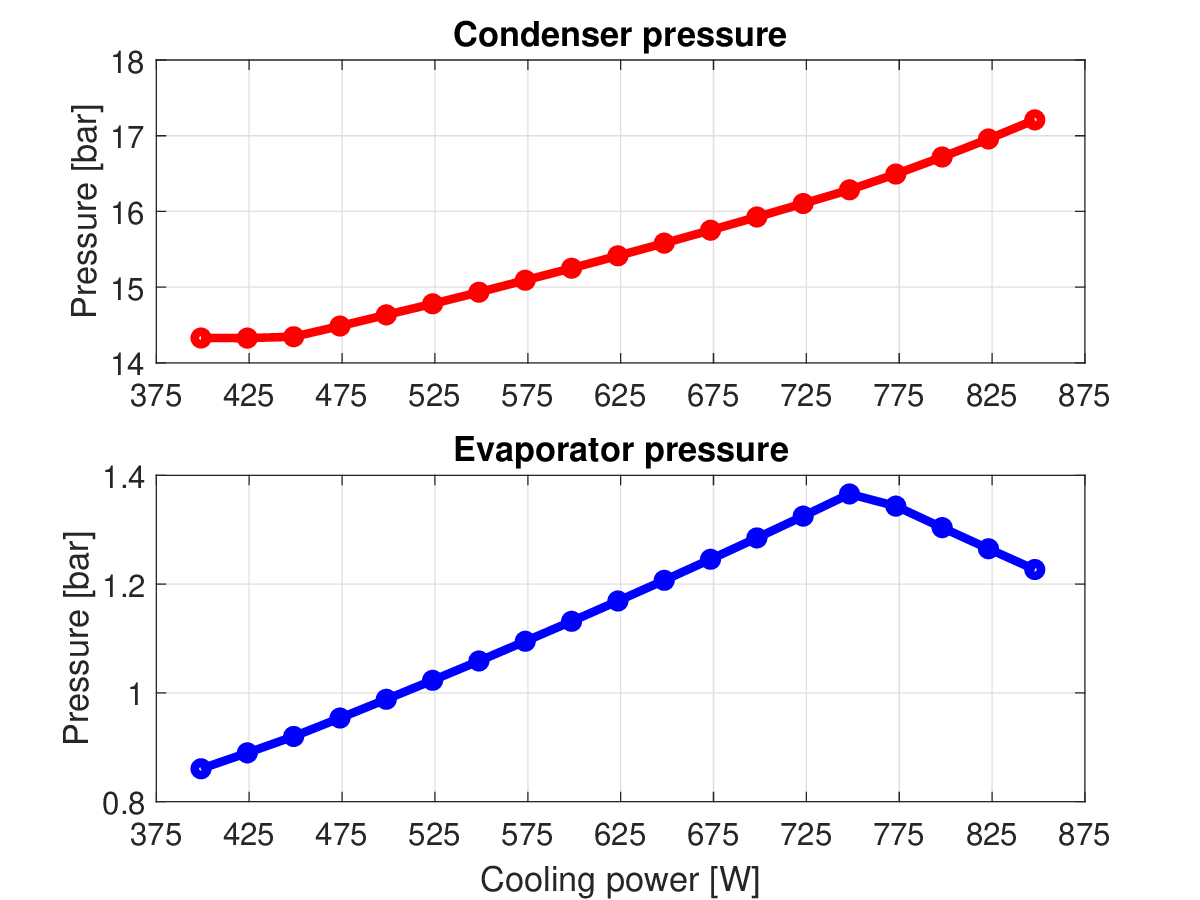}
		\label{fig_optimizacion_b}
	}
	\subfigure[Refrigerant flow or \emph{active charge} \cite{jensen2007optimal}]{
		\includegraphics[width=6.5cm] {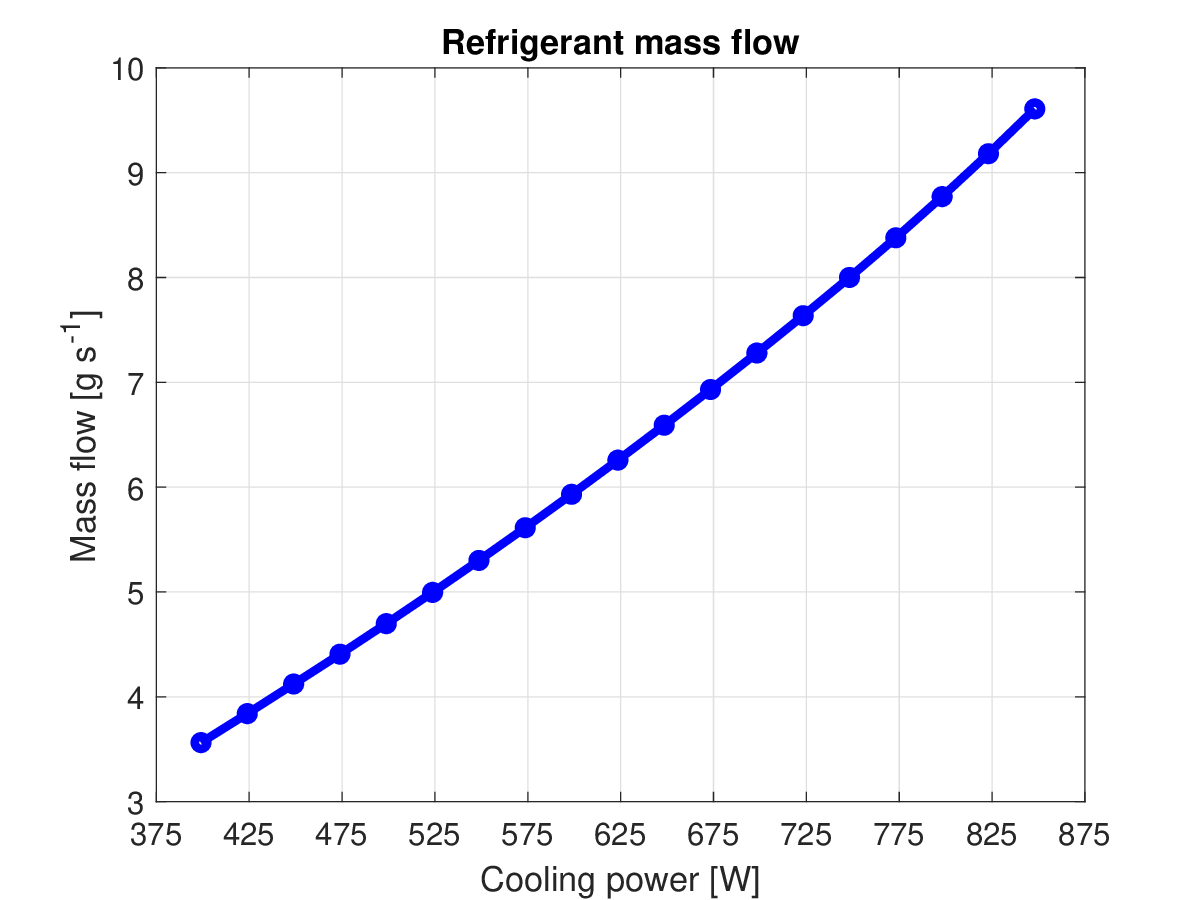}
		\label{fig_optimizacion_c}
	}\subfigure[Steady-state control actions]{
		\includegraphics[width=6.5cm] {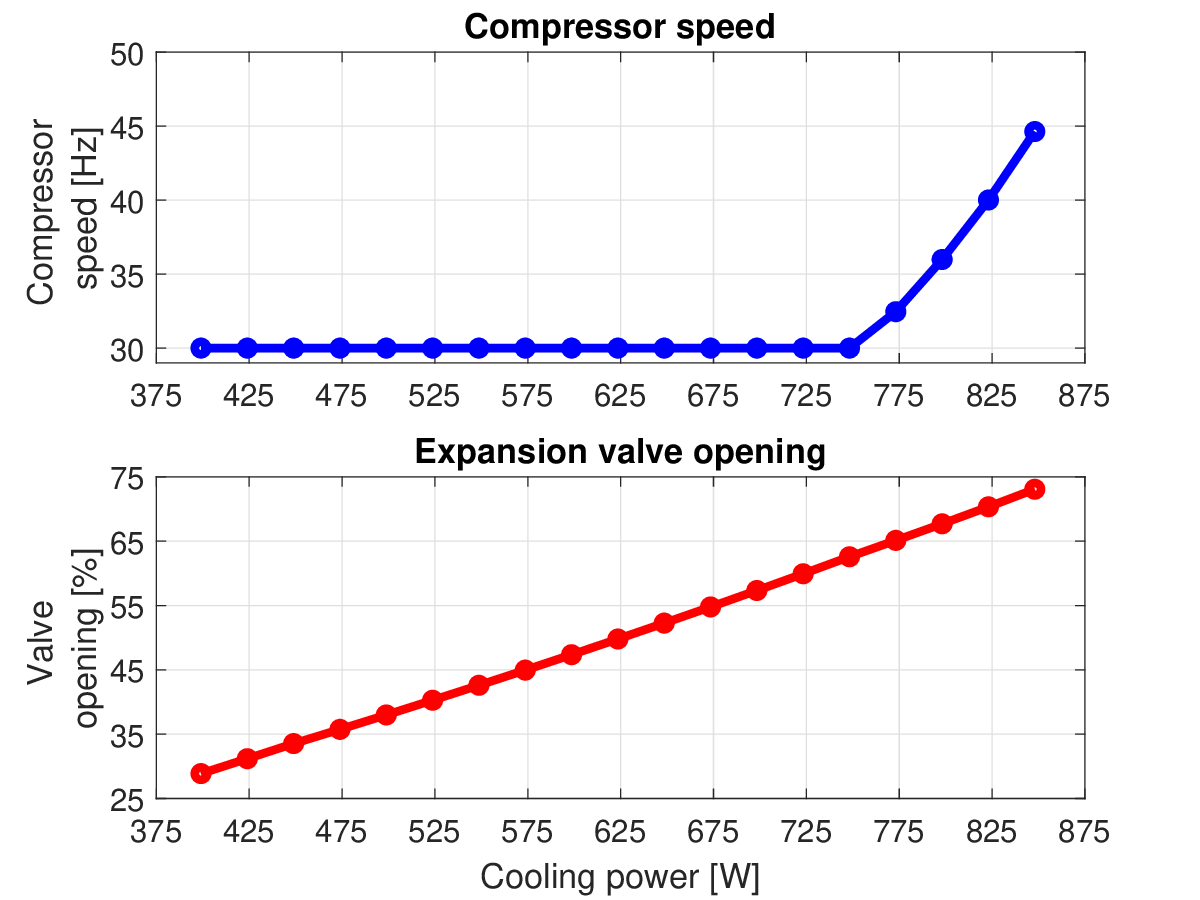}
		\label{fig_optimizacion_d}
	}
	\caption{Optimal cycle features for a range of cooling load}
	\label{fig_optimizacion} 
\end{figure*}

In this subsection some optimization results are represented to illustrate how the optimal cycles change while increasing the cooling demand, for the facility described in \cite{bejarano2015multivariable}. Some key variables, as the $COP$ (subplot \ref{fig_optimizacion_a}), $T_{SH}$ (subplot \ref{fig_optimizacion_a}), $P_e$ and $P_c$ (subplot \ref{fig_optimizacion_b}), $\dot{m}$ (subplot \ref{fig_optimizacion_c}), and finally the steady-state values of $N$ and $A_v$ (subplot \ref{fig_optimizacion_d}) are represented in Fig. \ref{fig_optimizacion}, for a feasible range of cooling power. The tightest constraints imposed are related to the compressor speed range, $N \in [N_{min}, N_{max}] = [30, 50]$ Hz, and the minimum admissible degree of superheating, $T_{SH,min} \geqslant 2$ K.

It can be observed in Fig. \ref{fig_optimizacion_a} that, for low and medium cooling power, the optimal $T_{SH}$ is much above the minimum imposed value $T_{SH,min}$. For those values of the cooling demand, if a lower $T{SH}$ was set as reference for the controller, the efficiency of the achieved cycle would be greatly reduced. This illustrates that, even though the conventional set point selection could lead to high energy efficiency for high cooling load, when considering low and medium cooling demand a model-based optimization is able to provide more efficient set points. Furthermore, when the cooling demand is high enough, the $COP$ is much reduced due to the compressor speed increase shown in Fig. \ref{fig_optimizacion_d}. In these cases, the cooling load can be no longer provided by keeping $N$ to the minimum value $N_{min}$ and merely increasing $A_v$ and thus the refrigerant mass flow $\dot{m}$, represented in Fig. \ref{fig_optimizacion_c}. Regarding the cycle pressures, $P_c$ and $P_e$ increase continuously as $\dot{m}$ does, but $P_e$ starts decreasing just when the $N$ starts increasing, which means that the pressure ratio is also greater for high cooling demand. Eventually, also concerning the cycle pressures, for low cooling demand the most restrictive constraint turns out to be the minimum refrigerant pressure at the condenser. The latter is determined by the ambient temperature (hot reservoir of the cycle) and the minimum temperature gradient between the refrigerant and
secondary fluid recommended by the condenser manufacturer. This kind of technological constraints are very important in practice, since they define the practical operation of the heat exchangers, however they have not been imposed in other similar optimizers \cite{jain2013thermodynamics,bejarano2017difficulty}. 

The global nature of the solutions provided in this subsection for the analysed range of cooling demand cannot be analytically proved, since the high non-linearity of the static model and the use of refrigerant-specific functions to evaluate several thermodynamic properties such as specific enthalpy, pressure, density, viscosity, etc., renders an analytical study impracticable. However, extensive simulations using several optimization methods (pattern search, genetic algorithms, multistart) have been performed and the results provided in Fig. \ref{fig_optimizacion} have been obtained for all methods, that allows to affirm that the cycles given by the optimizer are in this case globally optimal.

The analysed results are intended to be seen as an illustrative example of how efficiency might be improvable for a certain range of cooling load if optimal set point selection is addressed. It goes without saying that the particular numerical results certainly depend on the specific facility and static models, thus no general conclusions should be made. However, the performed simulations suggest consistently that the optimization stage could help achieve higher energy efficiency, especially when the cooling demand is time-varying and it covers a wide range, when compared to the conventional set point selection.


\section{Controllability study} \label{Controllability}

It has been shown in Section \ref{Optimization} that the conventional reference selection might not lead always to optimal operation and it may be necessary to lead the system from one random feasible point to the optimal cycle, given a certain cooling demand. This must be performed by the controller, but there are a number of issues that are worth considering.

Given that a feasible cycle is shown to be completely defined by three independent variables, and the reduced-order dynamic model presented in Section \ref{Modelling} is a third-order state-space model with states related to the condenser, it seems reasonable to study the controllability of such component, since the remaining models are assumed to be static. Moreover, the three independent variables defining the cycle could be selected to match the states of the reduced-order model, in such a way that the full controllability of the whole cycle could be ensured if the three states of the condenser model could be driven from a random point to the optimal operating point.

Firstly, the controllability of the condenser model is studied from the point of view of classical linear theory. The details are included in Appendix \ref{Appendix_LinearControllability}, where it is stated that the condenser model turns out to be underactuated, since only two manipulated variables are available and up to three states must be regulated. Moreover, according to linear control theory, the degree of controllability seems to be two, which suggests that there might be controllability problems when trying to move the full condenser state from a given operating point to another one.

Nevertheless, since the condenser dynamics are highly non-linear, a numerical study based on the phase portrait approach is also addressed. It is worth mentioning that an analytical controllability study could not be tackled, since the non-linear equations of the condenser model include tabulated thermodynamic functions for which common software tools such as CoolProp \cite{CoolProp} do not provide analytic expressions.

This non-linear study is shifted for the sake of readability to Appendix \ref{Appendix_NonLinearControllability}. This study confirms the conclusions provided by the linear analysis concerning the degree of controllability. That means that, from a random initial state, a desired operating point cannot be achieved by manipulating the available inputs, unless the latter falls within the controllable subspace. This constrains the number of reachable cycles for any control strategy, since controllability depends only on the system and not on the control technique.


\section{Control results} \label{Control}

\subsection{Problem statement} \label{Problem_statement}

The main objective of this section is to confirm, highlight, and illustrate the controllability problems found in Section \ref{Controllability} by applying two specific control strategies to the refrigeration cycle. Let assume that a certain cooling demand is intended to be met, for which the optimal cycle is computed using the global optimization procedure described in Section \ref{Optimization}. The controllability issues remarked by the study performed in the Section \ref{Controllability} are highlighted when applying two well-known control techniques: model predictive control and the feedback-plus-feedforward strategy proposed by Jain \cite{jain2013thermodynamics}.  

The facility analysed in Section \ref{Optimization} is considered and a certain cooling demand is selected within the range whose optimal cycles are previously obtained. Thus, cooling power of 600 W is demanded, which involves that $T_{e,sec,out}$ should be -22.4 ºC; the optimal cycle producing such cooling power is characterised by the variables gathered in Table \ref{tab_optimal_cycle}, where the decision variables (those defining the optimization set $\bm \chi_{cycle}$, in blue) and other cycle variables, including control inputs (in magenta) and measurable variables (in green), are gathered. The optimal $COP$ is also added, in brown font.

The decision variable vector $\bm \chi_{cycle}$ can be selected to include a set of pressures and temperatures which can be easily measured, since the only requirement is that it is made up of three independent cycle variables. Thus, $\bm \chi_{cycle}$ is redefined as the measurable variable set $\bm \phi_{cycle}$, as expressed in Eq. \eqref{eq_phi_cycle}.  

\begin{equation}
	\begin{aligned}
	\bm \phi_{cycle} = [P_e  \;\; P_c \;\; T_{e,sec,out}] ^T \\
	\end{aligned}
	\label{eq_phi_cycle}
\end{equation}

\subsection{Practical Non-linear Model Predictive Control} \label{PNMPC}

The PNMPC technique allows to simplify the adaptation of well-known linear MPC strategies to highly non-linear systems, making use of quadratic programming (QP) solvers \cite{plucenio2007practical}. The predicted outputs are intended to be linearly represented with respect to future increments of the control actions, but no equilibrium point is considered. It is assumed that, in a non-linear system, the predicted controlled variables depend on the past manipulated variables $\bm u_{pst}$, the past and current controlled variables $\bm y_{pst}$, and the future increments of the manipulated variables $\bm{\Delta u}$, as indicated in Eq. \eqref{eq_nonlinear}. Notice that  $\bm{NLF}$ is a non-linear function whose complexity is not limited. 

\begin{equation}
	\begin{aligned}
	\bm y_{predict} = \bm{NLF}(\bm y_{pst},\bm u_{pst}, \bm{\Delta u})\\
	\end{aligned}
	\label{eq_nonlinear}
\end{equation}  

In accordance with the approximation proposed in \cite{plucenio2007practical}, the predicted output vector can be rewritten as expressed in Eq. \eqref{eq_GPNMPC}.

\begin{equation}
	\begin{aligned}
	\bm y_{predict} &= \bm y_{fr}  + \bm y_{forc} = \bm y_{fr} + \bm G_{PNMPC} \; \bm{\Delta u} \\
	\bm y_{fr} &= \bm{NLF}(\bm y_{pst},\bm u_{pst}) \\
	\bm G_{PNMPC} &= \dfrac{\partial \bm y_{predict}}{\partial \bm{\Delta u}} \biggr\rvert_{\bm{\Delta u} = \bm 0} 
	\end{aligned}
	\label{eq_GPNMPC}
\end{equation}  

$\bm G_{PNMPC}$ refers to the Jacobian matrix of $\bm y_{predict}$, computed as the gradient of the predicted controlled variables with respect to future increments on the manipulated variables, computed for $\bm{\Delta u} = \bm 0$. In other words, the representation of the predicted controlled variables shown in Eq. \eqref{eq_GPNMPC} is a linearisation of the Taylor series for $\bm{\Delta u} = \bm 0$ where second- and higher-order terms are neglected. Moreover, an additional factor is added to the free response computed using the actual process outputs, in order to obtain a offset-free closed-loop system \cite{plucenio2007practical}. Once computed the predicted output vector in the form shown in Eq. \eqref{eq_GPNMPC}, the future increments on the control inputs can be computed using standard QP algorithms as if the system was linear. 

The PNMPC strategy is applied to the problem presented in Section \ref{Problem_statement}, trying to control the three measurable variables included in $\bm \phi_{cycle}$ and using only $N$ and $A_v$ as manipulated variables. The references for all controlled variables are computed by the optimizer presented in Section \ref{Optimization}. The reduced-order non-linear model described in Section \ref{Modelling} is used, discretised with a sampling time of 5 seconds. The cost function $J_{PNMPC}$ considered at time instant $j$ is shown in Eq. \eqref{PNMPC_J}.

\begin{flalign}
	\begin{aligned}
	&J_{PNMPC} (j) = \\
	&= \sum_{k=1}^{N_p}[\bm \phi_{cycle,ref}(j+k) - \bm \phi_{cycle}(j+k|j)]^T\;\bm Q\;[\bm \phi_{cycle,ref}(j+k) - \bm \phi_{cycle}(j+k|j)] + \\
	&+ \sum_{k=1}^{N_c}\bm{\Delta u}(j+k-1|j)^T \; \bm R\;\bm{\Delta u}(j+k-1|j)
	\end{aligned}
	\label{PNMPC_J}
\end{flalign}

Notice that matrices $\bm{Q}$ and $\bm{R}$ are user-defined tuning parameters. The constraints detailed in Eq. \eqref{eq_constraints_PNMPC} and \eqref{eq_values_constraints_PNMPC} are taken into account, while the tuning parameters are indicated in Eq. \eqref{eq_parameters_PNMPC}.

\begin{equation}
	\begin{array}{lllll}
	\bm{\Delta u}_{min}  & \leqslant &\bm{\Delta u} &\leqslant &\bm{\Delta u}_{max} \\
	\bm{u}_{min} & \leqslant &\bm{u} &\leqslant & \bm{u}_{max} \\
	\bm{\phi}_{cycle,min} & \leqslant &\bm{\phi}_{cycle} &\leqslant & \bm{\phi}_{cycle,max} \\
	\end{array}
	\label{eq_constraints_PNMPC}
\end{equation}

\begin{equation}
	\begin{array}{cccccc}
		\vspace{2mm}
		\bm{\Delta u}_{min} & = & \left[	
		\begin{array}{c}
		-2 \;\;  \textnormal{Hz} \\ 
		-1 \; \textnormal{\%} \\ 
		\end{array} \right] 
		& \bm{\Delta u}_{max} & = & \left[	
		\begin{array}{c}
		2 \;\; \textnormal{Hz} \\ 
		1 \; \textnormal{\%}  \\ 
		\end{array} \right] \\
		\vspace{2mm}
		\bm{u}_{min} & = & \left[	
		\begin{array}{c}
		30 \;\; \textnormal{Hz}  \\ 
		10 \; \textnormal{\%}  \\ 
		\end{array} \right] 
		& \bm{u}_{max} & = & \left[	
		\begin{array}{c}
		50 \;\; \textnormal{Hz}  \\ 
		100 \; \textnormal{\%}  \\ 
		\end{array} \right] \\
		\bm{\phi}_{cycle,min} & = & \left[	
		\begin{array}{c}
		0.8 \;\; \textnormal{bar}  \\ 
		14 \;\; \textnormal{bar}  \\ 
		-25 \; \textnormal{ºC}  \\
		\end{array} \right] 
		& \bm{\phi}_{cycle,max} & = & \left[	
		\begin{array}{c}
		2 \;\; \textnormal{bar}  \\ 
		25 \;\; \textnormal{bar}  \\ 
		-20 \; \textnormal{ºC}   \\
		\end{array} \right] \\
	\end{array}
	\label{eq_values_constraints_PNMPC}
\end{equation}

\begin{equation}
	\begin{aligned}
	&N_p = 10 
	&N_c = 3 \\
	&\bm Q_ =  \left[	
	\begin{array}{ccc}
	2 & 0 & 0 \\ 
	0 & 2 & 0 \\ 
	0 & 0 & 2 \\ 
	\end{array} \right] 
	&\bm R =  \left[	
	\begin{array}{cc}
	1 & 0 \\ 
	0 & 1 \\ 
	\end{array} \right] \\ 	
	\end{aligned}
	\label{eq_parameters_PNMPC}
\end{equation}	

The prediction horizon $N_p$ has been tuned considering the characteristic time constant of the system, whereas the control horizon $N_c$ has been tuned to maintain the computational cost of the whole strategy at a reasonable level, since online system linearisation must be performed at each sampling time and then QP optimization algorithms must be applied. Concerning the weighting matrices $\bm Q$ and $\bm R$, they have been tuned by  ¡initially considering the same priority for all three controlled variables (once normalised according to the intended optimal value), and a greater weight for tracking errors than for actuator wear. These tuning matrices will be modified later to highlight the different transient behaviour of the controlled variables when priority is considered.
 
Due to the system linearisation, standard solvers can be applied to perform the PNMPC optimization at every time step. In this case the MATLAB\textsuperscript{\textregistered} \emph{quadprog} tool is applied \cite{MatlabOTB}, in particular the \emph{interior-point-convex} algorithm. 
	
\subsection{Feedback-plus-feedforward control strategy} \label{FB+FF}
	
\begin{figure}[h] 
	\centering
	\centerline{\includegraphics[width=15cm,trim = 125 550 105 120,clip]{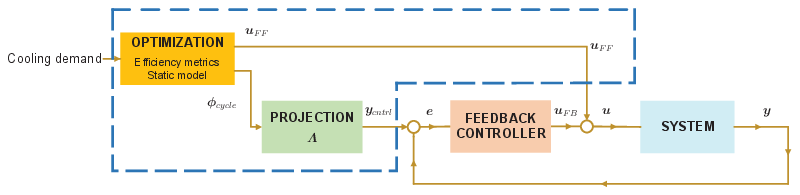}}
	\caption{Feedback-plus-feedforward controller by Jain \cite{jain2013thermodynamics}} 
	\label{fig_FFplusFB}
\end{figure}

The feedback-plus-feedforward (FB+FF) controller developed in \cite{jain2013thermodynamics} is described in Fig. \ref{fig_FFplusFB}. This controller benefits from the fact that not only does the optimizer provide the key variables of the optimal cycle, but also the steady-state values of the manipulated variables corresponding to the optimal cycle, as indicated in Table \ref{tab_optimal_cycle}. These steady-state values are suitable to be used in a feedforward block, which \emph{linearises} the system around the desired point and allows the feedback contribution to the control actions to effectively drive the system to the optimal desired point. 

However, as stated in Section \ref{Controllability}, the reduced controllability involves that the optimal point defined in the three-dof optimization space cannot be achieved in general. A solution for this problem is provided by Jain \cite{jain2013thermodynamics} by projecting the optimal point onto the two-dof control space by means of matrix $\bm \Lambda \in \mathbb{R}^{2x3}$ in Eq. \eqref{eq_Projection_y}.  

\begin{equation}
	\begin{aligned}
	\bm y_{cntrl} = \bm \Lambda \; \bm \phi_{cycle} \\
	\end{aligned}
	\label{eq_Projection_y}
\end{equation}

According to Eq. \eqref{eq_Projection_y}, the projection matrix $\bm \Lambda \in \mathbb{R}^{2x3}$ defines the linear combination between the variables included in $\bm \phi_{cycle} \in \mathbb{R}^3$ which generates the controlled variable vector $\bm y_{cntrl} \in \mathbb{R}^2$. In order to fairly compare this control strategy to the PNMPC, the matrix $\bm \Lambda$ presented in Eq. \eqref{eq_Projection_y_LAMBA} is selected. This forces the FB+FF controller to at least provide the same cooling power as the PNMPC strategy. It also involves that, with this projection matrix, the FB+FF controller gives up controlling the condenser pressure, whose steady-state value will be determined by the two-dof control space. The evaporator pressure has been selected as one of the controlled variables instead of the condenser pressure according to the estimated pressure drop along each heat exchanger. As indicated in Section \ref{Optimization}, a brazed-plate evaporator and a cross-flow tube condenser have been considered. Pressure drop is known to be greater in tube exchangers, therefore the measurement of the evaporator pressure might be more reliable.

\begin{equation}
	\begin{aligned}
	\bm y_{cntrl} & = 
	\left[	
	\begin{array}{c}
	y_{cntrl,1} \\ 
	y_{cntrl,2} \\
	\end{array} 
	\right] = 
	\bm \Lambda \; \bm \phi_{cycle} = 
	\left[
	\begin{array}{ccc}
	1 & 0 & 0 \\ 
	0 & 0 & 1 \\ 
	\end{array} 
	\right] \; 
	\left[
	\begin{array}{c}
	P_e \\ 
	P_c \\
	T_{e,sec,out} \\
	\end{array} 
	\right] = 
	\left[	
	\begin{array}{c}
	P_e \\ 
	T_{e,sec,out} \\
	\end{array} 
	\right]
	\end{aligned}
	\label{eq_Projection_y_LAMBA}
\end{equation}

After the application of the feedforward term, a two-input, two-output linear state-space model is identified, for which a linear quadratic regulator (LQR) is applied to provide the feedback contribution to the final control action. Since the identified model is type 0, the state vector is augmented with two tracking error integrators. The resulting gain matrix $\bm K_{FB+FF}$ is shown in Eq. \eqref{eq_K_FBplusFF}. Although it was not included in the original formulation \cite{jain2013thermodynamics}, an anti-windup strategy has been also included to deal with manipulated input saturation and the tracking error integrators.

\begin{equation}
	\begin{aligned}
	\bm K_{FB+FF} = 
	\left[
		\begin{array}{cccccc}
			-8.73 \cdot 10^{-5} & -1.1493 & -7.08 \cdot 10^{-6} & -0.0562 \\
			-1.01 \cdot 10^{-5} & -4.3996 & -8.18 \cdot 10^{-6} & -0.6346 \\
		\end{array}
	\right]
	\end{aligned}
	\label{eq_K_FBplusFF}
\end{equation}

\subsection{Simulation results} \label{Simulations}

Two random initial cycles IP\textsubscript{1} and IP\textsubscript{2} have been considered, featuring the variables gathered in Table \ref{tab_initial_conditions}. The aim of this section is to illustrate the controllability issues highlighted by the non-linear study performed in Section \ref{Controllability}. Starting at the described initial cycles, the control objective is to achieve the optimal operating point described in Table \ref{tab_optimal_cycle}. All simulations presented in this subsection are obtained using the reduced-order dynamic model presented in Section \ref{Modelling}.

\begin{figure*}[h]
	\centering
	\subfigure[Condenser and evaporator pressures]{
		\includegraphics[width=6.4cm,angle=0] {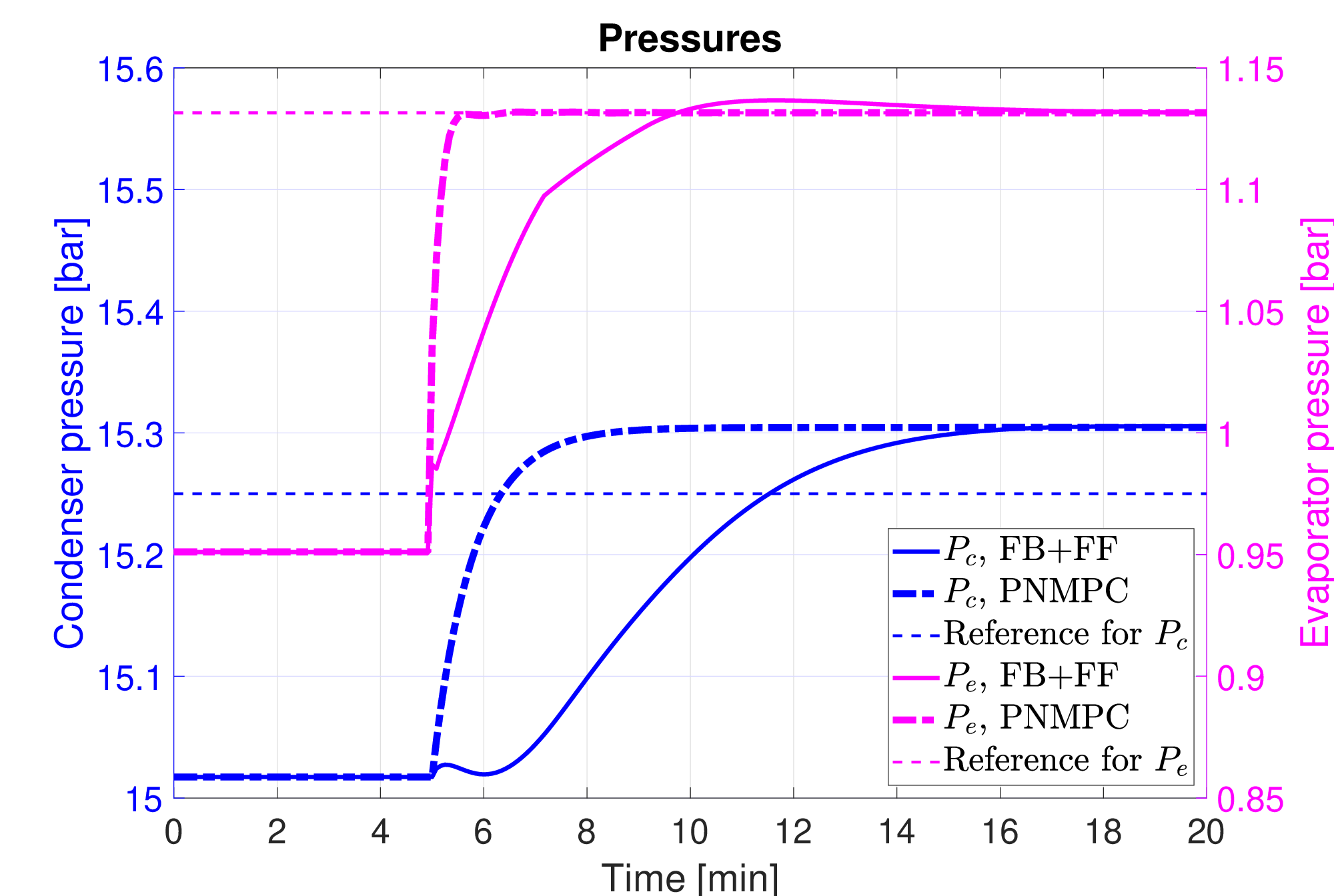}
		\label{fig_control_600W_IP1_a}
	}\subfigure[Outlet temperature of the evaporator secondary fluid]{
		\includegraphics[width=6.4cm,angle=0] {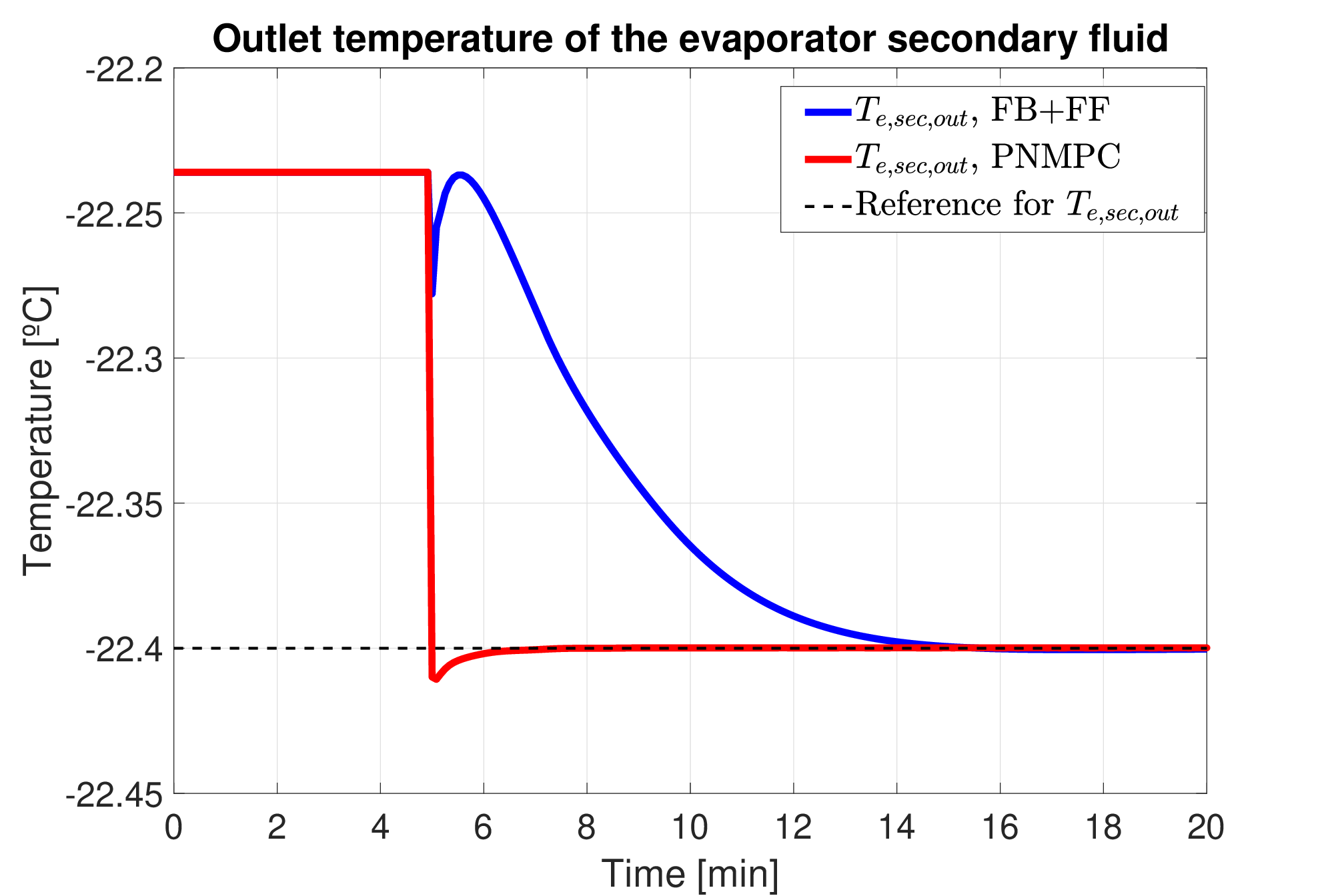}
		\label{fig_control_600W_IP1_b}
	}
	\subfigure[Degree of superheating and refrigerant mass flow]{
		\includegraphics[width=6.4cm,angle=0] {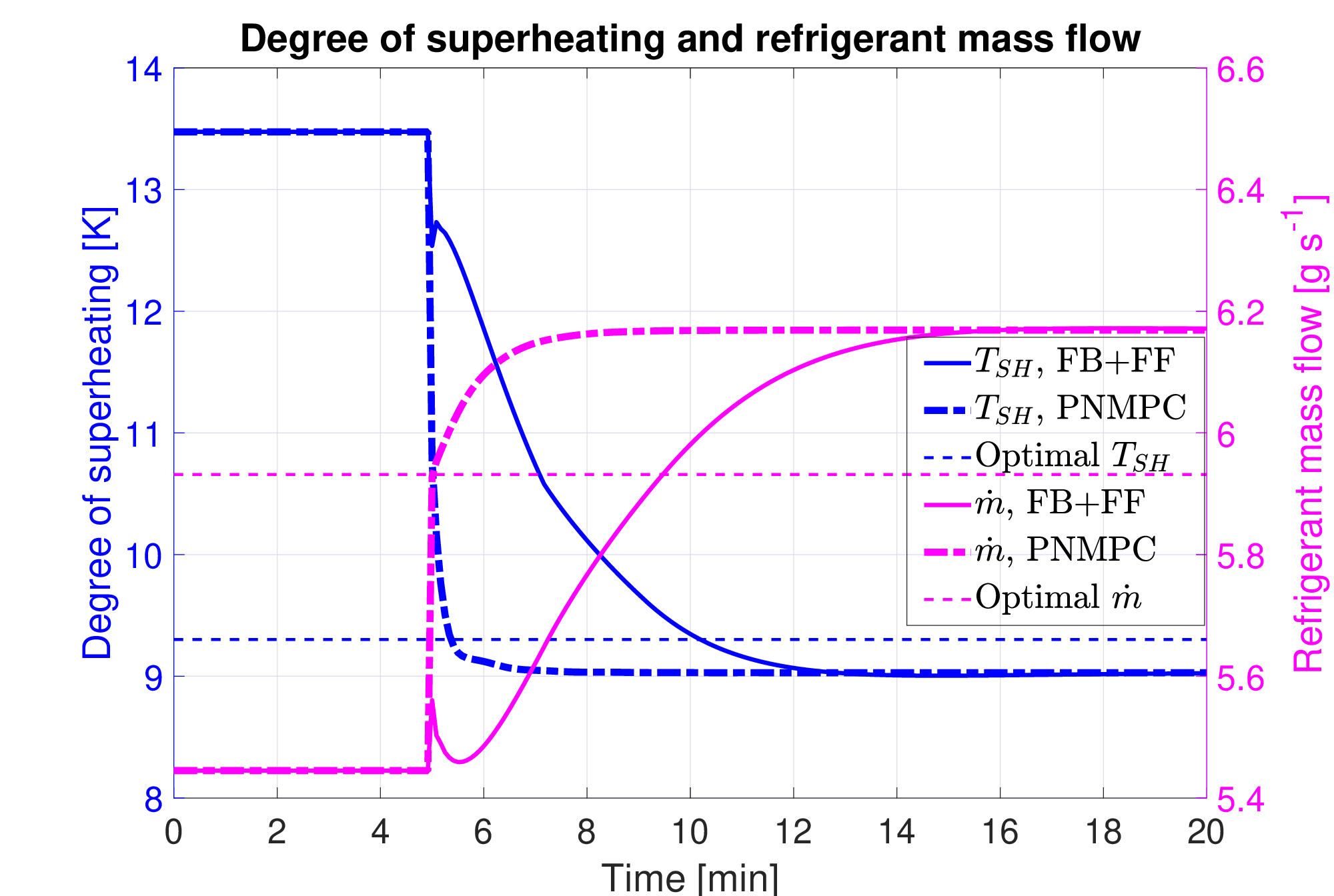}
		\label{fig_control_600W_IP1_c}
	}\subfigure[Manipulated variables]{
		\includegraphics[width=6.4cm,angle=0] {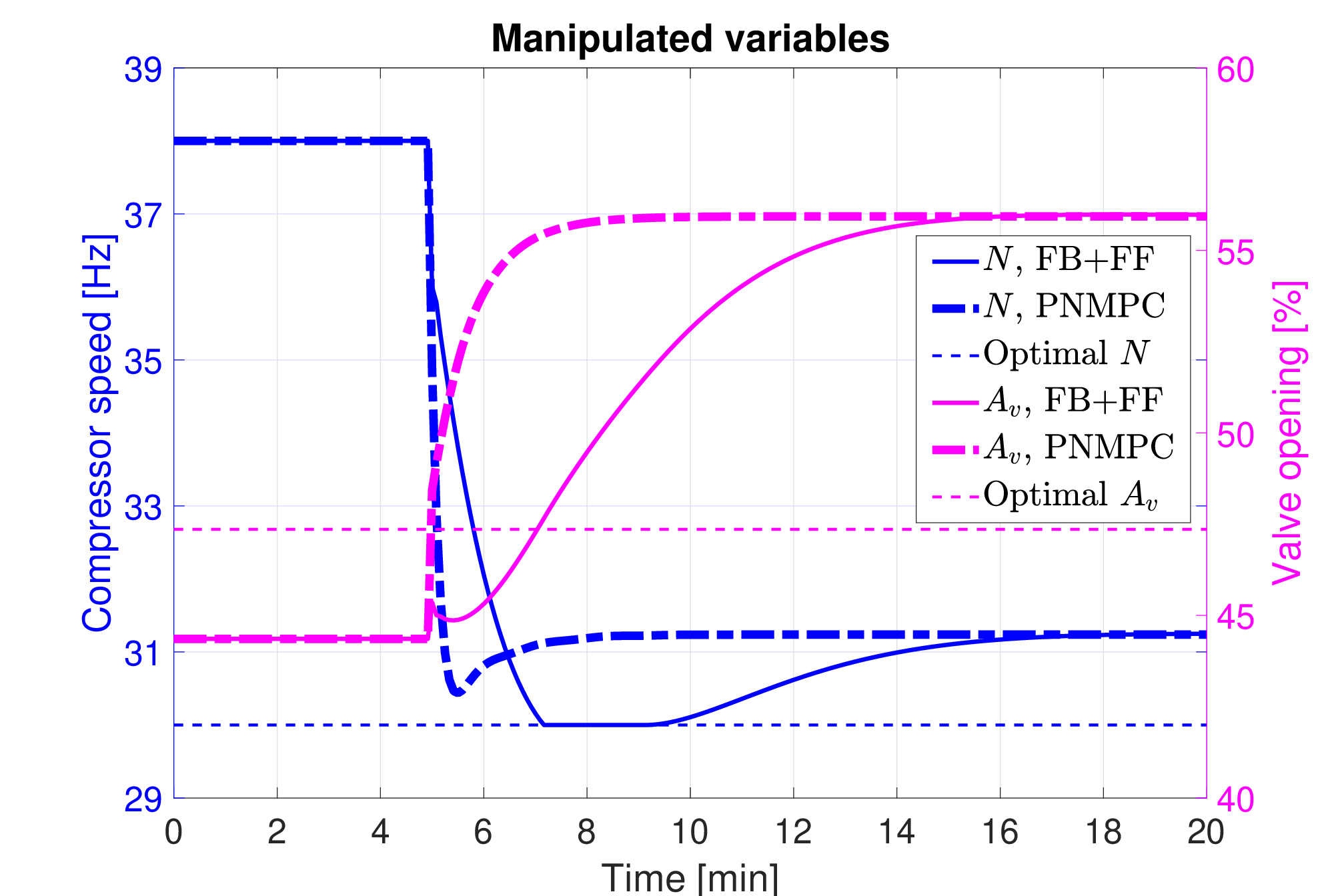}
		\label{fig_control_600W_IP1_d}
	}
	\caption{Control performance when trying to achieve optimal operation from initial point IP\textsubscript{1} using PNMPC and FB+FF control}
	\label{fig_control_600W_IP1}
\end{figure*}

\begin{figure*}[h]
	\centering
	\subfigure[Condenser and evaporator pressures]{
		\includegraphics[width=6.4cm,angle=0] {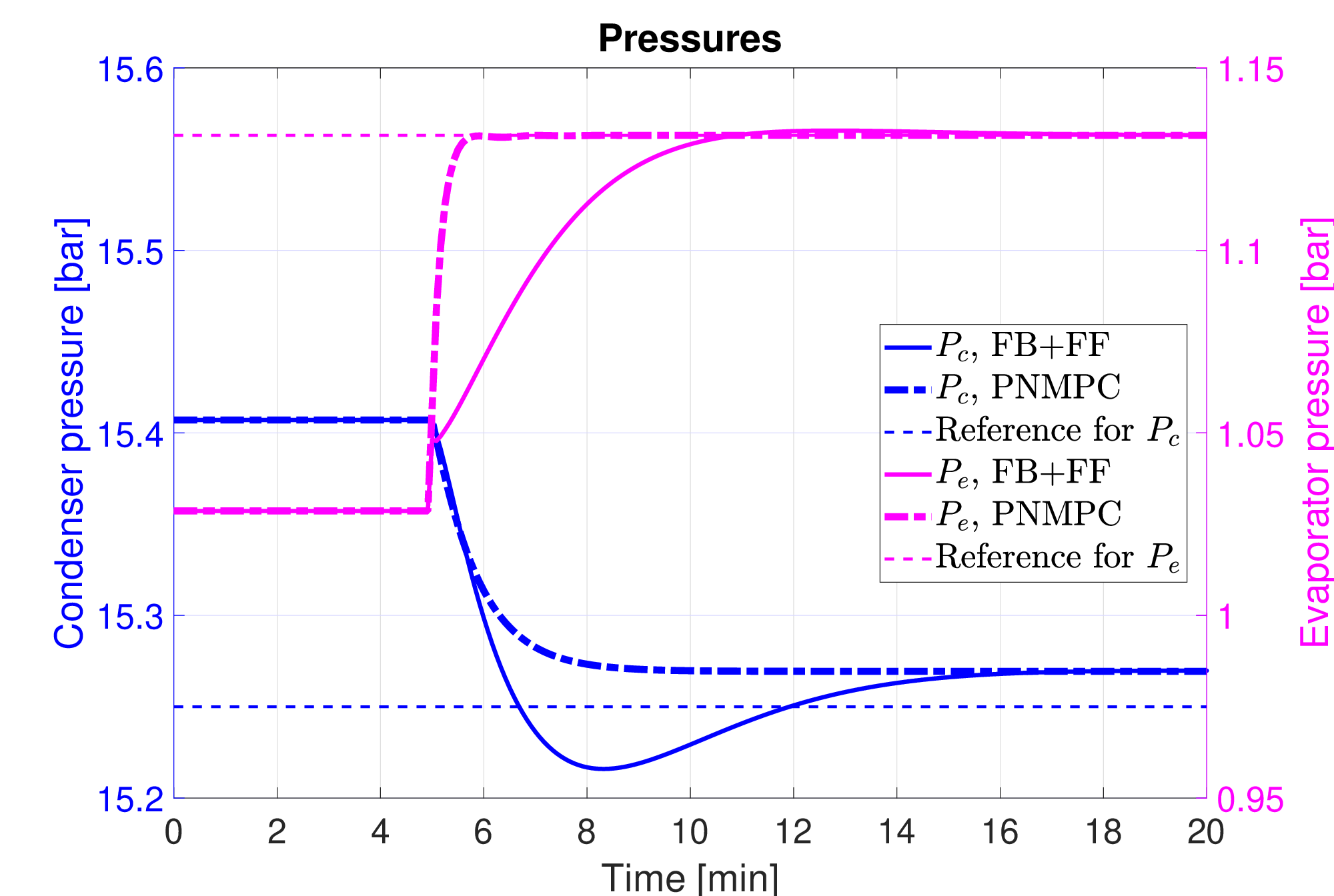}
		\label{fig_control_600W_IP2_a}
	}\subfigure[Outlet temperature of the evaporator secondary fluid]{
		\includegraphics[width=6.4cm,angle=0] {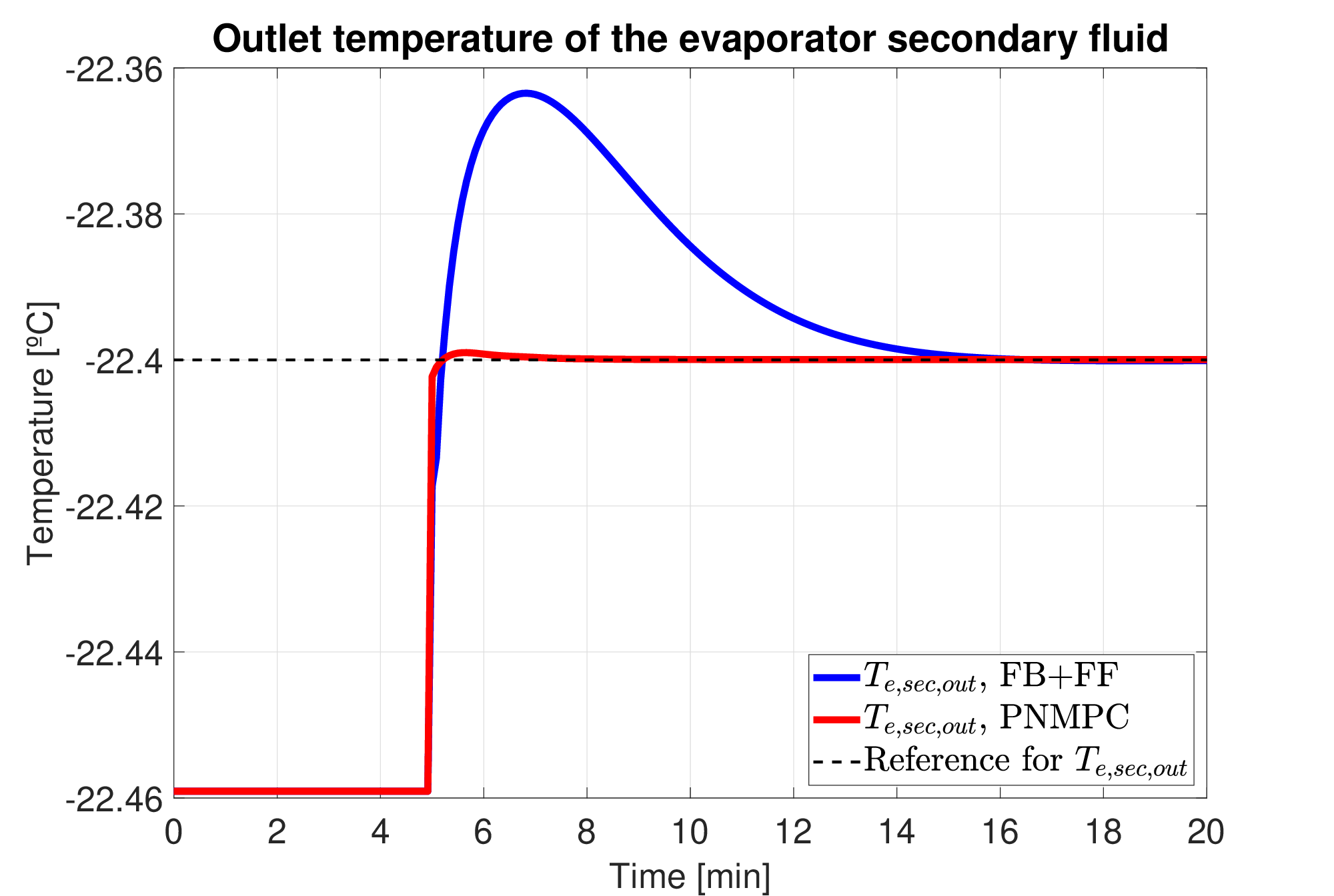}
		\label{fig_control_600W_IP2_b}
	}
	\subfigure[Degree of superheating and refrigerant mass flow]{
		\includegraphics[width=6.4cm,angle=0] {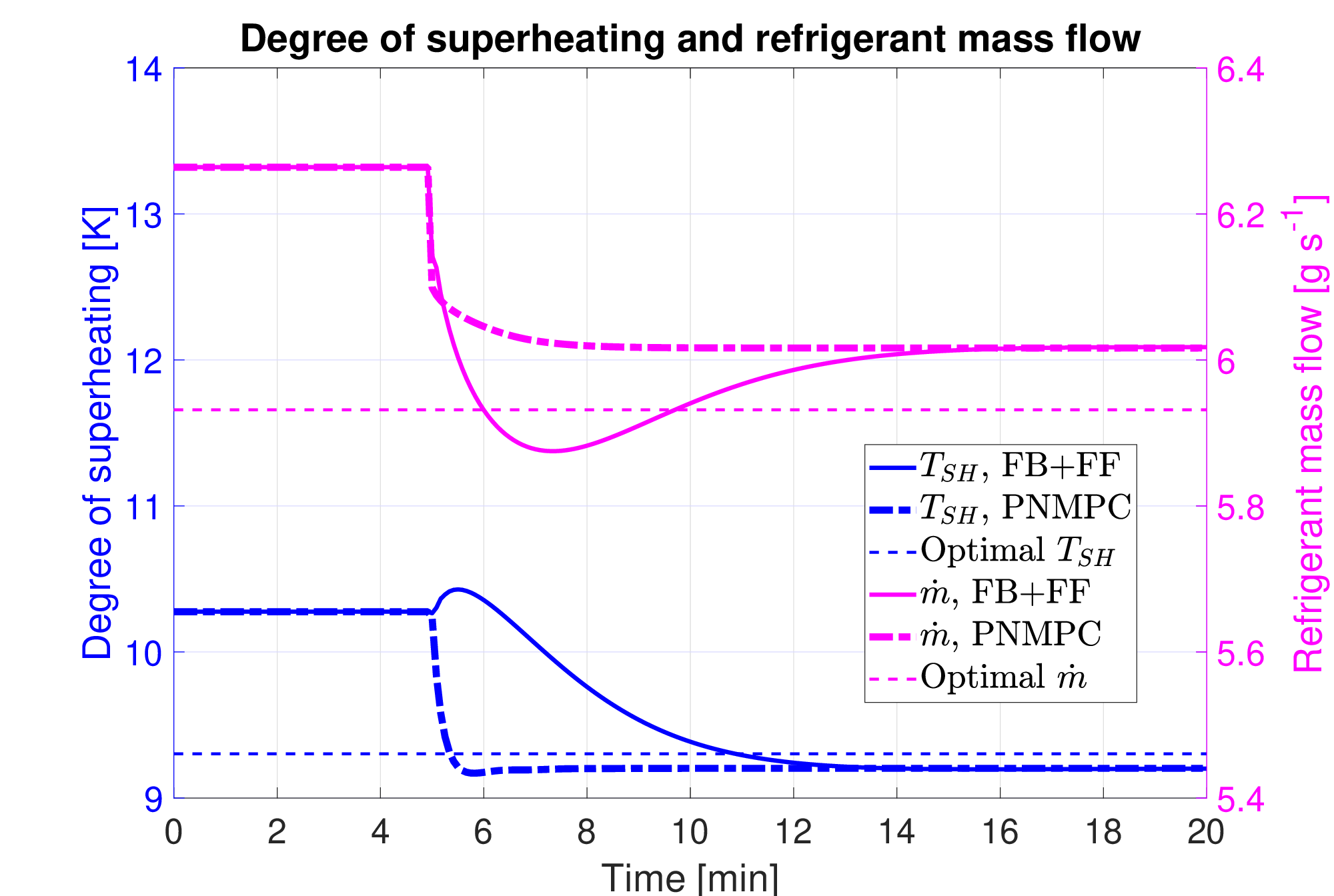}
		\label{fig_control_600W_IP2_c}
	}\subfigure[Manipulated variables]{
		\includegraphics[width=6.4cm,angle=0] {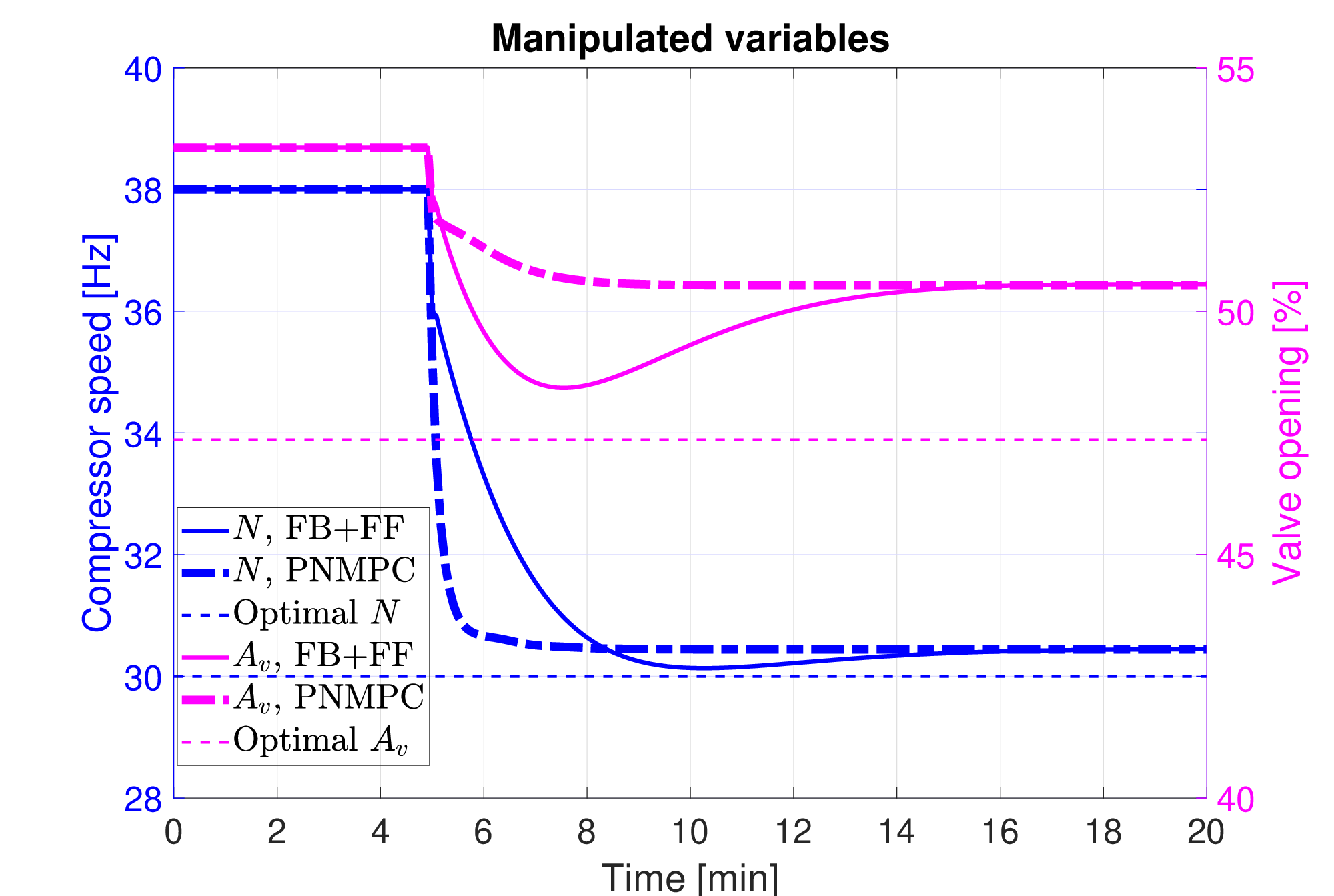}
		\label{fig_control_600W_IP2_d}
	}
	\caption{Control performance when trying to achieve optimal operation from initial point IP\textsubscript{2} using PNMPC and FB+FF control}
	\label{fig_control_600W_IP2}	
\end{figure*}

The performance of both controllers is analysed in Fig. \ref{fig_control_600W_IP1} and \ref{fig_control_600W_IP2} for both initial cycles. The behaviour of the pressures $P_e$ and $P_c$ (subplots \ref{fig_control_600W_IP1_a} and \ref{fig_control_600W_IP2_a}), the outlet temperature of the evaporator secondary fluid $T_{e,sec,out}$ (subplots \ref{fig_control_600W_IP1_b} and \ref{fig_control_600W_IP2_b}), the degree of superheating $T_{SH}$ and the refrigerant flow $\dot{m}$ (subplots \ref{fig_control_600W_IP1_c} and \ref{fig_control_600W_IP2_c}), and finally the control inputs $N$ and $A_v$ (subplots \ref{fig_control_600W_IP1_d} and \ref{fig_control_600W_IP2_d}) are compared for both control strategies. Open-loop behaviour is represented before $t$ = 5 min, when the control loop is closed. Although the original formulation did not include constraints, the maximum and minimum values of the manipulated variables, as well as the rate limits detailed in Eq. \eqref{eq_values_constraints_PNMPC}, are also applied to the FB+FF controller.

As shown in Fig. \ref{fig_control_600W_IP1} and \ref{fig_control_600W_IP2}, both control strategies reach the operating point, which involves a $COP$ of 1.1944 when starting at IP\textsubscript{1}, and a $COP$ of 1.2266 when applying control from IP\textsubscript{2}, far away from the optimum indicated in Table \ref{tab_optimal_cycle}. However, the PNMPC strategy shows better dynamic performance regarding settling time and rise time. The PNMPC reaches the equilibrium much faster than the FB+FF control strategy, which implies that the suboptimal cycle is achieved earlier and it can bring a certain improvement in energy consumption, at least in transient. The enhanced performance of the PNMPC controller against the FB+FF strategy is due to several reasons. Firstly, the PNMPC technique does not require a linearised model around the optimum, since it is not based on the equilibrium point concept; therefore the controller is more likely to work well around points far away from the optimum, like IP\textsubscript{1} and IP\textsubscript{2}. Secondly, as optimal cycles involve $N$ = 30 Hz for a wide range of the required cooling power, as Fig. \ref{fig_optimizacion} depicts, the PNMPC is able to better deal with proximity to control input and rate limits.

Regarding steady-state performance, it is observed in Fig. \ref{fig_control_600W_IP1} and \ref{fig_control_600W_IP2} how the achieved values of $P_e$ and $T_{e,sec,out}$  match those corresponding to the optimal cycle, but $P_c$ does not match the optimal reference, neither starting at IP\textsubscript{1} nor at IP\textsubscript{2}, although the value achieved starting at IP\textsubscript{2} is closer to the optimal one; that is the reason why the $COP$ achieved in this case is greater and therefore closer to the global optimum. Notice also that neither $\dot{m}$, nor $T_{SH}$, nor the control inputs achieve the steady-state optimal values, therefore the achieved points are clearly suboptimal. The steady-state values of the control inputs show to be greater than those provided by the optimizer, which means that the cooling demand is satisfied with higher refrigerant mass flow and compressor speed than necessary, that in turn involves that the $COP$ is reduced. 

\begin{figure*}[h]
	\centering
	\subfigure[Condenser and evaporator pressures]{
		\includegraphics[width=6.4cm,angle=0] {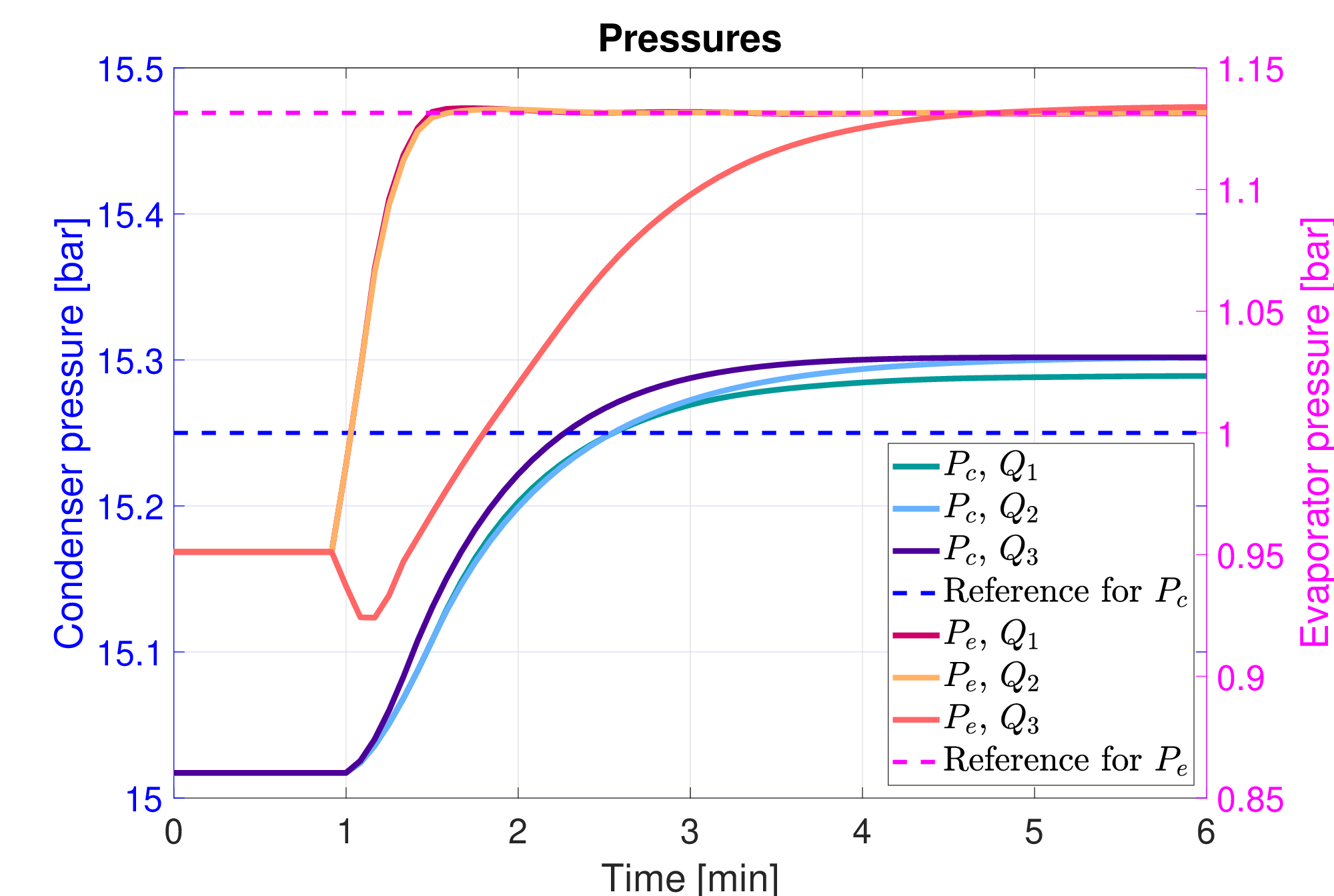}
		\label{fig_control_600W_IP1_Qi_a}
	}\subfigure[Outlet temperature of the evaporator secondary fluid]{
		\includegraphics[width=6.4cm,angle=0] {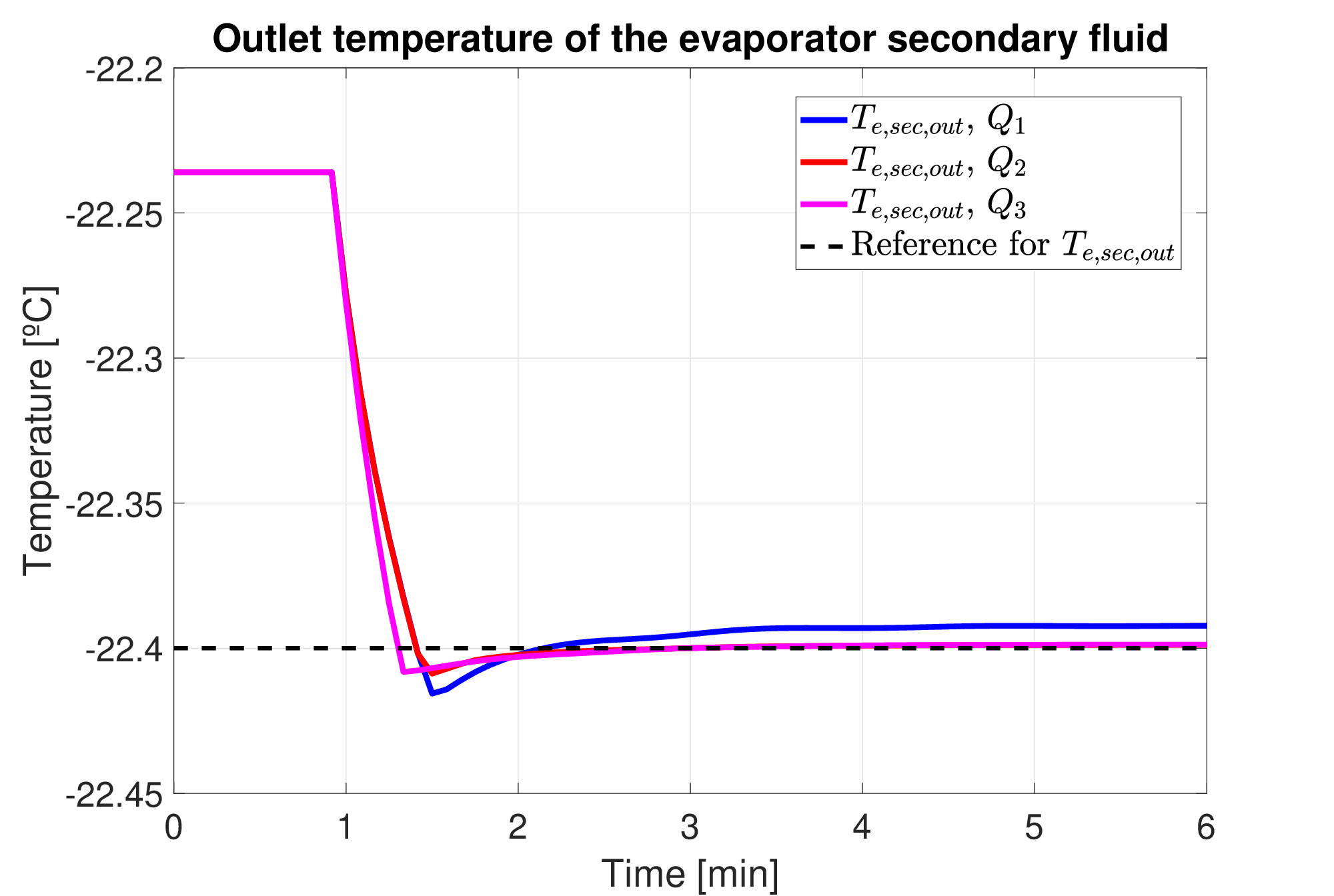}
		\label{fig_control_600W_IP1_Qi_b}
	}
	\subfigure[Degree of superheating and refrigerant mass flow]{
		\includegraphics[width=6.4cm,angle=0] {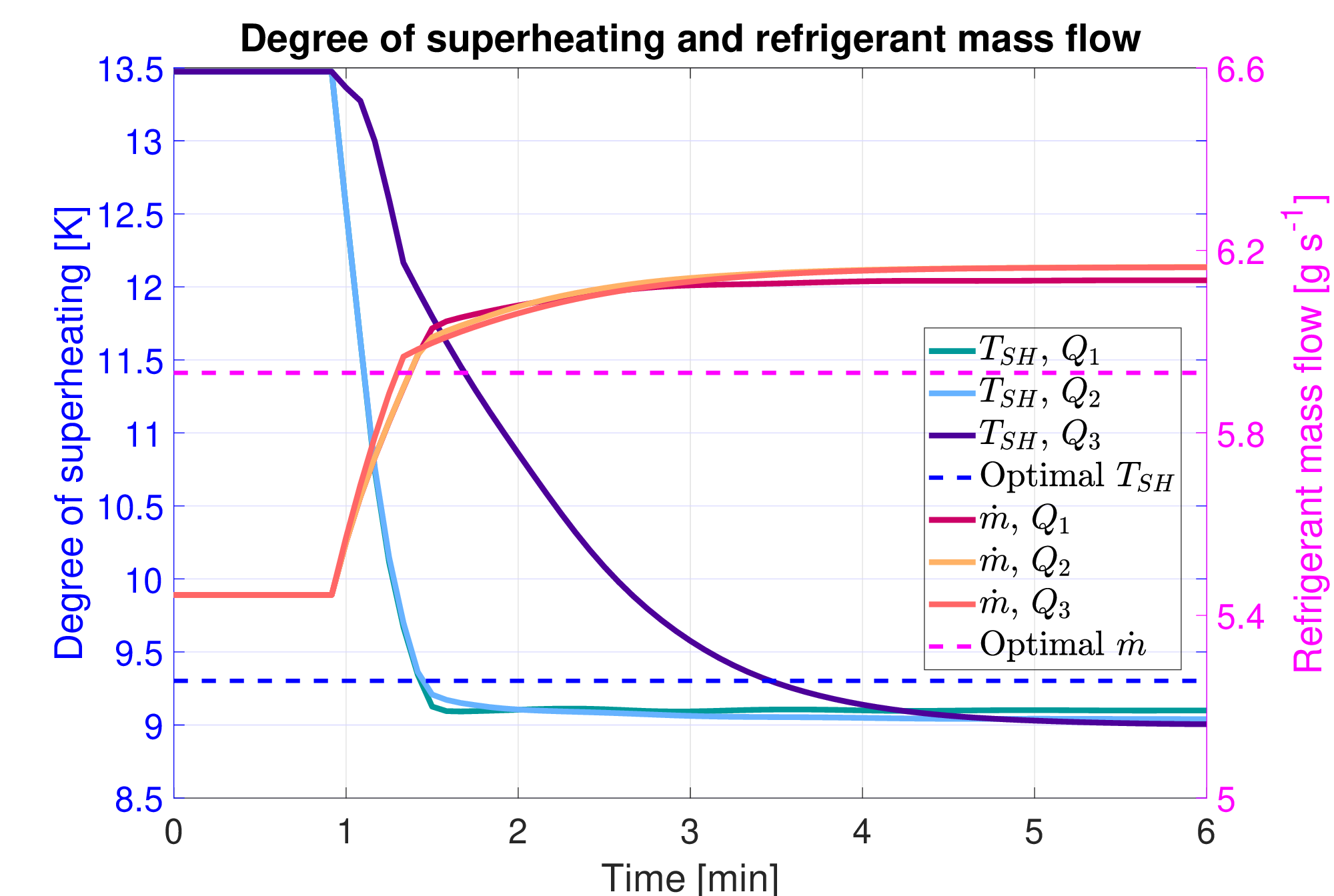}
		\label{fig_control_600W_IP1_Qi_c}
	}\subfigure[Manipulated variables]{
		\includegraphics[width=6.4cm,angle=0] {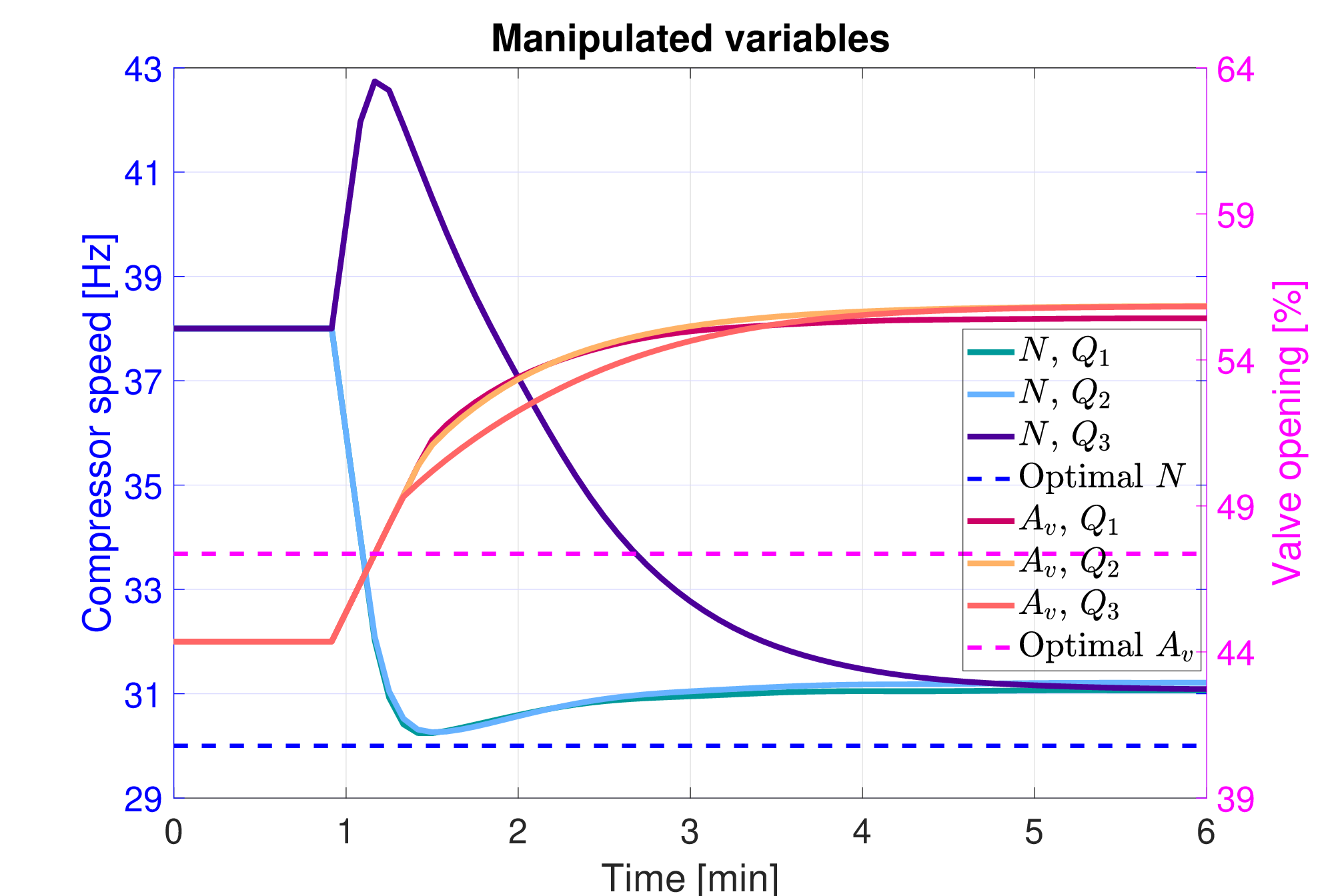}
		\label{fig_control_600W_IP1_Qi_d}
	}
	\caption{Control performance when trying to achieve optimal operation from initial point IP\textsubscript{1} and applying different PNMPC tuning}
	\label{fig_control_600W_IP1_Qi} 
\end{figure*}

\begin{figure*}[h]
	\centering
	\subfigure[Condenser and evaporator pressures]{
		\includegraphics[width=6.4cm,angle=0] {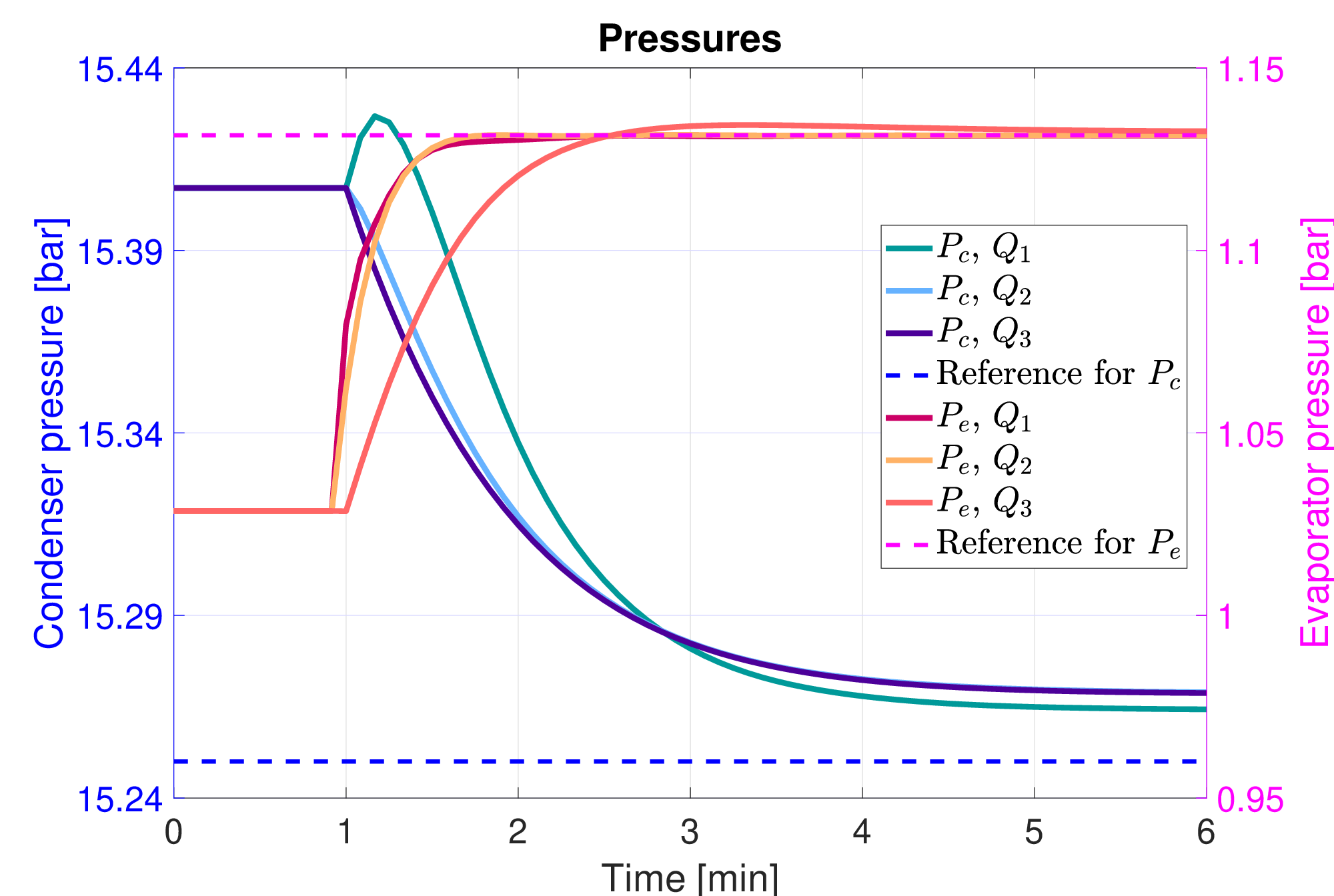}
		\label{fig_control_600W_IP2_Qi_a}
	}\subfigure[Outlet temperature of the evaporator secondary fluid]{
		\includegraphics[width=6.4cm,angle=0] {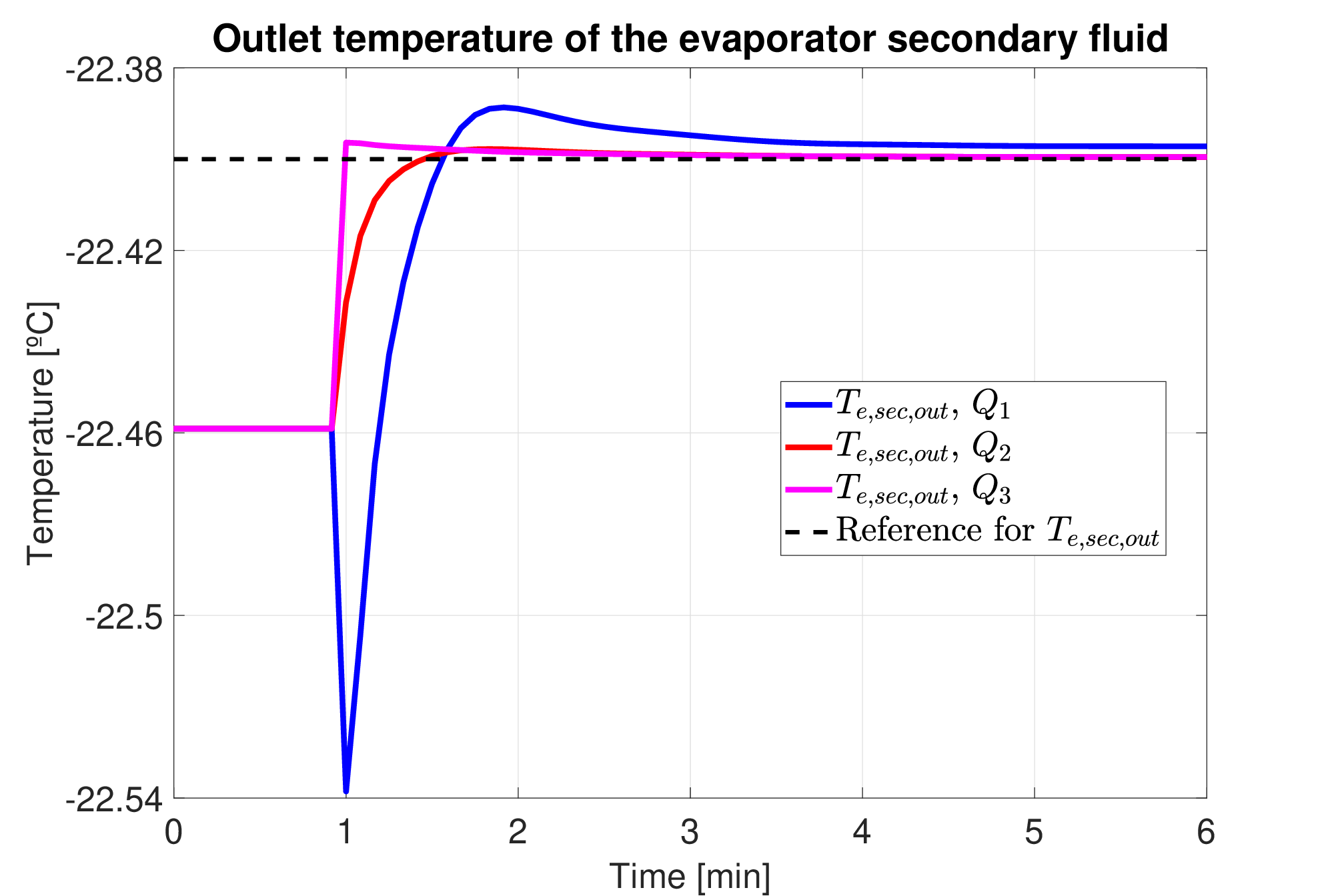}
		\label{fig_control_600W_IP2_Qi_b}
	}
	\subfigure[Degree of superheating and refrigerant mass flow]{
		\includegraphics[width=6.4cm,angle=0] {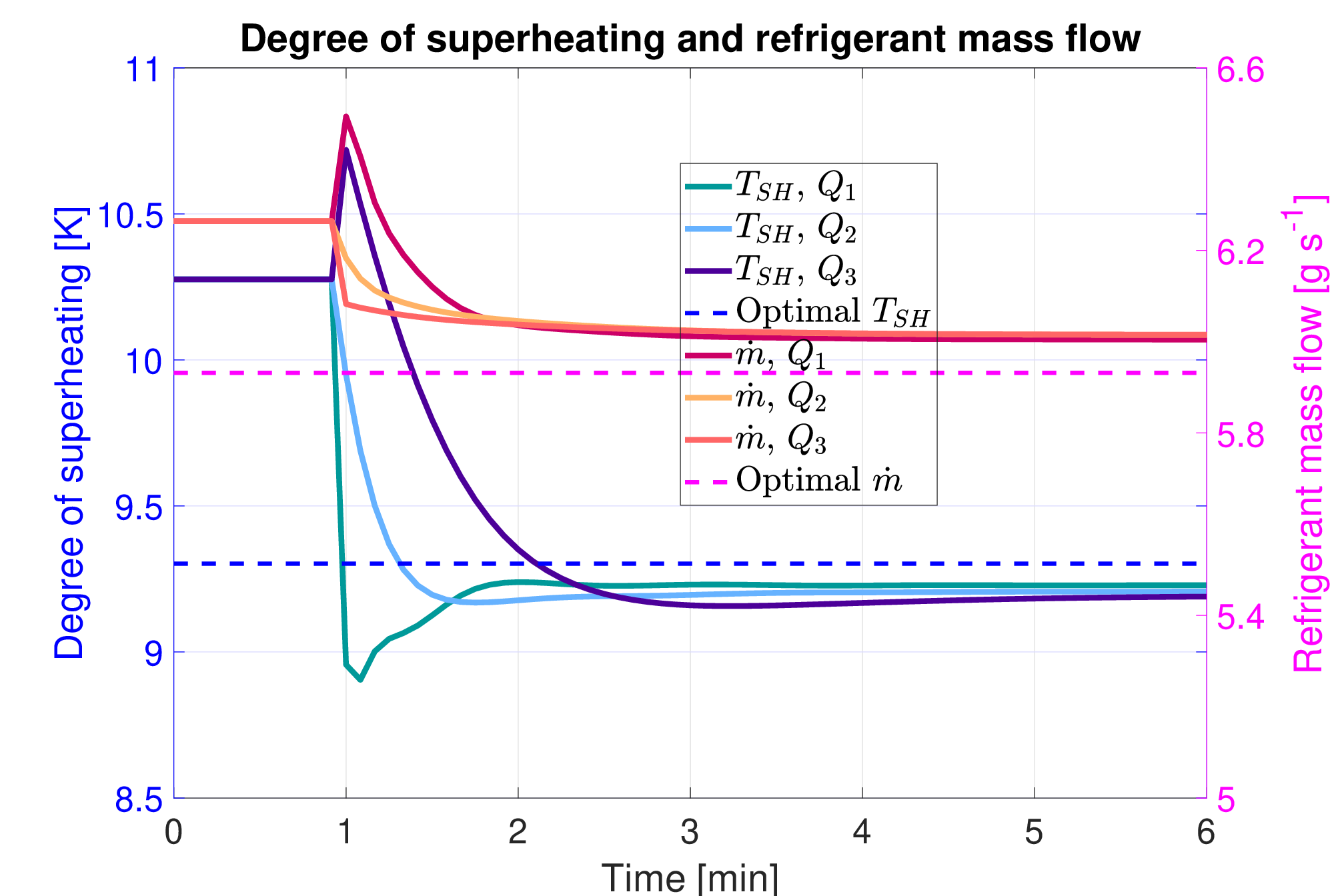}
		\label{fig_control_600W_IP2_Qi_c}
	}\subfigure[Manipulated variables]{
		\includegraphics[width=6.4cm,angle=0] {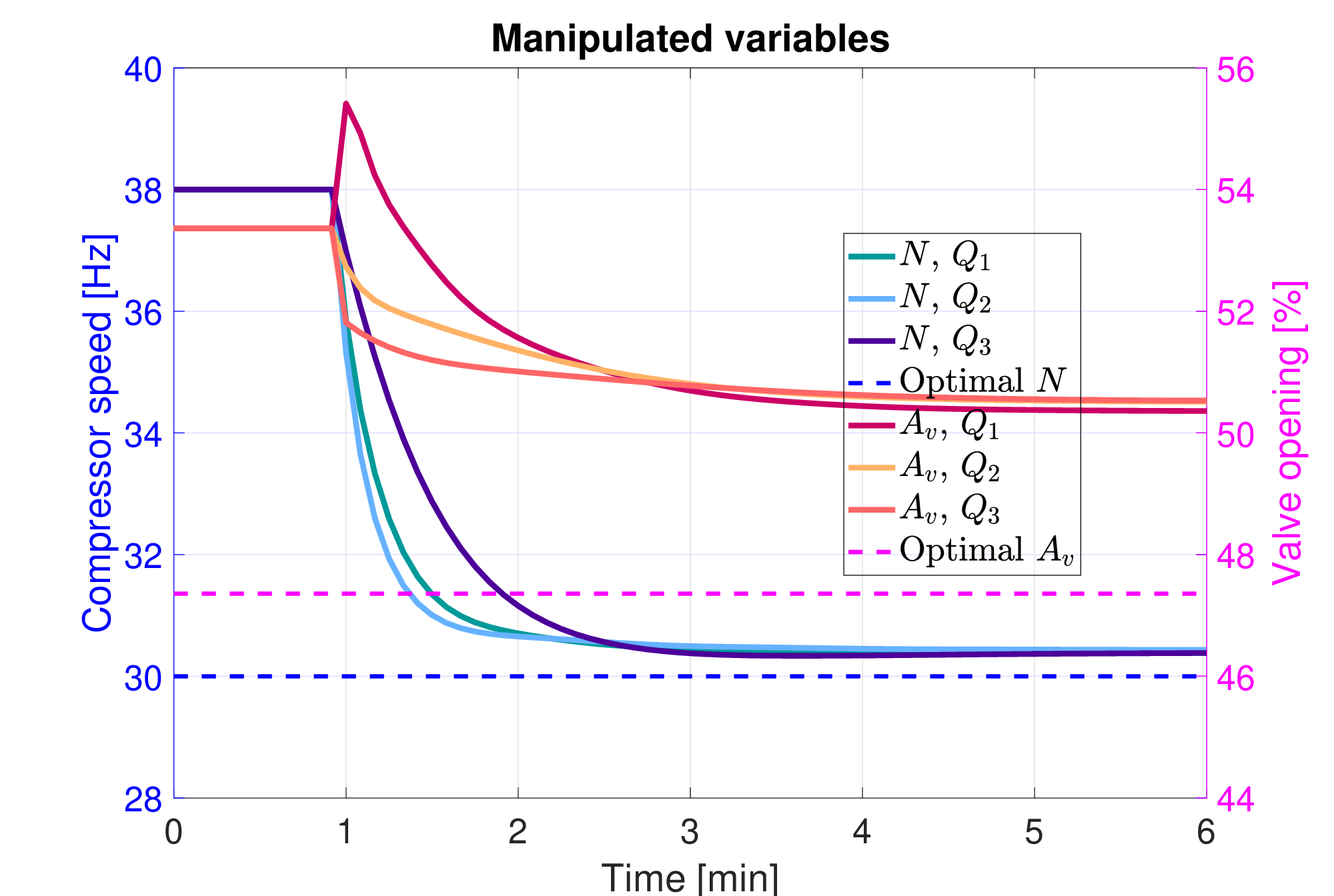}
		\label{fig_control_600W_IP2_Qi_d}
	}
	\caption{Control performance when trying to achieve optimal operation from initial point IP\textsubscript{2} and applying different PNMPC tuning}
	\label{fig_control_600W_IP2_Qi} 
\end{figure*}

Some tests on the PNMPC tuning are also performed. In Eq. \eqref{eq_Qi} three different error weighting matrices $\bm Q$ are gathered. They have been chosen to highlight the different transient behaviour of the controlled variables when one specific output is set as priority, in this case $T_{e,sec,out}$, i.e. the measurement of the cooling power supplied by the cycle. The remaining tuning parameters (prediction and control horizons, weighting matrix $\bm R$, and sampling time) are held constant and match the values indicated in Section \ref{PNMPC}. The control results are depicted when applying the different PNMPC tuning from IP\textsubscript{1} and IP\textsubscript{2}, respectively, in Fig. \ref{fig_control_600W_IP1_Qi} and \ref{fig_control_600W_IP2_Qi}. Notice that cool-coloured lines refer to the blue left X-axis, while warm-coloured lines refer to the magenta right Y-axis in both two-Y-axis charts.

\begin{equation}
	\begin{aligned}
	 \bm Q_{1} & =  \left[	
	 \begin{array}{ccc}
	 2 & 0 & 0 \\ 
	 0 & 2 & 0 \\ 
	 0 & 0 & 2 \\ 
	 \end{array} \right] \;\;\;\;
	 \bm Q_{2} & =  \left[	
	 \begin{array}{ccc}
	 2 & 0 &  0 \\ 
	 0 & 2 &  0 \\ 
	 0 & 0 & 20 \\ 
	 \end{array} \right] \;\;\;\;
	 \bm Q_{3} & =  \left[	
	 \begin{array}{ccc}
	 0.2 &  0 &   0 \\ 
	 0   & 20 &   0 \\ 
	 0   &  0 & 200 \\ 
	 \end{array} \right] \;\;\;\;	
	 \end{aligned}
	 \label{eq_Qi}
\end{equation}

As observed in Fig. \ref{fig_control_600W_IP1_Qi} and \ref{fig_control_600W_IP2_Qi}, the achieved cycle in steady state seems not to depend to a great extent on the controller tuning, but much more on the initial conditions. Little differences in the output values achieved in steady state are due to weighting matrices $\bm Q$, which lead the optimal cost function $J_{PNMPC}$ to different minimum values, but the achieved $COP$ mainly depends on the starting point, as shown in Table \ref{tab_COP}. However, the transient behaviour of the controlled variables behaves as expected, according to the weights imposed on the different matrices $\bm Q$ shown in Eq. \eqref{eq_Qi}. For instance, the settling time of $T_{e,sec,out}$ is much smaller when a greater tracking error weight is imposed using $\bm Q_{3}$.


\section{Conclusions} \label{Conclusions}

This work has provided a number of contributions to address optimal control of refrigeration cycles. The dynamic modelling of the heat exchangers has been assessed using the well-known \emph{moving boundary} approach and, considering (1) the frequency features of the condenser and evaporator intrinsic dynamics and (2) the fastest state variables among the condenser internal states, a third-order state-space model has been proposed for the complete refrigeration cycle. Simulation results have shown that the reduced-order model is able to provide reasonable accuracy when compared to the original high-order model, while the computational load is much reduced.

Moreover, a global optimization to provide references for the control stage has been proposed using the $COP$ as energy efficiency metrics. The canonical cycle is shown to be unequivocally defined by three independent variables, in the case that only the expansion valve opening and the compressor speed are considered as manipulated inputs and the secondary flows at both heat exchangers are seen as measurable disturbances. A highly coupled non-linear static model is included within the optimization to define the feasibility of the decision variables and ease the introduction of relevant technological and operating requirements. Some simulations have been presented comparing the different optimal cycles for a wide range of cooling demand. As a result, the conventional strategy consisting in holding a constant low value of the degree of superheating for the whole range of cooling demand has proven not to be always optimal, and it would be desirable to apply a control strategy able to move the system to the optimal solution provided by the optimizer.

Since a canonical cycle may be completely defined by three independent variables, these could be selected to match the three state variables of the reduced-order model previously proposed. Therefore, the controllability of the whole cycle is studied by carrying out a study on the controllability of the reduced-order model of the condenser. In the light of the results, the system turns out not to be completely controllable. It involves that, from a random initial cycle, the optimal solution is not reachable by exclusively manipulating the compressor and the expansion valve. 

This issue has been illustrated by applying the PNMPC strategy to drive the whole cycle to the optimal operating point for a certain cooling demand, computed using the proposed optimizer. Furthermore, the PNMPC has been compared with a feedback-plus-feedforward controller. It has been shown how the optimal cycle is not reached by any of the applied controllers, whereas the dynamic behaviour of the PNMPC shows to be faster and less sensitive to operating points. The achieved suboptimal cycles seem not to depend to a great extent on the controller tuning, but most on the initial point, which matches the conclusions of the controllability study.

As future work, some control issues such as recursive feasibility, stability, and robustness of the PNMPC strategy are intended to be studied in depth. Moreover, given the problems arising when trying to achieve the optimal solution, a suboptimal control strategy is intended to be developed to compute, among the achievable solutions, that optimizing a multi-objective cost function, where not only the steady-state efficiency but also the transient energy consumption are considered.


\section*{Acknowledgement}

The authors would like to acknowledge Spanish MCeI (Grant RTI2018-101897-B-I00) and Brazilian CNPq (project number 305785/2015-0) for funding this work.




\pagebreak

\section*{List of table titles}

\noindent Table I: Optimal cycle variables for a given cooling demand \\

\noindent Table II: Cycle variables of IP\textsubscript{1} and IP\textsubscript{2} \\

\noindent Table III: $COP$ achieved for initial points IP\textsubscript{1} and IP\textsubscript{2} and the different PNMPC tuning \\

\pagebreak

\begin{table} [h]
	\centering
	\caption{Optimal cycle variables for a given cooling demand}
	\label{tab_optimal_cycle}
	\scalebox{0.7}[0.7]{
		\begin{tabular}{c c c} 
			\midrule
			\textbf{Variable} & \textbf{Unit} & \textbf{Value} \\
			\midrule
			{\color{blue} $h_{e,out}$} & {\color{blue} kJ kg\textsuperscript{-1}} & {\color{blue} $349.01$} \\
			{\color{blue} $h_{c,out}$} & {\color{blue} kJ kg\textsuperscript{-1}} & {\color{blue} $248.10$} \\
			{\color{blue} $\dot{m}$} & {\color{blue} g s\textsuperscript{-1}} & {\color{blue} $5.93$} \\
			\midrule
			{\color{green} $P_e$} & {\color{green} bar} & {\color{green} $1.132$} \\
			{\color{green} $P_c$} & {\color{green} bar} & {\color{green} $15.25$} \\
			{\color{green} $T_{SH}$} & {\color{green} K} & {\color{green} $9.3$} \\
			\midrule
			{\color{magenta} $N$} & {\color{magenta} Hz} & {\color{magenta} $30$} \\
			{\color{magenta} $A_v$} & {\color{magenta} \%} & {\color{magenta} $47.36$} \\
			\midrule
			{\color{brown} $COP$} & {\color{brown} —} & {\color{brown} $1.2453$} \\
			\midrule
		\end{tabular}}
\end{table}

\pagebreak

\begin{table} [H]
	\centering
	\caption{Cycle variables of IP\textsubscript{1} and IP\textsubscript{2}}
	\label{tab_initial_conditions}
	\scalebox{0.7}[0.7]{
		\begin{tabular}{c c c c} 
			\midrule
			\textbf{Variable} & \textbf{Unit} & \parbox[c][1cm][c]{2.0cm}{\centering \textbf{\textbf{IP\textsubscript{1}}}} & \parbox[c][1cm][c]{2.0cm}{\centering  \textbf{\textbf{IP\textsubscript{2}}}} \\
			\midrule
			{\color{black} $h_{e,out}$} & {\color{black} kJ kg\textsuperscript{-1}} & {\color{black} $350.01$} & {\color{black} $348.51$}\\
			{\color{black} $h_{c,out}$} & {\color{black} kJ kg\textsuperscript{-1}} & {\color{black} $247.59$} & {\color{black} $250.60$} \\
			{\color{black} $\dot{m}$} & {\color{black} g s\textsuperscript{-1}} & {\color{black} $5.45$} & {\color{black} $6.26$} \\
			\midrule
			{\color{black} $P_e$} & {\color{black} bar} & {\color{black} $0.951$} & {\color{black} $1.029$} \\
			{\color{black} $P_c$} & {\color{black} bar} & {\color{black} $15.02$} & {\color{black} $15.41$} \\
			{\color{black} $T_{SH}$} & {\color{black} K} & {\color{black} $13.5$} & {\color{black} $10.3$} \\
			\midrule
			{\color{black} $N$} & {\color{black} Hz} & {\color{black} $38$} & {\color{black} $38$} \\
			{\color{black} $A_v$} & {\color{black} \%} & {\color{black} $44.36$} & {\color{black} $53.36$} \\
			\midrule
	\end{tabular}}
\end{table}

\pagebreak

\begin{table} [h]
	\centering
	\caption{$COP$ achieved for initial points IP\textsubscript{1} and IP\textsubscript{2} and the different PNMPC tuning }
	\label{tab_COP}
	\scalebox{0.7}[0.7]{
		\begin{tabular}{c c c c} 
			\midrule
			\textbf{Initial point} & \textcolor{black}{$\bm Q_1$} & \textcolor{black}{$\bm Q_2$} & \textcolor{black}{$\bm Q_3$} \\
			\midrule
			\textbf{IP\textsubscript{1}} & \textcolor{black}{$1.1977$} & \textcolor{black}{$1.1950$} & \textcolor{black}{$1.1963$} \\
			\midrule
			\textbf{IP\textsubscript{2}} & \textcolor{black}{$1.2273$} & \textcolor{black}{$1.2268$} & \textcolor{black}{$1.2269$} \\
			\midrule
	\end{tabular}}
\end{table}

\pagebreak

\section*{List of figure captions}

\noindent Figure 1: Mechanical compression system \\

\noindent Figure 2: Classical optimization and control structure \\

\noindent Figure 3: Open-loop simulation results when simple input profiles are applied to both the original SMB model and the reduced-order one \\

\noindent Figure 4: Pressure-specific enthalpy chart of a canonical mechanical compression cycle \\

\noindent Figure 5: Integration of the static model of all cycle components into the optimization procedure \\

\noindent Figure 6: Optimal cycle features for a range of cooling load \\

\noindent Figure 7: Feedback-plus-feedforward controller by Jain \cite{jain2013thermodynamics} \\

\noindent Figure 8: Control performance when trying to achieve optimal operation from initial point IP\textsubscript{1} using PNMPC and FB+FF control \\

\noindent Figure 9: Control performance when trying to achieve optimal operation from initial point IP\textsubscript{2} using PNMPC and FB+FF control  \\

\noindent Figure 10: Control performance when trying to achieve optimal operation from initial point IP\textsubscript{1} and applying different PNMPC tuning \\

\noindent Figure 11: Control performance when trying to achieve optimal operation from initial point IP\textsubscript{2} and applying different PNMPC tuning \\

\noindent Figure B.1: Controllable subspace for constant values $d_1$ and $d_2$ \\

\noindent Figure B.2: Graphic phase portrait of the dynamic system described in Eq. \eqref{eq_modelo_reducido_dinamico} \\

\pagebreak

\pagebreak

\pagebreak

\pagebreak

\pagebreak

\pagebreak

\pagebreak

\pagebreak

\pagebreak

\pagebreak

\pagebreak

\pagebreak

\renewcommand\thefigure{B.\arabic{figure}}
\setcounter{figure}{0}

\begin{figure}[h] 
	\centering
	\centerline{\includegraphics[width=8cm]{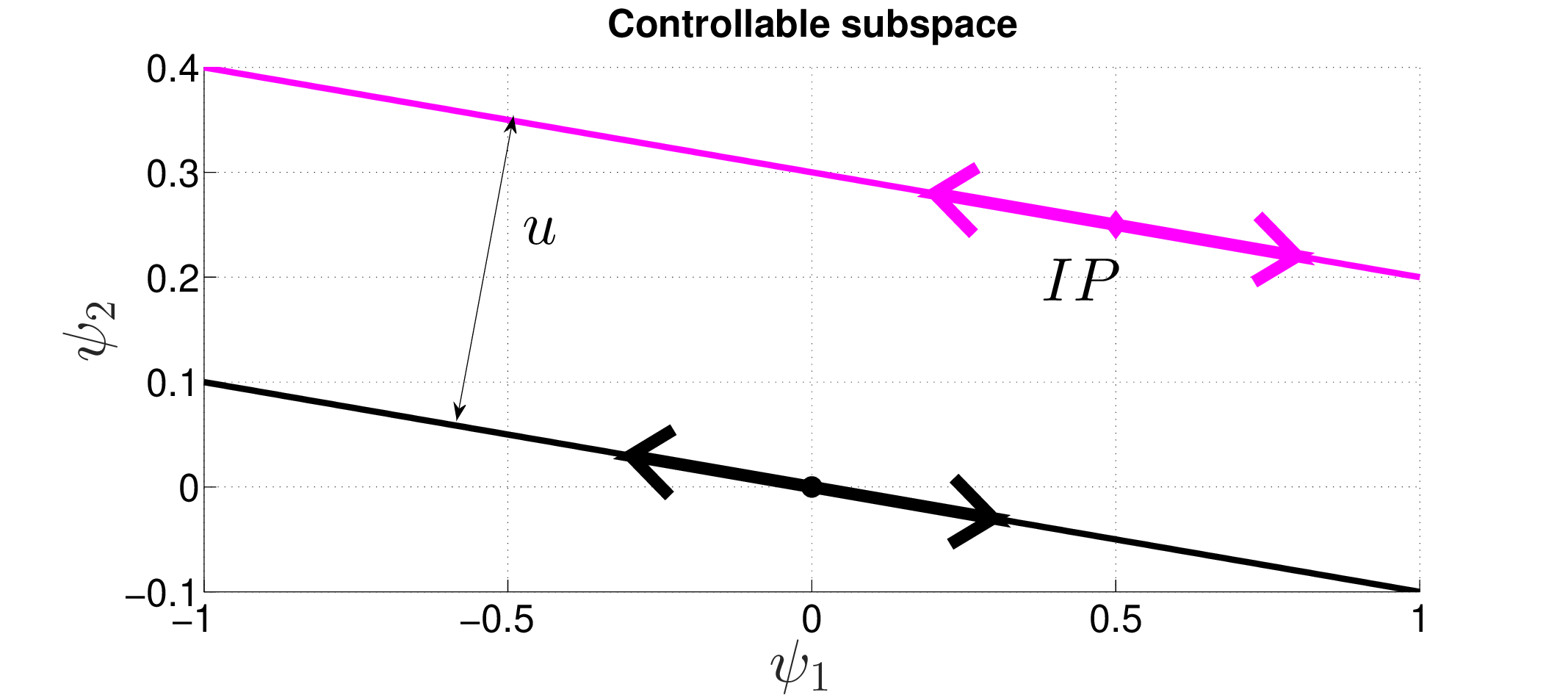}}
	\caption{Controllable subspace for constant values $d_1$ and $d_2$} 
	\label{fig_recta_movimiento}
\end{figure}

\pagebreak

\begin{figure}[h] 
	\centering
	\centerline{\includegraphics[width=8cm]{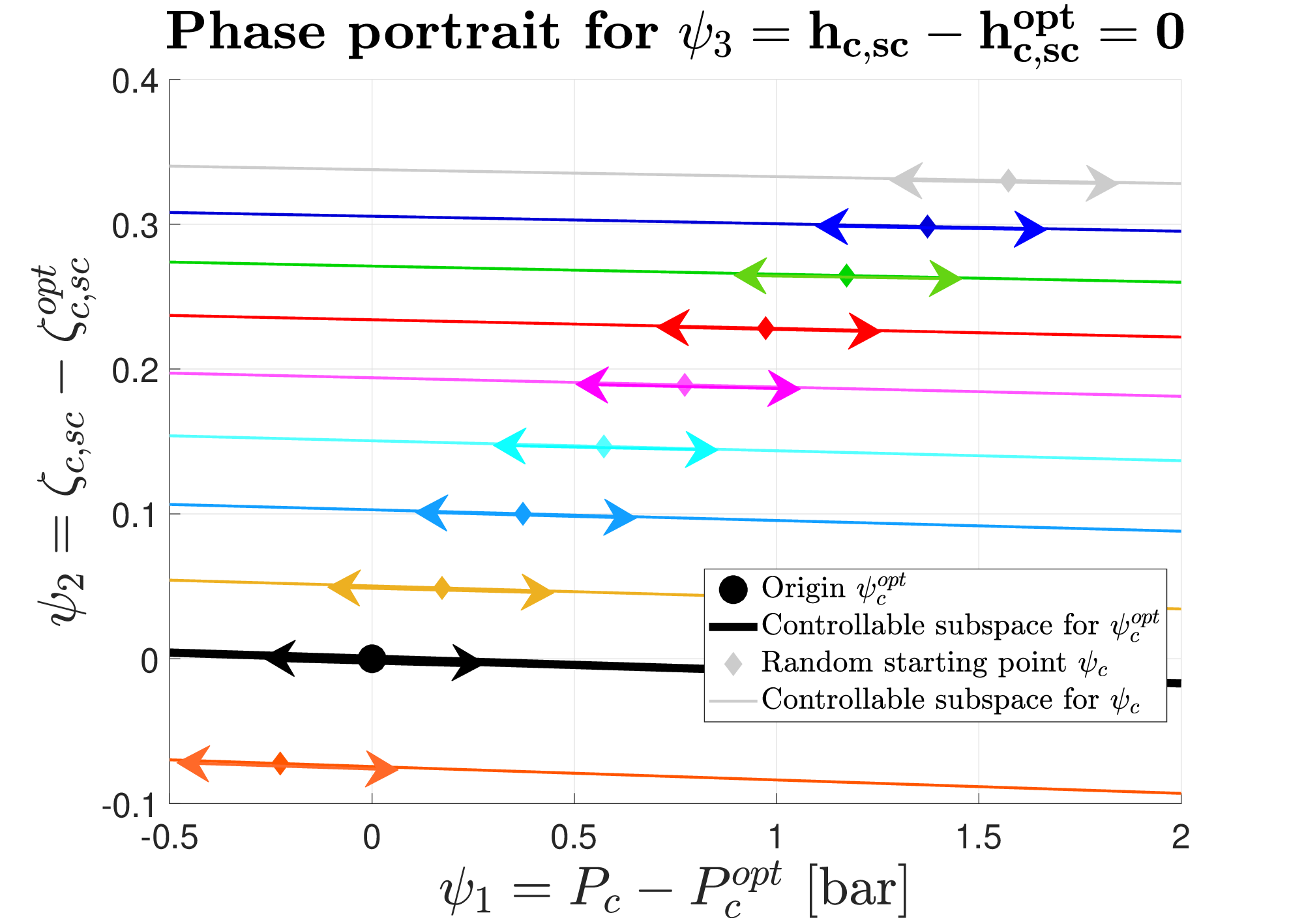}}
	\caption{Graphic phase portrait of the dynamic system described in Eq. \eqref{eq_modelo_reducido_dinamico}} 
	\label{fig_phase_plane}
\end{figure}

\pagebreak


\begin{appendices}	
	
	\setcounter{equation}{0}
	\renewcommand\theequation{A.\arabic{equation}}
	
	\section{Study of controllability based on linear theory} \label{Appendix_LinearControllability}
	
	The reduced-order state-space model of the condenser is shown in Eq. \eqref{eq_modelo_condensador}, whereas the coefficient matrix $\bm Z_{c}$ and the force function $\bm f_{c}$ are presented in Eq. \eqref{eq_Zm1_fm1}. Moreover, the elements $f_{c,i} \; \forall i = 1,2,3$ and $z_{c,ij} \; \forall i,j = 1,2,3$ are gathered in Eq. \eqref{eq_vector_fm1} and \eqref{eq_matrix_Zcm1}, for operating \emph{mode} 1, that turns out to be the most common one \cite{BINLI}. 
	
	\begin{equation}
		\begin{aligned}
			\bm Z_{c}(\bm x_{c},\bm w_{c}) \; \bm{\dot{x}}_{c} & = \bm f_{c} (\bm x_{c},\bm w_{c})\\
			\bm x_{c} & = [P_c \;\; \zeta_{c,sc} \;\;  h_{c,sc}]^T\\
			\bm w_c & = [\dot{m}_{c,sec} \;\; T_{c,sec,in} \;\; \dot{m} \;\; 
			h_{c,in}]^T\\
		\end{aligned}
		\label{eq_modelo_condensador}
	\end{equation}
	
	\begin{equation}
		\begin{aligned}
			\bm Z_{c}(\bm x_{c},\bm w_{c}) & =  \left[	
			\begin{array}{ccc}
				z_{c,11} & -1 & z_{c,13} \\ 
				z_{c,21} &  0 & z_{c,23} \\
				z_{c,31} & -1 & z_{c,33} \\
			\end{array} \right] \\
			\bm f_{c}(\bm x_{c},\bm w_{c}) & = [f_{c,1} \;\; f_{c,2} \;\;
			f_{c,3}]^T\\
		\end{aligned}
		\label{eq_Zm1_fm1}
	\end{equation}
	
	\begin{equation}
		\begin{aligned}
			f_{c,1}(\bm x_{c},\bm w_{c}) & = \dfrac{1}{\rho_{c,tp}\;V_R} \; \dfrac{1}{h_{c,f} - h_{c,sc}} \left(
			\dot{Q}_{c,sc} + \dot{m} (h_{c,f} - h_{c,sc})
			\right) \\
			f_{c,2}(\bm x_{c},\bm w_{c}) & = \dfrac{1}{\rho_{c,tp} \;\zeta_{c,tp} \; V_R} \left(
			\dot{Q}_{c,tp} + \dot{m} (h_{c,g} - h_{c,tp}) +  \dfrac{h_{c,f} - h_{c,tp}}{h_{c,f} - h_{c,sc}} \;\dot{Q}_{c,sc}
			\right) \\
			f_{c,3}(\bm x_{c},\bm w_{c}) & = \dfrac{1}{\rho_{c,sc}\;V_R} \dfrac{1}{h_{c,f} - h_{c,sc}} \left(
			\dot{Q}_{c,sc} + \dot{m} (h_{c,f} - h_{c,sc})
			\right) \\
		\end{aligned}
		\label{eq_vector_fm1}
	\end{equation}
	
	\begin{equation}
		\begin{aligned}
			z_{c,11} (\bm x_{c},\bm w_{c}) & = 
			\dfrac{1}{2}\dfrac{\partial h_{c,g}}{\partial P_c} \left( 
			\dfrac{\zeta_{c,sh}}{\rho_{c,sh}}\dfrac{\partial \rho_{c,sh}}{\partial h_{c,sh}} - \dfrac{\zeta_{c,sh}}{h_{c,g} - h_{c,sh}}
			\right) + \\
			& + \dfrac{\zeta_{c,sh}}{\rho_{c,sh}}\dfrac{\partial \rho_{c,sh}}{\partial P_c} + 
			\dfrac{\zeta_{c,sh}}{\rho_{c,sh}(h_{c,g} - h_{c,sh})} + 
			\dfrac{\zeta_{c,tp}}{\rho_{c,tp}}\dfrac{\partial \rho_{c,tp}}{\partial P_c} - \\
			& - \dfrac{\zeta_{c,sh}}{\rho_{c,tp}(h_{c,g} - h_{c,sh})} - \dfrac{\zeta_{c,sc}}{\rho_{c,tp}(h_{c,f} - h_{c,sc})} + \\
			& + \dfrac{1}{2}\dfrac{\partial h_{c,g}}{\partial P_c} \left( 
			\dfrac{\rho_{c,sh}}{\rho_{c,tp}}\dfrac{\zeta_{c,sh}}{h_{c,g} - h_{c,sh}} +
			\dfrac{\partial \bar{\gamma}_c}{\partial P_c} \dfrac{\zeta_{c,tp}}{\rho_{c,tp}}\dfrac{\partial\rho_{c,tp}}{\partial \bar{\gamma}_c}
			\right) \\
			z_{c,13} (\bm x_{c},\bm w_{c}) & = \dfrac{\rho_{c,sc}}{\rho_{c,tp}}\dfrac{\zeta_{c,sc}}{h_{c,f} - h_{c,sc}} \\ 
			z_{c,21} (\bm x_{c},\bm w_{c}) & = \dfrac{\partial h_{c,tp}}{\partial P_c} - \dfrac{1}{\rho_{c,tp}} + \dfrac{\partial \bar{\gamma}_c}{\partial P_c} \dfrac{\partial h_{c,tp}}{\partial \bar{\gamma}_c} - \dfrac{\zeta_{c,sh}}{\zeta_{c,tp}}\dfrac{1}{\rho_{c,tp}}\dfrac{h_{c,g} - h_{c,tp}}{h_{c,g} - h_{c,sh}} - \\ 
			& -\dfrac{\zeta_{c,sc}}{\zeta_{c,tp}}\dfrac{1}{\rho_{c,tp}}\dfrac{h_{c,f} - h_{c,tp}}{h_{c,f} - h_{c,sc}} + \dfrac{1}{2}\dfrac{\partial h_{c,g}}{\partial P_c} 	\dfrac{\zeta_{c,sh}}{\zeta_{c,tp}} \dfrac{\rho_{c,sh}}{\rho_{c,tp}} \dfrac{h_{c,g} - h_{c,tp}}{h_{c,g} - h_{c,sh}} \\
			z_{c,23} (\bm x_{c},\bm w_{c}) & = \dfrac{\zeta_{c,sc}}{\zeta_{c,tp}} \dfrac{\rho_{c,sc}}{\rho_{c,tp}} \dfrac{h_{c,f} - h_{c,tp}}{h_{c,f} - h_{c,sc}} \\
			z_{c,31} (\bm x_{c},\bm w_{c}) & = - \dfrac{\zeta_{c,sc}}{\rho_{c,sc}}\dfrac{1}{h_{c,f} - h_{c,sc}} \\ 
			z_{c,33} (\bm x_{c},\bm w_{c}) & = - \dfrac{\zeta_{c,sc}}{\rho_{c,sc}} \dfrac{\partial \rho_{c,sc}}{\partial h_{c,sc}} + \dfrac{\zeta_{c,sc}}{h_{c,f} - h_{c,sc}} \\
		\end{aligned}
		\label{eq_matrix_Zcm1}
	\end{equation}
	
	By using matrix manipulation, a new model representation may be obtained, which is characterised by a zero element in vector $\bm{\hat{f}}_{c}$, as shown in Eq. \eqref{eq_modelo_subactuado}.
	
	\begin{equation}
		\begin{aligned}
			\bm{\hat{Z}}_{c} \; \bm{\dot{x}}_{c} & = \bm{\hat{f}}_{c} \\
			\bm{\hat{Z}}_{c} & =  \left[	
			\begin{array}{ccc}
				z_{c,11} & -1 & z_{c,13} \\ 
				z_{c,21} & 0 & z_{c,23} \\
				z_{c,31} - \dfrac{\rho_{c,tp}}{\rho_{c,sc}} z_{c,11} 
				& -1 + \dfrac{\rho_{c,tp}}{\rho_{c,sc}} 
				& z_{c,33} - \dfrac{\rho_{c,tp}}{\rho_{c,sc}} z_{c,13} \\
			\end{array} 
			\right] \\
			\bm{\hat{f}}_{c} & = [f_{c,1} \;\;f_{c,2} \;\; 0]^T \\
		\end{aligned}
		\label{eq_modelo_subactuado}
	\end{equation}
	
	The input variables to the condenser are $\dot{m}$ and $h_{c,in}$, according to Eq. \eqref{eq_modelo_condensador}. These variables are equivalent to $A_v$ and $N$, respectively, since the steady-state models of such components could be used to define the values of $A_v$ and $N$ leading to some desired values of $\dot{m}$ and $h_{c,in}$. Similarly, the variables $f_{c,1}$ and $f_{c,2}$ in Eq. \eqref{eq_modelo_subactuado} could be also considered as virtual manipulable inputs, since, for desired values of the latter, the values of $\dot{m}$ and $h_{c,in}$ could be directly computed. 
	
	The condenser model turns out to be underactuated, since only two manipulated variables are available and up to three states must be regulated, thus an underactuated control problem at the condenser arises. The structure of the matrix $\bm{\hat{Z}}_{c}$ makes it impossible to be singular, considering the typical values of the thermodynamic properties of the refrigerant at the pressure range under study. Therefore it can be inverted and a new model representation arises, as indicated in Eq. \eqref{eq_modelo_subactuado_invertido}, that adopts a linear structure, although all elements in $\bm{\hat{Z}}_{c,inv}$ and $\bm{\hat{f}}_{c}$ depend on the state, the actual manipulable inputs, and the disturbances in a non-linear manner.
	
	\begin{subequations}
		\begin{equation}
			\begin{aligned}
				\bm{\hat{Z}}_{c} \; \bm{\dot{x}}_{c} & = \bm{\hat{f}}_{c} \\
				\bm{\hat{Z}}_{c}^{-1} & \equiv \bm{\hat{Z}}_{c,inv} \\
			\end{aligned}
		\end{equation}	
		\begin{equation}
			\begin{aligned}
				\bm{\dot{x}}_{c} = \bm{\hat{Z}}_{c,inv} \; \bm{\hat{f}}_{c} &= 
				\left[	
				\begin{array}{cc}
					\hat{z}_{c,inv,11} & \hat{z}_{c,inv,12} \\ 
					\hat{z}_{c,inv,21} & \hat{z}_{c,inv,22} \\
					\hat{z}_{c,inv,31} & \hat{z}_{c,inv,32} \\
				\end{array} 
				\right] 
				\left[
				\begin{array}{c}
					f_{c,1} \\ 
					f_{c,2} \\
				\end{array} 
				\right]
			\end{aligned}
			\label{eq_modelo_lineal}
		\end{equation}
		\label{eq_modelo_subactuado_invertido}
	\end{subequations} 
	
	Although the system is highly non-linear, the controllability matrix $\bm{\mathbb{C}}$ could provide some insight, since the structure of the model shown in Eq. \eqref{eq_modelo_subactuado_invertido} is essentially linear. According to Eq. \eqref{eq_controlabilidad}, the degree of controllability of the system seems to be two. However, the conclusions made from this study should be considered as only indicative, since for instance matrix $\bm B_c$ in Eq. \eqref{eq_modelo_subactuado_invertido} is highly non-linear and this has not been taken into consideration when computing the controllability matrix. Anyway, this preliminary study suggests that there might be controllability problems when trying to move the full condenser state vector from an arbitrary point to the sought one.
	
	\begin{subequations}
		\begin{equation}
			\bm{\dot{x}}_{c} = \bm{A}_{c}\;\bm x_{c} + \bm{B}_{c}\;\bm u_{c}
		\end{equation}	
		\begin{equation}
			\bm{A}_{c} = 
			\left[	
			\begin{array}{ccc}
				0 & 0 & 0 \\ 
					0 & 0 & 0\\
					0 & 0 & 0 \\
			\end{array} 
			\right] \;\;\;\;
			\bm{B}_{c} = \bm{\hat{Z}}_{c,inv} =  
			\left[	
			\begin{array}{cc}
				\hat{z}_{c,inv,11} & \hat{z}_{c,inv,12} \\ 
				\hat{z}_{c,inv,21} & \hat{z}_{c,inv,22} \\
				\hat{z}_{c,inv,31} & \hat{z}_{c,inv,32} \\
			\end{array} 
			\right] \;\;\;\;
			\bm u_{c} = \bm{\hat{f}}_{c} = 
			\left[
			\begin{array}{c}
				f_{c,1} \\ 
				f_{c,2} \\
			\end{array} 
			\right]
		\end{equation}
		\begin{equation}
			\bm{\mathbb{C}} = 
			\left[	
			\bm{B}_{c} \;\;\;\; \bm{A}_{c}\;\bm{B}_{c} \;\;\;\; \bm{A}_{c}^{2}\;\bm{B}_{c}
			\right] \;\;\;\;\;\;\;\;\;\;\;\;
			rank(\bm{\mathbb{C})} = rank(\bm{B}_{c}) \leqslant 2
		\end{equation}
		\label{eq_controlabilidad}
	\end{subequations} 
	
	\setcounter{equation}{0}
	\renewcommand\theequation{B.\arabic{equation}}
	
	\section{Non-linear controllability study} \label{Appendix_NonLinearControllability}
	
	\begin{statement}
		The controllable subspace of a canonical mechanical compression refrigeration system is two-dimensional, in the case that only the expansion valve opening and the compressor speed/power act as the control actions, that involves one degree of underactuation.
	\end{statement}
	
		
		This conclusion is made applying the phase portrait approach. The choice of this graphic method stems from the complexity of the dynamic equations of the condenser model detailed in Appendix \ref{Appendix_LinearControllability}, that include a number of refrigerant-specific functions to evaluate several thermodynamic properties. It renders an analytical non-linear study intractable.
		
		Given the condenser model shown in Eq. \eqref{eq_modelo_lineal}, the deviations of all state variables with respect to their optimal values are gathered in vector $\bm{\psi}_c$, as indicated in Eq. \eqref{eq_modelo_error}, where it is assumed that the optimal state $\bm{x}_c^{opt}$ does not change to get the dynamic model of $\bm{\psi}_c$.
		
		\begin{equation}
			\begin{aligned}
				\bm{\psi}_c &\equiv \bm{x}_c - \bm{x}_c^{opt} = \left[	
				\begin{array}{c}
					P_{c} - P_{c}^{opt} \\ 
					\zeta_{c,sc} - \zeta_{c,sc}^{opt} \\
					h_{c,sc} - h_{c,sc}^{opt} \\
				\end{array}
				\right] \\
				\bm{\dot{\psi}}_{c} &= \bm{\dot{x}}_{c} - \bm{\dot{x}}_c^{opt} = 
				\bm{\dot{x}}_{c} = 
				\left[	
				\begin{array}{cc}
					\hat{z}_{c,inv,11} & \hat{z}_{c,inv,12} \\ 
					\hat{z}_{c,inv,21} & \hat{z}_{c,inv,22} \\
					\hat{z}_{c,inv,31} & \hat{z}_{c,inv,32} \\
				\end{array} 
				\right] 
				\left[
				\begin{array}{c}
					f_{c,1} \\ 
					f_{c,2} \\
				\end{array} 
				\right]
			\end{aligned}
			\label{eq_modelo_error}
		\end{equation}
		
		To avoid excessive notation complexity, the elements of the dynamic matrix of $\bm{\psi}_c$ are renamed as indicated in Eq. \eqref{eq_redefine}.
		
		\begin{equation}
			\begin{aligned}
				\bm{\dot{\psi}}_{c} = \left[	
				\begin{array}{cc}
					\dot{\psi}_1 \\ 
					\dot{\psi}_2 \\
					\dot{\psi}_3 \\
				\end{array} 
				\right] =
				\left[	
				\begin{array}{cc}
					\hat{z}_{c,inv,11} & \hat{z}_{c,inv,12} \\ 
					\hat{z}_{c,inv,21} & \hat{z}_{c,inv,22} \\
					\hat{z}_{c,inv,31} & \hat{z}_{c,inv,32} \\
				\end{array} 
				\right] 
				\left[
				\begin{array}{c}
					f_{c,1} \\ 
					f_{c,2} \\
				\end{array} 
				\right] \equiv \left[	
				\begin{array}{cc}
					n_{11} & n_{12} \\ 
					n_{21} & n_{22} \\
					n_{31} & n_{32} \\
				\end{array} 
				\right] 
				\left[
				\begin{array}{c}
					f_{c,1} \\ 
					f_{c,2} \\
				\end{array} 
				\right] 
			\end{aligned}
			\label{eq_redefine}
		\end{equation}
		
		The problem dimension is again reduced to increase the comprehensibility of the graphic results of the non-linear controllability study. Then, among the two manipulated inputs regarding the condenser ($\dot{m}$ and $h_{c,in}$), it is assumed that one of them is devoted to ensure that a state variable (for instance, $h_{c,sc}$) matches its optimal value for a given cooling demand. This allows to represent the phase portrait only in two dimensions, corresponding to the remaining state variables ($P_c$ and $\zeta_{c,sc}$).
		
		If one of the manipulated inputs ensures that $h_{c,sc}$ matches its optimal value $h_{c,sc}^{opt}$, then $\dot{\psi}_3$ = 0 J kg\textsuperscript{-1} s\textsuperscript{-1}. It allows to reduce the dimension of the problem as shown in Eq. \eqref{eq_modelo_reducido}.
		
		\begin{subequations}
			\begin{equation}
				\dot{\psi}_3 = 0 \Rightarrow n_{31} f_{c,1} +  n_{32} f_{c,2} = 0  \Rightarrow f_{c,1} = - \dfrac{n_{32}}{n_{31}} f_{c,2} 
				\label{eq_f1}
			\end{equation}
			\begin{equation}
				\left[	
				\begin{array}{c}
					\dot{\psi}_1 \\ 
					\dot{\psi}_2 \\
				\end{array} 
				\right] = \left[	
				\begin{array}{cc}
					n_{11} & n_{12} \\
					n_{21} & n_{22} \\
				\end{array} 
				\right] 
				\left[
				\begin{array}{c}
					f_{c,1} \\ 
					f_{c,2} \\
				\end{array} 
				\right] = \left[	
				\begin{array}{c}
					n_{12} - \dfrac{n_{32}}{n_{31}} n_{11} \\ 
					n_{22} - \dfrac{n_{32}}{n_{31}} n_{21} \\
				\end{array} 
				\right] f_{c,2} \equiv
				\left[	
				\begin{array}{c}
					d_{1} \\ 
					d_{2} \\
				\end{array} 
				\right] f_{c,2}
				\label{eq_modelo_reducido_dinamico}
			\end{equation}
			\label{eq_modelo_reducido}
		\end{subequations}
		
		Therefore, leading the remaining states $\psi_1$ and $\psi_2$ to zero by using only $f_{c,2}$ is actually the control problem, since $f_{c,1}$ must be computed through Eq. \eqref{eq_f1} to ensure that $h_{c,sc}$ matches its optimal value $h_{c,sc}^{opt}$ at any instant. Notice that $f_{c,1}$ and $f_{c,2}$ are the so-called virtual manipulated variables; in order to compute the actual input variables to the condenser model $\dot{m}$ and $h_{c,in}$, the non-linear system of equations generated by the two first ones of Eq. \eqref{eq_vector_fm1} should be solved, once the condenser reduced state vector $\bm x_c$ and the disturbances related to the secondary fluid are known.
		
		It is easy to check that the reduced control problem is not completely controllable if the elements $d_1$ and $d_2$ are constant and state-independent. The controllable subspace is the line defined by Eq. \eqref{eq_determinante} and represented for random values of $d_1$ and $d_2$ in Fig. \ref{fig_recta_movimiento}. Notice that the arrows represent that the controller would be able to move the system in both directions, but in all cases the actuators physically limit the magnitude of the movement.
		
		\begin{equation}
			\left|
			\begin{array}{cc}
				\psi_1 & d_1 \\
				\psi_2 & d_2 \\
			\end{array} 
			\right| = 0  \Rightarrow
			\psi_1 d_2  -  \psi_2 d_1 = 0 \Rightarrow
			\psi_2 = \dfrac{d_2}{d_1} \psi_1
			\label{eq_determinante}
		\end{equation}
		
		That means that the desired point is achievable only if the initial point is within the controllable subspace. Any control law will be unsuccessful in reaching the desired point when starting at any random state out of this subspace. Moreover, the distance $u$ to the controllable subspace holds for any control action $f_{c,2}$, as indicated in Eq. \eqref{eq_distancia}.
		
		\begin{equation}
			\begin{aligned}
				\left[	
				\begin{array}{c}
					\dot{\psi}_1 \\ 
					\dot{\psi}_2 \\
				\end{array} 
				\right] = \left[	
				\begin{array}{c}
					d_{1} \\ 
					d_{2} \\
				\end{array} 
				\right] f_{c,2} 
				\Rightarrow
				u \equiv \left|
				\begin{array}{cc}
					\psi_1 & d_1 \\
					\psi_2 & d_2 \\
				\end{array} 
				\right| \\
				\dot{u} = \dot{\psi}_1 d_2 - \dot{\psi}_2 d_1 = (d_1 d_2 - d_2 d_1) f_{c,2} = 0 \;\;\; \forall f_{c,2}  
			\end{aligned}
			\label{eq_distancia}
		\end{equation}
		
		Obviously, the control problem in the case under study described in Eq. \eqref{eq_modelo_reducido_dinamico} is not so simple, since not only do the elements $d_1$, $d_2$, and $f_{c,2}$ depend on the state variables, but they also depend on the actual control actions $\dot{m}$ and $h_{c,in}$. Then, since an analytical study turns out to be impracticable, a pointwise non-linear controllability study is proposed. 
		
		Firstly, diverse starting points are defined on the phase portrait: all of them are equilibrium points of the system and they all satisfy that $h_{c,sc} = h_{c,sc}^{opt}$, since the phase portrait is represented for the reduced control problem defined in Eq. \eqref{eq_modelo_reducido}, where $\psi_3$ = 0 J kg \textsuperscript{-1} at any instant. 
		
		Secondly, the non-linear equations defining the condenser model are used to compute, for each starting point, the specific values of $d_1$, $d_2$, $f_{c,1}$, and $f_{c,2}$, and therefore the slope of the line which constitutes its controllable subspace, when different feasible control inputs $\dot{m}$ and $h_{c,in}$ are applied. The pair \{${\dot{m},h_{c,in}}$\} must satisfy the hypothesis of $\psi_3$ = 0 J kg \textsuperscript{-1}. Each pair \{${\dot{m},h_{c,in}}$\} is expected to generate different values of $d_1$ and $d_2$, which causes the slope of the controllable subspace to slightly vary. Nevertheless, it can been checked that, for achievable control inputs, such slopes are so similar that the standard deviation with respect to the average value turns out to be smaller than 0.5\% of the average slope.

		Eventually, the controllable subspaces (lines with the average slope previously described) for all starting points are represented on the phase portrait shown in Fig. \ref{fig_phase_plane} for reasonable ranges of $\psi_1$ and $\psi_2$ around the origin. Every thin coloured line represents the one-dimension controllable subspace when starting at the same-coloured initial point, that is marked by a diamond. In order not to include all colours corresponding to several starting points in the plot legend, they have been clustered in the grey entry. The origin has been graphically emphasized, as well as its controllable subspace, in black. Notice that the arrows represent that, if starting at an arbitrary initial point, the movement is possible in both directions along the line representing its controllable subspace, but the magnitude of the movement is constrained by the limits on the actual control actions.   
		
		It can be observed on the phase portrait that the slopes of the lines which constitute the controllable subspaces are very similar, which makes these lines quasi-parallel. As a result, it is not possible in practice to achieve the origin if the initial condition is out of its controllable subspace. 
		
	
	\begin{remark}
		Similar conclusions can be derived if the phase portraits for $\psi_1$ = 0 bar and $\psi_2$ = 0 are depicted, which are not included for the sake of brevity. Furthermore, if the problem dimension was not reduced, a three-dimensional phase portrait would be obtained where the controllable subspaces would turn out to be quasi-parallel planes. Therefore, the degree of controllability turns out to be in practice two.
	\end{remark}
	
\end{appendices}

\end{document}